\documentclass[11pt,a4paper,reqno]{amsart}
\usepackage{amsfonts,amsthm, amscd, epsfig, amsmath, amssymb,enumerate}
\numberwithin{equation}{section}
\makeatletter
\def\captionfont@{\footnotesize}
\def\captionheadfont@{\scshape}

\long\def\@makecaption#1#2{%
  \vspace{2mm}
  \setbox\@tempboxa\vbox{\color@setgroup
    \advance\hsize-6pc\noindent
    \captionfont@\captionheadfont@#1\@xp\@ifnotempty\@xp
        {\@cdr#2\@nil}{.\captionfont@\upshape\enspace#2}%
    \unskip\kern-6pc\par
    \global\setbox\@ne\lastbox\color@endgroup}%
  \ifhbox\@ne 
    \setbox\@ne\hbox{\unhbox\@ne\unskip\unskip\unpenalty\unkern}%
  \fi
  \ifdim\wd\@tempboxa=\z@ 
    \setbox\@ne\hbox to\columnwidth{\hss\kern-6pc\box\@ne\hss}%
  \else 
    \setbox\@ne\vbox{\unvbox\@tempboxa\parskip\z@skip
        \noindent\unhbox\@ne\advance\hsize-6pc\par}%
\fi
  \ifnum\@tempcnta<64 
    \addvspace\abovecaptionskip
    \moveright 3pc\box\@ne
  \else 
    \moveright 3pc\box\@ne
    \nobreak
    \vskip\belowcaptionskip
  \fi
\relax } \makeatother
\def\writefig#1 #2 #3 {\rlap{\kern #1 truecm
\raise #2 truecm \hbox{#3}}}

\newtheorem{theorem}{Theorem}[section]
\newtheorem{lemma}[theorem]{Lemma}
\newtheorem{prop}[theorem]{Proposition}
\newtheorem{cor}[theorem]{Corollary}

\newtheorem{definition}[theorem]{Definition}
\newtheorem{remark}[theorem]{Remark}

\newcommand{\RR}{{\mathbb R}}

\newcommand{\ZZ}{{\mathbb Z}}
\newcommand{\NN}{{\mathbb N}}
\newcommand{\PP}{{\mathbb P}}

\newcommand{\EE}{{\mathbb E}}
\newcommand{\GG}{{\mathbb G}}
\newcommand{\FF}{{\mathbb F}}

\newcommand{\id}{{\mathbb I}}
\newcommand{\ep}{{\epsilon}}
\newcommand{\toromi}{{\mathbb T}^{d}_{\epsilon}} 
\newcommand{\toroma}{{\mathbb T}^{d}} 
\newcommand{\Av}{{\rm Av}}  
\newcommand{\Var}{{\rm Var}} 

\newcommand{\lstor}{\mathcal{L}}

\newcommand{\dstor}{\mathcal{D}}

\newcommand{\gstor}{\mathcal{G}}
\newcommand{\mstor}{\mathcal{M}}

\newcommand{\no}{\noindent}
\newcommand{\la}{ \lambda}
\newcommand{\La}{ \Lambda}

\newcommand{\De}{ \Delta}

\newcommand{\dto}{\downarrow}
\newcommand{\uto}{\uparrow}

\newcommand{\cA}{\ensuremath{\mathcal A}}
\newcommand{\cB}{\ensuremath{\mathcal B}}

\newcommand{\cD}{\ensuremath{\mathcal D}}
\newcommand{\cE}{\ensuremath{\mathcal E}}
\newcommand{\cF}{\ensuremath{\mathcal F}}
\newcommand{\cG}{\ensuremath{\mathcal G}}
\newcommand{\cH}{\ensuremath{\mathcal H}}

\newcommand{\cL}{\ensuremath{\mathcal L}}
\newcommand{\cM}{\ensuremath{\mathcal M}}
\newcommand{\cN}{\ensuremath{\mathcal N}}

\newcommand{\frL}{\ensuremath{\mathfrak L}}

\newcommand{\bbE}{{\ensuremath{\mathbb E}} }

\newcommand{\bbN}{{\ensuremath{\mathbb N}} }

\newcommand{\bbP}{{\ensuremath{\mathbb P}} }

\newcommand{\bbR}{{\ensuremath{\mathbb R}} }

\newcommand{\bbZ}{{\ensuremath{\mathbb Z}} }

\def\bF{{\Bbb F}}

\newcommand{\eps}{\epsilon} 

\let\a=\alpha \let\b=\beta  \let\c=\chi \let\d=\delta  \let\e=\varepsilon
 \let\g=\gamma \let\h=\eta    \let\k=\kappa  \let\l=\lambda
      \let\o=\omega      
  \let\s=\sigma \let\t=\tau   \let\th=\vartheta
  
\let\D=\Delta   \let\G=\Gamma  \let\L=\Lambda \let\The=\Theta
\let\O=\Omega      

\let\Z=\integer

\let\neper=e
\let\ii=i
\let\mmin=\wedge

\def\ie{\hbox{\it i.e.\ }}
\let\id=\identity

\let\sset=\subset

\def\nep#1{ \neper^{#1}}

\def\ov#1{{1\over#1}}

\def\tc{\thsp | \thsp}
\let\<=\langle
\let\>=\rangle

\def\gap{\mathop{\rm gap}\nolimits}

\def\ninf#1{ \| #1 \|_\infty }

\outer\def\nproclaim#1 [#2]#3. #4\par{\medbreak \noindent
   \talato(#2){\bf #1 \Thm[#2]#3.\enspace }%
   {\sl #4\par }\ifdim \lastskip <\medskipamount
   \removelastskip \penalty 55\medskip \fi}
\def\thmm[#1]{#1}
\def\teo[#1]{#1}
\def\sttilde#1{%
\dimen2=\fontdimen5\textfont0
\setbox0=\hbox{$\mathchar"7E$}
\setbox1=\hbox{$\scriptstyle #1$}
\dimen0=\wd0
\dimen1=\wd1
\advance\dimen1 by -\dimen0
\divide\dimen1 by 2
\vbox{\offinterlineskip%
   \moveright\dimen1 \box0 \kern - \dimen2\box1}
}

\def\smallno{\smallskip\noindent}
\def\medno{\medskip\noindent}

\def\\{\hfill\break}

\def\thsp{\thinspace}

\def\tthsp{\kern .083333 em}

\def\?{\mskip -10mu}

\begin{document}

\title[Hydrodynamic limit of a disordered lattice gas]
{Hydrodynamic limit of a disordered lattice gas }
\date{\today}
\author[A. Faggionato]{Alessandra Faggionato} \address{Fakult\"at II -
Mathematik und Naturwissenschaften, Technische
  Universit\"at Berlin, Strasse des 17. Juni 136, 10623 Berlin, Germany}
\email{faggiona\@@math.tu-berlin.de}
\author[F. Martinelli]{Fabio Martinelli} \address{Dipartimento di Matematica,
Universita' di Roma Tre, L.go S. Murialdo 1, 00146 Roma, Italy}
\email{martin\@@mat.uniroma3.it}

\vskip 1cm
\begin{abstract}
We consider a model of lattice gas dynamics in $\bbZ^d$ in the presence of
disorder. If the particle interaction is only mutual exclusion and if
the disorder field is given by i.i.d. bounded random variables, we prove
the almost sure existence of the hydrodynamical limit in dimension $d\ge
3$.
The limit equation is a non linear diffusion equation with diffusion
matrix characterized by a variational principle.
\vskip1cm
\noindent
{\em 2000 MSC: 60K40, 60K35, 60J27, 82B10, 82B20 }

\noindent
{\bf Key words and phrases}: hydrodynamic limit, disordered systems,
lattice gas dynamics, exclusion process.

\end{abstract}
\maketitle

\section{Introduction}
\noindent
Hopping motion of particles between spatially distinct locations is one
of the fundamental transport mechanisms in solids and it has been
extensively used in a variety of models, including electron conduction
in disordered systems under a tight binding approximation. The
interested
reader is refereed to \cite{bello} for a detailed physical review.  \\
From a mathematical point of view, hopping motion is often modeled as an
interacting particle system in which each particle performs a random
walk over the sites of an ordered lattice like $\bbZ^d$, with jump rates
depending, in general, on the interaction with the nearby particles and,
possibly, on some external field. Typically the interaction between the
particles is assumed to be short range with an hard core exclusion rule
(multiple occupancy of any site is forbidden) and only jumps between
nearest neighbors sites are allowed. In the conduction models the hard
core exclusion condition reflects the underlying Pauli exclusion
principle for electrons. The main focus of the mathematical and physics
literature on hopping motion models has been the understanding of
transport properties and particularly of the collective diffusive
behavior (see for instance \cite{Spohn}).
\\
In this paper we consider an interacting particle system related to
conduction of free electrons in doped crystals that can be described as
follows. A particle sitting on a site $x$ of the cubic lattice $\bbZ^d$
waits an exponential time and then attempts to jump to a neighbor site
$y$.  If the site $y$ is occupied then the jump is canceled otherwise it
is realized with a rate $c_{xy}^\a$ depending only on the values
$(\a_x,\,\a_y)$ of some external quenched disorder field
$\{\a_x\}_{x\in\bbZ^d}$ that, for simplicity, is assumed to be a
collection of i.i.d. bounded random variables.  Our assumptions on the
transition rates are quite general. We require them to be translation
covariant, strictly bounded and positive (to avoid trapping phenomena),
and to satisfy the detailed balance condition w.r.t. to the (product)
Gibbs measure $\mu^\a(\h) \propto \nep{-H^\a(\h)}$,
$H_\a(\h)=-\sum_{x}\a_x\h_x$, where $\h_x$ is the particle occupation
number at site $x$. These requirements are general enough to include
some popular models like the Random Trap and the Miller--Abrahams
models, but not other models like the Random Barrier Model in which the jumps rates
between $x,y$ is assumed to depend only on the unoriented bond $[x,y]$
\cite{KW1}. For a detailed derivation of the Hamiltonian $H^\a$ in the
tight-binding approximation and a
discussion of the regime of its validity we refer to \cite{bello}. \\

Since in the linear--response regime the conductivity in a solid is
linked to the diffusion matrix via the Einstein relation (see
\cite{Spohn}), our main target has been the study of the bulk diffusion
of the disordered lattice gas discussed above. Our main result states that, for $d\ge
3$, for almost any realization of the random
field $\a$, the diffusively rescaled system has hydrodynamical limit
given by a non linear differential equation
$$
\partial_t m =\nabla\cdot( D(m)\nabla m)
$$ where $m(t,\theta)$ denotes the macroscopic density function at time
$t$ at the point $\theta$ of the $d$--dimensional torus in $\bbR^d$ with
unit volume and the non random matrix $D(\cdot)$ is the diffusion
matrix. Moreover, we give a variational characterization of the matrix
$D(m)$ in terms of the distribution of the random field $\a$ similar to
the usual Green--Kubo formula and we prove that $\inf_m D(m) >0$ and
that $D(\cdot)$ is continuous in the open interval $(0,1)$.\\ We remark
that the above result without the restriction on the dimension $d$, was
already announced in \cite{Q} several years ago together with some
sketchy ideas for its proof.  However the details of the proof have
never since been published and some of the technical estimates suggested
in \cite{Q} turned out to be troublesome even in the absence of disorder
(symmetric simple exclusion model) as explained in \cite{AF}, chapter 6.
Therefore we decided to tackle again the problem but we were forced to
take a different route w.r.t. that indicated in \cite{Q}.  \\We also
observe that the problem of collective behavior in disordered lattice
gas has been discussed mathematically in other papers, but, to the best
of our knowledge, only for models with either homogeneous equilibrium
measures (see for example \cite{N}, \cite{F} for the one--dimensional
Random Barrier model and its Brownian version) or with periodicity in
the random field $\a$ allowing to solve directly the generalized Fick's
law (see \cite{R} and \cite{W} for the one--dimensional Random Trap
model having random field $\a$ of period 2) or finally for models
satisfying the so called ``gradient condition'' (see below) \cite{KK}.
From the physical point of view, diffusion of lattice gases in systems
with site disorder has been studied mainly by means of simulations and
more or less rough approximations like mean field . We refer the
interested reader to \cite{KP}, \cite{KW1}, \cite{KW2}, \cite{KPW},
\cite{K} and to \cite{GP} for an iterative procedure to compute
corrections to the mean--field approximation.

\\
The main technical features of the model considered here are the absence
of translation invariance (for a given disorder configuration) and the
non validity of the so called gradient condition. This condition
corresponds to the Fick's law of fluid mechanics according to which the
current can be written as the gradient of some function. Since the
continuity equation states that $\partial_t m=\nabla\cdot J$, $J$ being
the macroscopic current, the main problem is to derive $J$ from the
family of microscopic instantaneous currents
$j^\a_{x,y}(\h):=c_{x,y}^\a(\h) \bigl (\h_x-\h_y\bigr)$, defined as the
difference between the rate at which a particle jumps from $x$ to $y$
and the rate at which a particle jumps from $y$ to $x$. The gradient
condition (the Fick's law) is satisfied if, for each disorder
configuration $\a$, there exists a local function $h^\a(\h)$ such that
$j^\a_{x,x+e}(\eta) =\t_{x+e} h^\a(\h)- \t_x h^\a(\h)$ for any
$x\in\bbZ^d$, where $\t_x h^\a(\h) := h^{\t_x\a}(\t_x\h)$ and $\t_x
\h,\,\t_x\a$
denote the particle and disorder configurations $\h,\,\a$ 
translated by the vector $x$. \\
If the system satisfies the gradient condition, the derivation of $J$ is
not too difficult (see \cite{KL} and reference therein). It is however
simple to check (as in \cite{Spohn}, p. 182) that our system never
satisfies the gradient condition except for constant disorder field
$\a$. We thus have to appeal to the methods developed by Varadhan
\cite{V}, Quastel \cite{QQ} and Varadhan-Yau \cite{VY} (see also
\cite{KL} and references therein) for studying the hydrodynamic limit of
non disordered non gradient systems. There the main idea is to prove a
generalized Fick's law of the form
\begin{equation}\label{sonntag}
j^\a_{0,e}\approx \sum_{e'\in\cE}D_{e,e'}(m_\ell)(\eta_{e'} -\eta_0)+\cL^\a g
\end{equation}
for a suitable non random matrix $D(m)$, where $m_\ell$ is the particle
density in a cube centered in the origin of mesoscopic side $\ell$,
$g(\a,\h)$ is a
local function and $\cL$ is the generator of the dynamics.\\
One (among many others) main difficulty in proving such an approximation
for a disordered system is due to the fact that the disorder itself
induces strong fluctuations in the gradient density field as it is
easily seen by taking, for any fixed disorder configuration $\a$, the
average w.r.t.  to the Gibbs measure $\mu^\a$ of (\ref{sonntag}). By
construction the current $j_{0,e}^\a$ and the fluctuation term $\cL g$
have in fact zero average while the average of $\eta_{e'} -\eta_0$ (we
neglect the factor $D(m_\ell)$ for simplicity) is in general $O(1)$
because of the disorder. However, and this is a key input, the average
over the disorder of the Gibbs average of $\mu^\a\bigl(\eta_{e'}
-\eta_0\bigr)$ vanishes and therefore one can hope to tame the disorder
induced fluctuations in the gradient of the density field by first
smearing them out using suitable spatial averages and then by appealing
to the ergodic properties of the disorder field $\a$, at least in high
enough dimension. It turns
out that the above sketchy plan works as soon as $d\ge 3$.\\

\smallno We conclude this short introduction with a plan of the paper.
In section 2 we fix the notation, describe the model and state the main
results. In section 3 and section 4 we discuss most of the ``high
level'' technical tools (entropy estimates, perturbation theory,
spectral gap bounds) and complete the proof of the main theorems
following the standard route of non gradient systems, modulo some key
technical results. In section 5 we discuss in detail the problem of the
fluctuations of the gradient density field induced by the disorder.
Section 6 is devoted to the proof of several technical bounds while in
section 7 we discuss at length central limit variance, closed and exact
forms in our context together with our own interpretation of the long
jump method described in \cite{Q}. Finally some very technical
estimates are collected in an appendix at the end. \\
We finish by saying that most of the material presented here is based on
the unpublished thesis \cite{AF} written by one of us (A.F) where an
expanded version of several of the arguments used in this paper can be found.


\subsection*{Acknowledgments}
Part of this work was done while both authors were visiting the
Institute H. Poincar\'e during the special semester on ``Hydrodynamic
limits''. We would like to thank the organizers F. Golse and S. Olla for
their kind invitation and the stimulating scientific atmosphere
there. We are also grateful to J. Quastel for providing unpublished
notes on the problem and for sharing his insight of the subject. We are
also grateful to S. Olla, C. Landim, G.B. Giacomin for many
enlightening discussions and to P. Caputo for his proof of the spectral
gap bound.

\section{Notation, the model and main results}
In this section we fix the notation, we define the model and state our
main result.

\subsection{Notation}
\label{notation}

\\
\emph{Geometric setting.} We consider the $d$ dimensional lattice $\Z^d$
with sites $x = \{x_1, \ldots, x_d \}$, canonical basis $\cE$ and norm
$|x|=\max\{|x_1|,\dots, |x_d|\}$. The bonds of $\ZZ^d$ are non oriented
couple of adjacent
sites and a generic bond will be denoted by $b$. \\
The cardinality of a finite subset $\L\sset \bbZ^d$ is denoted by $|\L|$
and $\bF$ denotes the set of all nonempty finite subsets of $\Z^d$. \\
Given $\ell\in\NN$ we denote by $\L_{\ell}$ the cube centered at the
origin of side $2\ell +1$. If $\ell=2j+1$ we also set $Q_{\ell} =
\L_{j}$. The same cubes centered at $x$ will be denoted by $\L_{x,\ell}$
and $Q_{x,\ell}$ respectively. More generally, for any $V\sset \bbZ^d$
and $x\in \bbZ^d$, we will set $V_x := V+x$.
\\
Next, given $e\in \cE$ and $\ell=2\ell'+1$ with $\ell'\in \NN$, we let
\begin{equation}
  \label{blocchetto}
\L_\ell^{1,e}:=\L_{-(\ell'+1)e,\ell'}\,,\quad \L_{\ell}^{2,e}:=\L_{\ell'
  e,\ell'}\,,
\quad \L_{\ell}^e := \L_\ell^{1,e} \cup \L_\ell^{2,e}.
\end{equation}
Finally, given $\ep\in (0,1)$ such that $\ep^{-1}\in \bbN$, we
define the discrete torus of spacing $\ep$ by $\toromi :=\ZZ^d
/\ep^{-1} \ZZ^{d}$.  The usual $d$--dimensional torus
$\RR^d/\ZZ^d$ (with unite volume)  will instead be denoted by $\toroma$.
$\cM_1(\toroma)$ will denote the set of positive Borel measures on
$\toroma$ with total mass bounded by $1$, endowed of the weak
topology, while $\cM_2\sset \cM_1$ will denote the set of measures
in $\cM_1$ which are absolutely continuous w.r.t. the Lebesgue
measure with density $\rho$ satisfying $\ninf{\rho}\leq 1$.

\\
\emph{Spatial averages.} We will make heavy use of spatial averages and
it is better to fix from the beginning some handy notation. Given
$\L\in \bF$ and $\ell\in \bbN$, the spatial
average of $\{f_x\}_{x\in \bbZ^d}$ in $\L\cap \ell\,\bbZ$ will be denoted by
$\Av^{(\ell)}_{x\in \L}f_x$. When $\ell=1$ we will simply write
$\Av_{x\in \L}\,f_x$.

\smallno Next, given $e\in \cE$ and two odd integers $\ell=2\ell'+1$,
$s=2s'+1$ such that $\frac{s}{\ell}\in \bbN$, we let $Q^{(\ell)}_{s}:=
\ell\,\bbZ^d\cap Q_s$. Notice that, if we divide the cube $\L_s^{1,e}$
in cubes of side $\ell$, the centers of these cubes form the set
$Q^{(\ell)}_{x,s}$ with $x= -(s'+1)e$.
\\
With these notation we define the $(\ell,s,e)$ spatial average around
$y\in \bbZ^d$ by
\begin{equation}
\label{dinoce}
\Av^{\ell,s}_{z,y}\,f_z := \frac{1}{(s/\ell)}\sum_{i=0}^{(s/\ell)-1}
\Av_{x\in Q^{(\ell)}_{s}}\,f_{y+x+(\ell' + i\ell -s') e}.
\end{equation}
The motivation of introducing such  a spatial average will be discussed
in subsection \ref{urca}.

\\
\emph{The disorder field.}
We assume the disorder to be described by a collection of real i.i.d random
variables $\a:=\{\a_x\}_{x\in \bbZ^d}$ such that $\sup_x |\a_x|\le B$
for some finite constant $B$. The corresponding product measure on
$\O_{\rm D}:=[-B,B]^{\bbZ^d}$
will be denoted by $\bbP$. Expectation w.r.t. $\bbP$ will be denoted by $\bbE$.

\smallno
Notice that, for any given $\ep\in (0,1)$ such that $\ep^{-1}$ is an odd
integer, the random field $\a$ induces in a natural way a random field
on $\toromi$ via the identification of $\toromi$ with the cube
$Q_{1/\ep}$. For notation convenience the induced random field will
always be denoted by $\a$.\\
Finally, given $\a\in \O_{\rm D}$ and $\L \sset \Z^d$, we define
$\a_\L:=\{\a_x\}_{x\in \L}$.

\\
\emph{The particle configuration space.} Our particle
configuration space is $\O = S^{\Z^d}$, $S=\{0,1\}$ endowed with
the discrete topology, or $\O_\L = S^\L$ for some $\L\in \bF$.
When $\L= \toromi$ we will simply write $\O_\ep$.  Given $\h\in
\O$ and $\L \sset \Z^d$ we denote by $\h_\L$ the natural
projection over $\O_\L$.  Given two sites $x,y \in \ZZ^d$ and a
particle configuration $\h$ we denote by $\h^{x,y}$ and $\h^x$ the
configurations obtained from $\h$ by exchanging the values of $\h$
at $x,y$ and by ``flipping'' the value of $\h$ at $x$
respectively. More precisely,
\begin{equation*}
  \bigl(\h^{x,y}\bigr) _z: =
\begin{cases}
\h _y & \text{if} \; z=x \\
\h _x & \text{if} \; z=y \\
\h _z & \text{otherwise}
\end{cases}, \quad \quad (\h ^x)_z :=
\begin{cases}
1-\h_x & \text{if} \; z=x \\
\h _z & \text{otherwise}.
\end{cases}
\end{equation*}
Sometimes we will write $\h^{x,y}:= S_{x,y}\h$ and call $S_{x,y}$ the
\emph{exchange operator} between $x$ and $y$.
Finally, given a probability measure $\mu$
on $\O_\L$, we will denote by $\Var_\mu(\xi)$ the variance of the random
variable $\xi$ w.r.t. $\mu$, by $\mu(\xi;\xi')$ its covariance  with the random
variable $\xi'$ and by $\mu(\xi,\xi')$ the scalar product between
$\xi$ and $\xi'$ in the Hilbert space $L^2(\O_\L,d\mu)$.

\\
\emph{Local functions.}  If $f$ is a measurable function on $\tilde\O:=
\O_{\rm D}\times \O$, the support of $f$, denoted by $\D_f$, is the
smallest subset of $\Z^d$ such that $f(\a,\h)$ depends only on
$\a_{\D_f},\s_{\D_f}$ and $f$ is called {\it local} if $\D_f$ is finite.
By $\|f\|_\infty$ we mean the supremum norm
of $f$. Given  two sites $x,y\in\ZZ^d$ we define
\begin{align*}
\nabla_{x,y}f(\a,\h)&:=f(\a,\h^{x,y})-f(\a,\h),\\
 \nabla_x f(\a,\h)&:=f(\a,\h^x)-f(\a,\h).
\end{align*}
We write $\GG$ for the set of measurable, local and bounded
functions $g$ on $\tilde\O$ and for any $g\in \GG$ we introduce the formal series $\underline
g$
\begin{equation*}
\underline g:=\sum_{x\in\ZZ^d}\t_x g
\end{equation*}
where $\t _x f(\a,\h):=f(\t_x\a,\t_x\h)$ and
$\t_x\a$ and $\t_x\h$ are the disorder and particle
configurations translated by $x\in\bbZ^d$ respectively:
$$
(\t_x\a)_z:=\a_{x+z}, \quad (\t_x\h)_z:=\h_{x+z}.
$$
Although the above series is only formal, by the locality of $g$, the gradient
$\nabla_{x,y}\,\underline g$ is meaningful for any $x,y\in\ZZ^d$.

\\
\emph{Limits.} Given $n$ parameters $\ell_1,\dots\ell_n$ we use the compact notation
$\lim_{\ell_n\rightarrow\ell_n ',\dots,\,\ell_1\rightarrow\ell_1'}$ for
the ordered limits $\lim_{\ell_n\rightarrow\ell
  _n'}\dots \lim_{\ell_1\rightarrow\ell_1'}$.
The same convention is valid when ``$\lim$'' is
replaced by ``$\limsup$'' or ``$\liminf$''.

\subsection{The model}
\label{model}
In this subsection we describe the lattice gas model at the microscopic scale
$\ep$ for a given disorder configuration $\a$.

\\
\emph{Gibbs measures.}
Given an external chemical potential
$\l\in \bbR$, the Hamiltonian of the system in the set $\L\sset \bbZ^d$ is defined  as
$$
H^\a_\L(\h)=-\sum_{x\in\L}(\a_x+ \la)\h_x
$$
and the corresponding grand canonical Gibbs measure on $\O_\L$, denoted by
$\mu^{\a,\l}_\L$, is simply the product measure
\begin{equation}
  \label{grand}
  \mu_\L^{\a,\l}(\h) := \frac{1}{Z_\L^{\a,\l}}\exp(-H^\a_\L(\h))
\end{equation}
where $Z_\L^{\a,\l}$ is such that  $\mu_\L^{\a,\l}(\O_\L)=1$.\\
 For our purposes it is important to introduce also the
canonical measures $\nu_{\L,m}^\a$. Let $N_\L(\h) = \sum_{x\in
\L}\h_x$ and let $m\in [0,\frac{1}{|\L|},\dots ,1]$. Then
\begin{equation}
  \label{canonical}
  \nu_{\L,m}^\a(\cdot) = \mu_\L^{\a,\l}(\cdot \tc N_\L=m|\L|)
\end{equation}
The random variable $N_\L$ will usually be referred to as the number of
particles and \hbox{$m_\L := N_\L/|\L|$} as the particle density or
simply the density.  The set of all canonical measure $\nu_{\L,m}^\a$ as
$m$ varies in $[0,\frac{1}{|\L|},\dots ,1]$ will be denoted by
$\cM^\a(\L)$. Notice that $\nu_{\L,m}^\a$ does not depend on the
chemical potential $\l$.  However, as it is well known \cite{CM1}, the
canonical and grand canonical Gibbs measures are closely related if the
chemical potential $\l$ is canonically conjugate to the density $m$ in
the sense that the \emph{average density} w.r.t. $\mu_\L^{\a,\l}$ is
equal to $m$.  With this in mind, for any $m\in [0,1]$, we define the
\emph{empirical chemical potential} $\l_{\L}(\a,m)$ as the unique value
of $\l$ such that $\mu^{\a,\l}_\L(N_\L)=m|\L|$, the \emph{annealed
  chemical potential} $\l_0(m)$ as the unique $\l$ such that
$ \EE\bigl[\mu^{\a,\l}(\h_0)\bigr]=m$ and the corresponding static
compressibility $\chi(m)$ as
$\c(m)=\EE\bigl[\mu^{\a,\l_0(m)}(\h_0;\h_0)\bigr]$.  Since
$\frac{\partial}{ \partial \l}\mu^{\a,\l}_\L(f)=\mu
^{\a,\l}_\L(f;N_\L)$ for any local function $f$, we get the
following thermodynamic relations:
$$\frac{\partial}{\partial m}\l_{\L}(\a,m) = \bigl [\,\mu_\L^{\a,
\l_\L (\a,m)} (m_\L;N_\L) \bigr]^{-1}\quad\text{ and } \quad
\frac{\partial}{\partial m} \l_0(m)= \chi(m)^{-1}.
$$

\medno
{\bf Notation warning.}
\emph{ From now on, in order to keep the notation to an
acceptable level, we need to adopt the following shortcuts whenever no confusion
arises.
\begin{enumerate}[i)]
\item Most of the times the label $\a$ will be omitted. That
  means that quantities like $\mu_\L^{\l}(f)$ will actually be random
  variables w.r.t the disorder $\a$. Moreover, the label $\l$ of the chemical potential will be omitted   when $\l=0$.
\item If the region $\L$ on which the Gibbs measures or, later,
  the generator of the dynamics are defined coincides with $\toromi$, then the suffix
  $\L$ will be simply replaced by $\ep$ while if
$\L=\bbZ^d$ it will simply be dropped (i.e. $\mu_\ep:=\mu ^\a_{\toromi}$).
\item The symbol $\mu_\L^{\l(m)}$ will always denote the grand
  canonical Gibbs measure on $\O_\L$ with empirical chemical potential
  $\l_\L(\a,m)$.
\item The letter $c$ will denote a generic positive constant depending
only on $d$ and $B$ that may vary from estimate to estimate.
\end{enumerate}
}

\\
\emph{The dynamics.} The lattice gas dynamics we are interested in is the  continuous time
Markov chain on $\O_\ep$ described by the Markov generator
$\ep^{-2}\cL_\ep$ where $\cL_\ep:= \cL_{\toromi}$ and for any $\L\sset
\bbZ^d$
\begin{equation*}
\cL_\L f(\h)= \sum_{b\sset \L} \cL_b f(\h)
\end{equation*}
where, for any bond $b=\{x,y\}$,
\begin{equation*}
  \cL_{x,y}f(\h) := c_{x,y}^\a(\h)\nabla_{x,y}f(\h)
\end{equation*}
The non-negative real quantities $c^\a_{x,y}(\h)$ are the
\emph{transition rates} for the process. They are defined as
$$
c^\a_{x,x+e}(\h)=f_e(\a_x,\h_x,\a_{x+e},\h_{x+e}) \quad \forall x\in\ZZ^d,
\;e\in\cE
$$
where  $f_e$ is a generic  bounded function on
$\bigl([-B,B]\times S\bigr)^2$ such that
$f_e(a,s,a',s')=f_e(a',s',a,s)$ and $f_e\geq c>0$ for a suitable constant $c$. Thanks to this definition the
transition rates are \emph{translation covariant}, i.e.
\begin{equation*}
 c_{x+z,y+z}^\a(\h)=c_{x,y}^{\tau_z\a}(\t_z\h) \quad \quad \forall z\in\ZZ^d.
\end{equation*}
The key hypothesis on the transition rates is the \emph{detailed balance
  condition} w.r.t the Gibbs measures $\mu^\l_\L$, $\L\sset \bbZ^d$ and
$\l\in \bbR$, \ie
\begin{equation*}
f_e(a,s,a',s')=f_e(a,s',a',s)e^{-(s'-s)(a'-a)}\qquad \forall
e\in\cE,\;a,a'\in [-B,B], \,s,s'\in S
\end{equation*}
which implies that the generator $\lstor_\L$ becomes a selfadjoint
operator on $L^2(\mu_\L^\l)$ for any $\l$. Actually, since the moves of
the Markov chain generated by $\cL_\L$ do not change the number of
particles, for any canonical
Gibbs measure $\nu\in \cM(\L)$ the operator $\cL_\L$ is selfadjoint on $L^2(\nu)$
with a positive spectral gap
\begin{equation}
  \label{eq:def_gap}
 \gap(\cL_\L,\nu) := \inf \Bigl\{ \, \frac{\nu(f,-\cL_\L f)}{\Var_\nu(f)};\,\Var_\nu(f)\neq
0\,\Bigr\}
\end{equation}
and the corresponding Markov chain is irreducible on $\{\h\in \O_\L:\;
N_\L(\h)=n\}$ for any $n\in [0,1,\dots, |\L|]$.

\smallno Given $g\in\GG$ we denote by $\cL g$ the function $\sum
_{b\sset\ZZ^d}\cL_b\, g$.  Given $\D\sset \L$ and a probability measure
$\mu$ on $\O_\L$, for any $f$ with support inside $\L$ we will set
$$
\cD_\D(f;\mu):=\ov2\sum_{b\sset\D} \mu\bigl(c_b
(\nabla_bf)^2\bigr).
$$
Notice that, if $\D=\L$ and $\mu$ is
either a grand canonical or a canonical measure on $\L$, then the above
expression is nothing but the Dirichlet form of the Markov chain
generated by $\cL_\L$ w.r.t. $\mu$.

\smallno Finally, given a probability measure $\mu$ on $\O_\ep$
and $T>0$, we denote by $\PP_{T}^{\,\a,\mu}$ the distribution at
time $T$ of the Markov chain on $\toromi$ with generator
$\ep^{-2}\cL^\a _\ep$ and initial distribution $\mu$, and by
$\PP^{\,\a,\mu}$ the induced probability measure on the Skorohod
space $D([0,T],\O_\ep)$ (see  \cite{Billingsley}).
 The expectation w.r.t. $\PP^{\,\a,\mu}$ will be denoted by
$\EE^{\a,\mu}$.  Notice that, in turn, $\PP^{\,\a,\mu}$ induces a
probability measure $Q^{\,\a,\mu}$ on $D([0,T],\cM_1)$ by the
formula $ \PP^{\,\a,\mu}\circ \pi_\ep^{-1}$, where
\begin{equation*}
\pi_\ep(\h):=\Av_{x\in\toromi}\,\h_x\,\d_{\ep x}\in \cM _1(\toroma)
\end{equation*}
denotes the empirical measure.\\

\smallno
{\bf Warning} In all the above measures, the crucial dependence on the parameter $\ep>0$
does not appear in the various symbols in order to keep the notation to
an acceptable level.

\subsection{Main results}
Our first result concerns the existence and regularity of the diffusion
matrix $D(m)$ corresponding to the usual Green-Kubo matrix (see
 \cite{Spohn}, proposition $2.2$ page $180$).
\begin{theorem}
\label{matrix}
Let $d\geq 3$. Then for any density $m\in(0,1)$ there exists a unique symmetric
$d\times d$ matrix $D(m)$, such that
\begin{equation}
\bigl(a,D(m)a\bigr) =\frac{1}{2\chi(m)}
  \inf_{g \in \GG}\,\sum_{e \in \cE}\,
 \EE\Bigl[\,\mu^{\a,\la_0(m)}\Bigl(
 c_{0,e}^\a\bigl(a_e (\h_e - \h_0)
+ \nabla_{0,e} \underline g\bigr)^2 \Bigr)\,\Bigr] \quad \forall a\in\RR^d.
\end{equation}
Moreover $D(m)$ is continuous in the open interval $(0,1)$ and
$$
0<c^{-1}\id \leq D(m)\leq c\id \quad \forall m \in (0,1)
$$
for some positive
constant $c$.
\end{theorem}
\begin{remark}
  We actually expect the matrix $D$ to be continuously extendable to the
  closed interval $[0,1]$. In particular we expect that $D(m)$ converges
  to the diffusion matrix of the random walk of a single particle in the
  random environment $\a$ as $m$
  goes to zero, as confirmed by simulations (see \cite{KPW}).
\end{remark}

In order to state the next main result we need the following
definition.
\begin{definition}
Given a Lebesgue absolutely continuous measure
$m(\theta)d\theta\in \cM_2(\toroma)$, a sequence of probability measures
$\mu^\ep$ on $\O_\ep$ is said to correspond to the macroscopic profile
$m(\cdot)$ if, under $\mu^\ep$, the random variable $\pi_\ep$ in $\cM_1(\toroma)$
converges in probability to $m(\theta)d\theta$ as $\ep \dto 0$, i.e. for
any smooth function $H$ on $\toroma$ and any $\d>0$
\begin{equation*}
\lim_{\ep\dto 0}\mu^\ep\bigl(\bigl|\,\Av_{x\in\toromi}H(\ep x)\eta_x
-\int_{\toroma}H(\theta)m(\theta)d\theta\,\bigr|>\d\bigr)=0\,.
\end{equation*}
\end{definition}

\smallno With the above definition the existence of the hydrodynamical
limit for almost all disorder configurations reads as follows.
\begin{theorem}\label{HL}
Let $d\geq 3$, let $T>0$ and assume that $D(m)$ can be continuously extended to the closed
interval $[0,1]$. Then almost all disorder configurations $\a$ satisfy
the following property. Let $m_0(\theta)\in\cM_2$ and suppose that the Cauchy problem
\begin{equation}\label{eden}
\begin{cases}
 \partial_t m (t,\theta)=\nabla_\theta\Bigl(D\bigl(m(t,\theta)\bigl)
\nabla_{\theta}\,m(t,\theta)\Bigr)\\
m(0,\theta)=m_0(\theta)
\end{cases}
\end{equation}
has a unique weak solution $m\in C([0,T],\cM_2)$  satisfying the
energy estimate
\begin{equation}\label{kadewe}
\int_0^T ds \int_{\toroma} d\theta\, |\nabla_\theta m (t,\theta) |^2 <\infty.
\end{equation}
Let also $\{\mu^\ep\}_{\ep>0}$ be a sequence of
probability measures on $\O_\ep$ corresponding to the macroscopic
density profile $m_0(\theta)$. Then the measure $Q^{\a,\,\mu^\ep}$
converges weakly to the probability measure on $D([0,T],\cM_1)$
concentrated on the path $\{m(t,\theta)d\theta\}_{t\in[0,T]}$.
In particular, for any $0\le t\le T$, the sequence of time dependent
probability measures $\{\PP^{\a,\mu^\e}_t\}_{\ep>0}$ corresponds to the
macroscopic density profile $m(t,\theta)$, i.e. for any smooth function
$H$ on $\toroma$ and any $\d>0$
\begin{equation}
\label{babbonatale100}
\lim_{\ep\dto 0}\PP^{\a,\mu_\ep}_t\bigl(\bigl|\,\Av_{x\in\toromi}H(\ep x)\eta_x
-\int_{\toroma}H(\theta)m(t,\theta)d\theta\,\bigr|>\d\bigr)=0.
\end{equation}
The thesis remains valid also if $D(m)$ has no continuous extension
provided that one assumes instead that, for some fixed $\rho\in (0,1)$,
there exists a sequence of probability measures $\mu^{\ep}_\ast$ on
$\O_\ep$ such that
\begin{equation}
\label{condition_cut} H[\mu^\ep|\mu^\ep_\ast]=o(\ep^{-d})\qquad
\text{and}\qquad
\inf_{\ep}\inf_{x\in\toromi}\min\bigl(\mu^\ep_\ast(\h_x),1-\mu^\ep_\ast(\h_x)\bigr)\geq\rho,
\end{equation}
where $H[\cdot|\cdot]$ denotes the relative entropy.
\end{theorem}
\begin{remark}
Notice that condition (\ref{condition_cut}) becomes rather natural if the
initial profile $m_0(\cdot)$ satisfies
$\rho\leq m_0(\theta)\leq 1-\rho$ for any  $\theta \in \toroma$.
\end{remark}

\section{Plan of the proof of the two main theorems}
\label{Proof-main}
The proof of theorem \ref{matrix} will be given in section
\ref{biodegradabile}  and it
is based on more or less standard techniques. The proof of theorem
\ref{HL} is more involved and it can be divided into several steps that
we illustrate in what follows. In order to work in the simplest possible
setting, in the sequel we assume that the diffusion matrix $D$ can be
continuously extended to the closed interval $[0,1]$. Only at the end
(see subsection \ref{generale}) we will explain how to treat the other
case.

Let us begin with some remarks on the weak interpretation of (\ref{eden}) and
(\ref{kadewe}). Let $A(m)$, $m \in [0,1]$, be a $d\times d$
matrix such that
$A'(m)=D(m)$ so that
$$
\bigl(D(m(t,\theta))\nabla_\theta
m(t,\theta)\bigr)_e=\sum_{e'\in\cE}\partial_{\theta_{e'}}A_{e,e'}(m(t,\theta))\,,\quad
\forall e\in \cE.
$$
It is simple to prove (see appendix of \cite{AF}) that given $m\in
D([0,T],\cM_2)$ there is a measurable function $m(t,\theta)$ univocally
defined up to sets of zero Lebesgue measure such that
$m_t=m(t,\theta)d\theta$ for any $t\in [0,T]$ (see appendix of
\cite{AF}). In what follows, we will often identify $m$ with the funtion
$m(t,\theta)$. 
\\
A path $m\in D([0,T],\cM_2)$ is called a weak solution
of (\ref{eden}) if $m(0,\cdot)=m_0(\cdot)$ and
\begin{equation*}
\Phi(m,H)=0 \quad \forall H\in C^{1,2}([0,T]\times\toroma)
\end{equation*}
where
\begin{equation}
\label{esmeralda}
\begin{split}
 \Phi(m,H):&=\int_{\toroma} m(T,\theta) H(T,\theta)\,d\theta
-\int_{\toroma} m(0,\theta) H(0,\theta)\,d\theta -\int_0^T
\int_{\toroma}m(s,\theta)\partial_s H(s,\theta)\,d\theta\,ds\\
& - \sum_{e,e'}\int_0^T\int_{\toroma}
A_{e,e'}\bigl(m(s,\theta)\bigr)\,\partial^2_{\theta_e,\theta_{e'}}H(s,\theta)\,d\theta
\,ds.
\end{split}
\end{equation}
Moreover, $m \in D([0,T],\cM_2)$ satisfies the energy estimate
(\ref{kadewe}) if
\begin{equation}
\label{weak_energy_estimate} \sup_{e\in \cE}\sup_{H\in C^1([0,T]\times\toroma)}
\int_0^T\int_{\toroma}\bigl(2\,m\, \partial_{\theta_e} H-
H^2\bigr)d\theta\,ds<\infty.
\end{equation}

\smallno
{\bf Warning}. In what follows, we will introduce some other mesoscopic
scales in addition to the microscopic scale $\ep$. For example,
 we will introduce some positive  scale parameters $a,b$ and  consider
the  mesoscopic scales $\bigl [\frac{a}{\ep}\bigl ]$ and
$\bigl[\frac{b}{\ep}\bigl ]$, where $[\cdot]$ denotes the integer part. For simplicity of notation  these new scales will be denoted only by $\frac{a}{\ep}$ and $\frac{b}{\ep}$. Moreover, we will introduce the scale $n$ where $n$ is a positive  odd integer. The property of  $n$ to be odd will be always understood.

\subsection{Tightness}\label{tightness}
The first step toward the proof of theorem \ref{HL} is to show that, for
 all disorder configurations $\a$, if $\{\mu^\ep\}_{\ep>0}$ is a sequence of probability measures on $\O_\ep$ then  the sequence of measures on
$D([0,T],\cM_1)$, $\{Q^{\a,\,\mu^\ep}\}_{\ep>0}$, is relatively compact.
For this purpose it is enough to use the Garsia-Rodemich-Rumsey
inequality as done in \cite{KL}, chapter $7$, section $6$.

\subsection{Regularity properties of the limit points} \label{gatorade}
In the second step one proves that, for almost all $\a$, given a
sequence $\{\mu^\ep\}_{\ep>0}$ of probability measures on
$\O_\ep$, any limit point $Q^\a$ of the sequence
$\{Q^{\a,\,\mu^\ep}\}_{\ep>0}$ is concentrated on paths enjoying a
certain regularity property. For this purpose we first observe
that, for any $\a$, $Q^\a$ must satisfy
$Q^\a\bigl(\,C([0,T],\cM_2)\,\bigr)=1$, since for any
$\h\in\O_\ep$, $H\in C(\toroma)$ and $b\sset\toromi$
\begin{equation*}
\bigl|\,\pi_\ep(\h)[H]\,\bigr|\leq \Av _{x\in\toromi}|H(\ep x)|
\quad \text{and}\quad \bigl |\, \pi_\ep (\h^b)[H]-\pi_\ep(\h)[H]\,
\bigr|\leq 2\,\|H\|_\infty \ep ^d.
\end{equation*}
Moreover, if the sequence of $\{\mu^\ep\}_{\ep>0}$ corresponds to
the macroscopic profile $m_0(\theta)$, then necessarily
\begin{equation}\label{girotondo100}
Q^\a\bigl (\,m \in C([0,T],\cM_2)\,:\;m (0,\theta)= m_0
(\theta)\,\bigr)=1\quad \forall\a.
\end{equation}
The key result here, whose proof will be given later in section
\ref{gatorade11}, is the following. \\
Given a path $\h(\cdot) \in
D([0,T],\O_\ep)$, $x \in \toromi$ and $\ell \in \bbN$, let
$m_{x,\ell}(t)$ be the particle density of $\h(t)$ in the cube
$Q_{x,\ell}$. Then we have
\begin{lemma}
\label{teti}{\bf (Energy estimate).} Let $d\geq 3$, let $T>0$ and
assume that $D(m)$ can be continuously extended to the closed
interval $[0,1]$. Then almost
any disorder configurations $\a$ have the following property. For
any sequence $\{\mu^\eps\}_{\eps >0}$ of probability measures on
$\O_\ep$
 and any $e\in \cE$
\begin{equation}
\label{jerusalem}
 \sup_{b>0} \limsup_{a\dto 0,\ep\dto
0}\EE^{\a,\mu^\ep}\Bigl (\Av_{x\in\toromi}\int_0^T
\Bigl[\,\frac{m_{x+\frac{b}{\ep}e,\frac{a}{\ep}}(s)-m_{x,\frac{a}{\ep}}(s)}{b}
\,\Bigr]^2 \,ds \Bigr)<+\infty.
\end{equation}
Moreover any limit point $Q^\a$ of the sequence
$\{Q^{\a,\,\mu^\ep}\}_{\ep>0}$ satisfies
\begin{equation}
\label{mirta}
 Q^\a\Bigl\{\,m \in C([0,T],\cM_2)\,:\text{\rm
 l.h.s. of (\ref{weak_energy_estimate} )}  < \infty \Bigr\}=1.
\end{equation}
\end{lemma}

\subsection{Microscopic identification of the hydrodynamic equation}
\label{poka}
In the third step of the proof one identifies at the
microscopic level the hydrodynamic equation. It is convenient to introduce some more notation.
Given $e,e'\in \cE$, two positive
numbers $a,b$ and a
smooth function $H$
on $[0,T ]\times\toroma$, we set
\begin{equation}
\label{eq:F}
\begin{split}
  {\bar H}_{b,a,\ep} &:=\Av_{x\in\toromi}\Bigl[H(T,\ep x)\h_x(T)-
  H(0,\ep x)\h_x(0)-
  \int_{0}^{T}ds\,\h_x(s)\partial_s H(s,\ep x)\Bigr]  \\
  & \\
  &+\sum_{e,e'\in \cE}\int_0^T ds\,\Av_{x\in \toromi} \nabla^\ep_e
  H(s,\ep x) D_{e,e'}\bigl(m_{x,\frac{a}{\ep}}(s)\bigr)
\Bigl[\,
\frac{m_{x+\frac{b}{\ep}e',\frac{a}{\ep}}(s)-m_{x-\frac{b}{\ep}e',\frac{a}{\ep}}(s)}{2b}
\,\Bigr]
\end{split}
\end{equation}
where $\nabla^\ep_e H(s,\ep x) :=\ov \ep \bigl[H(s,\ep x+\ep e )-H(s,\ep x)\bigr]$.

\smallno
The following theorem, whose proof will be discussed in a little while,
corresponds to the microscopic identification of the hydrodynamical
equation.
\begin{theorem}
\label{zeus} Let $d\geq 3$, let $T>0$ and assume that $D(m)$ can
be continuously extended to the closed interval $[0,1]$. Then
almost all disorder configurations $\a$ have the following
property. For any sequence $\{\mu^\eps\}_{\eps >0}$ of probability
measures on $\O_\ep$, any $\d>0$ and any $H\in C^{1,2}([0,T
]\times\toroma)$
\begin{equation}
\label{era}
 \limsup_{b\dto 0,\,a\dto 0,\,\ep\dto 0}
\PP^{\a,\mu^\e}\bigl(\,|{\bar H}_{b,a,\ep}|>\d\,\bigr)=0.
\end{equation}
\end{theorem}

\noindent
The proof of theorem \ref{HL}, given Lemma \ref{teti} and theorem
\ref{zeus}, now follows by more or less standard arguments and it can be
found in section 1.5 of \cite{AF}.
\section{Proof of theorem \ref{zeus} modulo some technical steps}\label{pf_zeus}
In this section we prove theorem \ref{zeus} modulo
certain technical results that will be discussed in the remaining sections.
Following \cite{VY} the first main step is to reduce the proof
of the theorem to the eigenvalue estimates of certain symmetric
operators, via the entropy inequality and the Feynman--Kac formula. To this aim we define
$j_{x,x+e}$ as the instantaneous current through the oriented bond
$\{x,x+e\}$, i.e. as the difference between the rate at which a particle
jumps from $x$ to $x+e$ and the rate at which a particle jumps from $x
+e$ to $x$. It is simple to check that
\begin{equation*}
j_{x,x+e}(\h)=c_{x,x+e}(\h)(\h_x-\h_{x+e})
\end{equation*}
and
\begin{equation*}
  \cL_\ep\h_x=\sum_{e\in\cE}(-j_{x,x+e}(\h)+j_{x-e,x}(\h)).
\end{equation*}
In particular (see lemma $5.1$, appendix $1$ in \cite{KL}, or
\cite{AF}), for any smooth $H(t,x)$, integration by parts and
stochastic calculus show that
\begin{equation}
\label{eq:martingale1}
\begin{split}
& \Av_{x\in\toromi}\bigl[H(T,\ep x)\h_x(T)- H(0,\ep
x)\h_x(0)\bigr]=\\ \Av_{x\in\toromi}
 \int_{0}^{T}\partial_s &H(s,\ep x)\h_x(s)ds+\ep^{-1}\sum_{e\in\cE}
 \Av_{x\in\toromi}\int_0^T \nabla^\ep_e
H(s,\ep x) j_{x,x+e} ds + M(T)
\end{split}
\end{equation}
where $M(\cdot)$ is a martingale w.r.t
$\PP^{\mu^\ep}$ satisfying
\begin{equation}
\label{eq:martingale2}
 \PP^{\mu^\ep}\bigl[\,|M(T)|>\d\,\bigr]\leq
c(H)\,\d^{-2}\ep^d \qquad \forall\d>0.
\end{equation}
In order to benefit of the ergodicity of the system, it is convenient to
replace the current $j_{x,x+e}$ in (\ref{eq:martingale1}) by its
local average around $x$. To this aim let us
introduce a new scale parameter $\ell$, that will be sent to $\infty$ after the
limit $\ep\dto 0$. Then, because of the smoothness of the
function $H$, for any $\ell \gg 1$ one can safely replace in the r.h.s.
of (\ref{eq:martingale1}) the current $j_{x,x+e}$ by a local average
$\Av_{y:|y-x|\leq \ell _1}j_{y,y+e}$, $\ell_1:=\ell-\sqrt{\ell}$, in the
sense that, for any $\d>0$
\begin{equation}
  \lim _{\ep\dto 0}
\PP^{\mu^\ep}\Bigl[\,\bigl|{\ep}^{-1}\Av_{x\in\toromi}\int_{0}
^{T}\nabla_e^\ep H(s,\ep x)\bigl[j_{x,x+e}-\Av_{y:|y-x|\leq \ell _1}j_{y,y+e}\bigr]
\,ds
\,\bigr|>\d\,\Bigr]=0.
\label{local}
\end{equation}
The key observation in the theory of non-gradient systems is that,
thanks again to  stochastic calculus,
\begin{equation}
 \label{eq:martingale3}
\lim _{\ep\dto 0}
\PP^{\mu^\ep}\bigl[\,\bigl|{\ep}^{-1}\Av_{x\in\toromi}\int_{0}
^{T}\nabla_e^\ep H(s,\ep x)\t_x\mathcal{L}g\,ds
\,\bigr|>\d\,\bigr]=0\qquad \forall \d>0,\, \forall g\in\GG
\end{equation}
and similarly for $\Av_{y:|y-x|\leq \ell _1}\t_y\cL g$ in place of
$\t_x\mathcal{L}g$.

\smallno
In conclusion,
thanks  to (\ref{eq:martingale1}), (\ref{eq:martingale2}), (\ref{local}) and
(\ref{eq:martingale3}), in order to prove (\ref{era}) it is enough to show that
for almost all disorder configuration $\a$ and for any $e\in \cE$
\begin{equation}
\label{aiace}
\begin{split}
 \inf_{g\in\GG}\limsup_{b\dto 0,\,a\dto 0,\,\ell\uto\infty,\,\ep\dto0}\EE^{\mu^\ep}&\Bigl(\,\bigl|\int_{0}
^{T}{\ep}^{-1}\Av_{x\in\toromi}\nabla_e^\ep H(s,\ep x)\Bigr[\Av_{y:|y-x|\leq \ell
  _1}(j_{y,y+e}+\t_y\cL g)\\
&+ \sum_{e'\in\cE}D_{e,e'}(m_{x,\frac{a}{\ep}})\bigl[
\frac{m_{x+\frac{b}{\ep}e',\frac{a}{\ep}}
-m_{x-\frac{b}{\ep}e',\frac{a}{\ep}}}{2b/\ep}\bigr] \Bigr]\,ds
\,\bigr|\Bigr)=0.
\end{split}
\end{equation}
We next  reduce (\ref{aiace}) to certain equilibrium eigenvalue
estimates by means of the entropy inequality and the Feynman-Kac
formula (see proposition \ref{fk}). Let us recall the former:
given two probability measures $\pi,\pi'$ on the same probability space, for
any $\b>0$ and any bounded and measurable function $f$,
\begin{equation}
\label{ent_ineq}
\pi(f)\leq\b^{-1}\bigl\{H(\pi\tc\pi')+\ln\bigl(\pi'(e^{\b f})\bigr)\bigr\}
\end{equation}
where $H(\pi\tc \pi')$ denotes the entropy of $\pi$ w.r.t. $\pi'$. It
is simple to verify that, for any  initial distribution
$\mu$ on $\O_\ep$, the relative entropy between the path measure $\PP^\mu$
starting from $\mu$ and the equilibrium path measure $\PP^{\mu_\ep}$ starting from
the grand canonical measure $\mu_\ep$ with zero
chemical potential, satisfies
$$
H\bigl(\PP^{\mu}\tc \PP^{\mu_\ep}\bigr)\leq c\,\ep^{-d}.
$$
Therefore, for any $\g>0$ and any function $h$ on $[0,T]\times\O_\ep$
\begin{equation}
\label{eq:entropy1000}
 \EE^{\mu}\bigl(\,\bigl|
 \int_0^T h(s,\h(s)) ds\bigr|\,\bigr)\leq
\frac{c}{\g}+\frac{\ep^d}{\g}\ln\EE^{\mu_\ep} \Bigl(
\exp\bigl \{\g\ep^{-d}\bigl|\int_0^T h(s,\h(s))
ds\bigr|\bigr\}\Bigr).
\end{equation}
The Feynman--Kac formula (see proposition \ref{fk}) now shows that,
\begin{gather}
\label{eq:fk1000}
 \frac{\ep^d}{\g}\ln \EE^{\mu_\ep}\Bigl( \exp\bigl
\{\g\ep^{-d}\bigl(\pm \int_0^T h(s,\h(s))
ds\bigr)\bigr\}\Bigr)\leq \nonumber \\ \int_0^T \sup spec _{L^2(\mu_\ep)}
\bigl\{\pm h(s,\cdot)+\g^{-1}\ep^{d-2}\cL_\ep\bigr \}\, ds \;.
\end{gather}
We now apply the above reasoning to the function $h(s,\h)= \text{
  integrand of (\ref{aiace})}$. Since for any $\ep >0$
$\sup_{s\in [0,T]}\sup_{x\in \toroma} |\nabla_e^\ep H(s,\ep x)| \le
  c(H)$, after a suitable reparametrization of $\g$,
in order to prove (\ref{aiace}) it is enough to prove the
following key eigenvalue estimate.
\begin{prop}   
\label{wannsee}
Let $d\geq3$. Then, almost all disorder configurations $\a$ have the
following property. For all $\g >0$
\begin{equation}
\label{faust}
 \inf_{g\in\GG} \limsup_{b\dto 0,\,a\dto 0,\,\ell\uto\infty,\,\ep\dto0}
\,\sup_J\,\sup spec _{L^2(\mu_\ep)}\bigl \{\ep^{-1}
\bar J_{b,a,\ell,\ep}^g+\g\ep^{d-2}\cL\bigr\}\leq 0
\end{equation}
where
\begin{equation}
  \label{eq:kalinka}
\begin{split}
&\bar J_{b,a,\ell,\ep}^g:=
\Av_{x\in\toromi} J(\ep x)\Bigr[\Av_{y:|y-x|\leq \ell
  _1}(j_{y,y+e}+\t_y\cL g)\\
&+ \sum_{e'\in\cE}D_{e,e'}(m_{x,\frac{a}{\ep}})\bigl[
\frac{m_{x+\frac{b}{\ep}e',\frac{a}{\ep}}
-m_{x-\frac{b}{\ep}e',\frac{a}{\ep}}}{2b/\ep}\bigr] \Bigr]
\end{split}
\end{equation}
and $J$ varies in $\{J\in C(\toroma)\,:\,\|J\|_\infty\leq 1\}$.
\end{prop}

\subsection{Some technical tools to bound the spectrum}\label{lete}
Before we turn to the proof of proposition \ref{wannsee}, let us
introduce some tools to deal with the
eigenvalue problem posed in (\ref{faust}). \\
We begin by recalling a useful subadditivity property of the supremum of the spectrum
of a selfadjoint operator and explain its role in the so--called
\emph{localization technique}.

\smallno
Given a finite family
$\{X_i\}_{i\in I}$ of self-adjoint operators on $L^2(\mu_\ep)$,
\begin{equation}
\sup spec _{L^2(\mu_\ep)} \{\sum_{i\in I} X_i \}
\leq \sum_{i\in I} \sup spec_{L^2(\mu_\ep)}\{ X_i \},
\label{koukaburra}
\end{equation}
and similarly with $\sum_{i\in I}$ replaced by $\Av_{i\in I}$.
The subadditivity property
allows one to exploits the \emph{localization} method which is best
explained by means of an example, although the underlying idea has a much
wider application.
 Let $\ep>0$, $\ell<\frac{1}{\ep}$ and for any
$x\in\toromi$ let $f_x$ be a local function with support in $\L_{x,\ell}$.
Recall the definition of $\cM(\L_{x,\ell})$ as the set of canonical
Gibbs measures on $\L_{x,\ell}$. Then
\begin{equation}
\label{cappy}
\begin{split}
\sup &spec_{L^2(\mu_\ep)}\bigl\{\Av_{x\in\toromi}f_x+ \ep^{d-2}
\cL_\ep\bigr\} \\ &\leq\Av_{x\in\toromi}\sup spec _{L^2(\mu_\ep)}\bigl\{f_x+
c\,\ep^{-2}\Av_{b\in\L_{x,\ell}}\cL_b\bigr\}\\
&\leq \Av_{x\in\toromi} \sup_{\nu\in\cM(\L_{x,\ell})}\sup spec _{L^2(\nu)}
\bigl\{ f_x+ c\,\ep^{-2}\Av_{b\in\L_{x,\ell}}\cL_b\bigr\}
\end{split}
\end{equation}
where the former inequality follows from $\ep^d\cL_\ep\leq
c\,\Av_{x\in\toromi}\Av_{b\in\L_{x,\ell}}\cL_b$ together with the
subadditivity
property, while the latter can be easily checked.


Next we state a very general result on $\sup spec _{L^2(\nu)}\{\frL+\b V\}$,
where $\frL$ is an ergodic reversible Markov generator on a finite set $E$ with
invariant measure $\mu$, and whose proof is based on perturbation theory
for selfadjoint operators (see e.g. \cite{KL}).

\begin{prop}
 \label{jamme}
Let $\gap(\frL,\mu)$ be the
 spectral gap of $\frL$ in $L^2(\mu)$ and let, for $\b >0$ and $V:E\mapsto \bbR$,
$$
 \l_\b:=\sup spec _{L^2(\mu)}\{\frL+\b V\}.
$$
Assume without loss of generality $\mu(V)=0$. If
\begin{equation*}
2\b\, \gap(\frL,\mu)^{-1}\ninf{V}<1
\end{equation*}
then
\begin{equation*}
 0\leq\l_\b\leq \frac{\b^2}{1-2\b \,\gap(\frL,\mu)^{-1}
 \|V\|_\infty}\,\mu\Bigl(V,(-\frL)^{-1}\,V\Bigl).
\end{equation*}
\end{prop}

\medno The above proposition suggests that in order to prove proposition
\ref{wannsee} we must be able to estimate:
\begin{enumerate}
\item the spectral gap of
the generator $\cL_{\L}$ in a generic box $\L$;
\item the $H_{-1}$ norm appearing above.
\end{enumerate}
We begin with the first one.
\begin{prop}\cite{Caputo}
Let $\L\sset \bbZ^d$ be a parallelepiped with longest side $\ell$. Then, for all
disorder configurations $\a$ and all $\nu \in \cM(\L)$,
\begin{equation}
\label{eq:gap}
\gap(\cL_\L;\nu) \ge c\,\ell^{-2}
\end{equation}
In particular, for all
disorder configurations and all $\nu\in \cM(\L)$, the following Poincar\'e inequality holds
\begin{equation}
  \label{eq:Poincare}
  \Var_{\nu}(f) \le c\,\ell^2 \cD_\L(f;\nu)
\end{equation}
\end{prop}
\begin{remark}
The key aspect of the above result is the uniformity in the disorder
configuration. Its proof is based on some clever technique developed
recently in \cite{Carlen} to deal with the Kac model for the Boltzmann
equation and extended in \cite{Caputo-Martinelli} and \cite{Caputo2} to
other kind of diffusions. For other models of lattice gas dynamics like
the dilute Ising lattice gas in the Griffiths regime the above uniformity will
no longer be available and a more sophisticated analysis is required
(see \cite{AF} for a discussion).
\end{remark}

\smallno
Let us now tackle with the $H_{-1}$ norm. Unfortunately that will
requires the proof of some technical bounds that, on a first reading,
can be just skipped.

Following the theory of non disordered non-gradient systems, we introduce
the space $\cG \sset \GG$ defined as
\begin{equation}
\label{vialattea}
\cG:=\{g\in\GG\,:\,\exists \,\L\in\FF \; \text{ such that, }
\;\forall \a \text{ and } \forall \nu\in \cM^\a(\L)\,, \;\nu(g)=0\, \}.
\end{equation}
Equivalently, $\cG$ can be defined as the set of functions $g\in\GG$
such that there exists $\L\in\FF$ and $h\in\GG$ with $g=\cL_\L h$. Since the
dynamics is reversible w.r.t. Gibbs measures, this second
characterization assures an integration by parts property that will play
an important role in the sequel. More precisely,
if $g=\cL_\L h$, then, for any $\L'$ containing $\L$ and
any $\nu\in \cM(\L')$, $\nu\bigl(g,f)=\nu(h,\cL_\L f)$.
Moreover, if $V$ and $\D$
are such that $\La_{x}\sset V$ for any $x\in\D$, then for
any $A>0$ and $\nu\in\cM(V)$,
\begin{equation}
\label{int_by_parts}
\bigl|\nu\bigl(\sum_{x\in\D}\t_x g,f\bigr)\bigr|
\leq c(g) |\D|^{\frac{1}{2}}\cD_V (f;\nu)
^{\frac{1}{2}} \leq A\, c(g)|\D|+A^{-1}\,c(g)\cD_V(f;\nu)
\end{equation}
where, for some suitable constant $c(\L)$,
\begin{equation}
\label{silenziosanotte}
c(g):=c(\L)\sup_\a\sup_{\nu\in\cM^\a(V)}\bigl(\nu(h^2)\bigr)
^{\ov 2}.
\end{equation}
A first simple consequence of integration by parts (see 
chapter 7 of \cite{KL} and section 1.16 of \cite{AF} for a proof) 
is the following bound.

\begin{prop}
\label{reduce}
 Let $g\in\cG$ have support included in  $\L_s$.  Then for any
disorder configuration $\a$, any $\g>0$  and any family of
functions $F:=\{f_x\}_{x\in\toromi}$ on $\O_\ep$,
\begin{equation}
 \label{mrfridge}
 \begin{split}
   \sup spec_{L^2(\mu_\ep)}&\{\ep^{-1}\Av_{x\in\toromi}(\t_xg\,f_x)+\g\ep^{d-2}\cL\}\leq\\
   & \ep^{-1}c(g,\ninf{F})\ninf{\nabla F}+\sup
   spec_{L^2(\mu_\ep)}\{c(g)\Av_{x\in\toromi}f_x^2+\ov 2\g\ep^{d-2}\cL\},
\end{split}
\end{equation}
where $\ninf{F}:=\sup_{x\in\toromi}\|f_x\|_\infty $ and
$\ninf{\nabla F}:=\sup_{x\in\toromi} \sup_{b\sset\L_{x,s}} \|\nabla_b f_x\|_\infty$.
\end{prop}

In the space $\cG$ it is also possible to introduce a $H_{-1}$ norm closely
related to that given by perturbation theory (see proposition
\ref{jamme} above). \\
Given positive integers $\ell,s$ with $s^2\leq\ell$ and
$f,g\in\cG$ with $\D_f,\D_g\sset\L_s$, for any canonical or grand
canonical Gibbs measure $\mu$ on $\L_\ell$ we define
\begin{equation}
 V_\ell(f,g;\mu) :=
(2l)^{-d}\mu\Bigl(\sum_{|x|\leq\ell_1}\t_x f,
\bigr(-\cL_{\L_\ell}\bigr)^{-1}\sum_{|x|\leq\ell_1}\t_x g\Bigr).
\label{V_ell}
\end{equation}
If $\L_\ell$ is replaced by $\L_{z,\ell}$ and the above sums are over $x\in\L_{z,\ell_1}$  we will simply write
$V_{z,\ell}(f,g;\mu)$ and if $f=g$ we write
$V_{\ell}(g;\mu)$ or $V_{z,\ell}(g;\mu)$.\\
It is simple to check that $V_\ell(g;\mu)$ can be variationally
characterized as follows:
\begin{equation}
\label{maria}
\begin{split}
V_\ell(g;\mu)&=(2l)^{-d}\sup_h
\bigl\{2\mu\bigl(\sum_{|x|\leq\ell_1}\t_x g,h\bigr)-\cD_{\L_\ell}
(h;\mu)\bigr\}\\
&= (2l)^{-d}
\sup_h\frac{\mu(\sum_{|x|\leq\ell_1}\t_x g,h)^2}{\cD_{\L_\ell}(h;\mu)}
\end{split}
\end{equation}
where $\sup_h$ is taken among the non constant functions with support
contained in $\L_\ell$.\\
The variational characterization allows one to derive some simple
bounds on $V_\ell(g;\mu)$. Let $\D$ be a box such that $\D_g\sset
\D\sset\L_s$ and for any $x \in \ZZ ^d$ let $\cF_x$ be the $\s$--algebra
generated by $m_{\D_x}$ and $\{\h_y\}_{y\notin\D_x}$.  Then, for any
function $h$,
\begin{gather*}
 \mu(\t_x g,h)=\mu\bigl(\,\mu(\t_x g;h\tc\cF_x)\bigr)\leq\mu
\bigl(\,\Var_\mu (\t_x g\tc\cF_x)^{\ov 2}\Var_\mu (h\tc \cF_x)^{\ov
  2}\bigr)\\
  \leq \mu\bigl(\,\Var_\mu (\t_x g\tc\cF_x)\bigr)^{\ov 2}
       \mu\bigl(\,\Var_\mu (h \tc\cF_x)\bigr)^{\ov 2}
\end{gather*}
which implies that
\begin{equation}
\label{myskin}
  \mu(\sum_{|x|\leq\ell_1}\t_x g,h)^2 \le
c  \sum _{|x|\leq\ell_1}\mu\bigl(\,\Var_\mu (\t_x
g\tc\cF_x)\bigr)
\sum_{|x|\leq\ell_1}\mu\bigl(\,\Var_\mu (h\tc\cF_x)\bigr).
\end{equation}
If we appeal now to the Poincar\'e inequality
$$
\Var_\mu (h\tc\cF_x) \le c s^2 \sum_{b\sset \D_x}\mu\bigl(c_b (\nabla_b h)^2
\tc \cF_x\bigr),
$$
the last sum in (\ref{myskin}) is bounded by $c\,s^{d+2}\cD_{\L_{\ell}}(h;\mu)$. Recalling (\ref{maria}), for any $\ell > s^2$
we finally get
\begin{equation}
\label{montag}
V_\ell(g;\mu)\leq c\,s^{d+2}\Av_{|x|\leq\ell_1}\mu\bigl(\,\Var_\mu(\t_x g\tc
\cF_x)\,\bigr).
\end{equation}
In particular
\begin{equation}
\label{aglaja}
V_\ell(g;\mu)\leq c\, s^{d+2}\|g\|_\infty^2.
\end{equation}
In order to benefit of the ergodicity of the random field, it is natural
to define, for any $m\in (0,1)$ and any $g\in \cG$,
\begin{equation}
V_m(g):=\lim _{\ell\uto\infty}
(2\ell)^{-d}\,\EE\Bigl[\mu^{\la_0(m)}\bigl(\,\sum_{|x|\leq \ell_1}\t_x g,
(-\cL_{\L_\ell})^{-1}
\sum_{|x|\leq \ell _1}\t_x g\bigr)\Bigr]
\label{V_m}
\end{equation}
where, we recall, $\l_0(m)$ is the annealed chemical potential
corresponding to the particle density $m$.
If $m=0,1$ we simply set $V_m(g)=0$ for any $g\in\cG$. In section
\ref{CLTV} we will prove, among other results, that the limit appearing
in (\ref{V_m})
exists finite and that it defines a semi--inner product on $\cG$ (see
theorem \ref{bingo} there). With this definition we have the following
result.
\begin{lemma}
\label{carlosmontoya}
Let $g\in \cG$. Then
\begin{equation}
\label{schutz}
\limsup _{\ell\uto\infty,\,\ep\dto 0} \Av_{|x|\leq \frac{1}{\ep}}
\sup_{\nu\in \cM(\L_{x,\ell})} V_{x,\ell}\bigl(g;\nu\bigr)\leq
\sup_{m\in[0,1]}V_m(g).
\end{equation}
\end{lemma}
\begin{proof}
As in \cite{KL}, chapter 7, lemma 4.3, we introduce a scale parameter $k$, with
$k\uto \infty$ after $\ell\uto\infty$,  and partition $\L_\ell$ in cubes of
side $2k+1$. More precisely, we define $\L^{(k)}_\ell:=\L_\ell \cap (2k+1)\ZZ^d$ and
write $\L_\ell=B_{k,l}\cup\bigl(\cup_{z\in \La_\ell^{(k)}} \La_{z,k}\bigr)$ where
 $B_{k,\ell}:=\L_\ell\setminus\cup_{z\in\L^{(k)}_\ell}\L_{z,k}$. Then,
by proceeding as in \cite{KL} and by using the variational
characterization (\ref{maria}) together with
the integration by parts formula (\ref{int_by_parts}),
for any $\nu \in \cM(\L_\ell)$ we get
\begin{equation}
\label{reichstag}
 V_\ell(g;\nu)\leq (2\ell)^{-d} \sup_{\underline{h}}\bigl\{ \sum_{z\in\L_\ell
^{(k)}}F_z(h_z;\nu)\bigr\}+c(g) \sqrt{k\ell^{-1}+k^{-\ov 2}}\,\bigr\}
\end{equation}
where $c(g)$ is as in (\ref{silenziosanotte}),
$
F_z(h_z;\nu):=2\sum_{y\in\L_{z,k_1}}\nu(\t_y g,h_z)-\cD_{\L_{z,k}}(h_z;\nu)
$
and the supremum $\sup_{\underline {h}}$ is taken over all families
$\underline{h}=\{h_z\}_{z\in\L_\ell^{(k)}}$ such that $h_z$ depends only
on $\h_{\L_{z,k}}$ and $\cD(h_z;\nu) \leq c(g) k^d$. \\
Actually it is simple to
check that in (\ref{reichstag}) we can restrict the supremum to families
$\underline{h}$ that satisfy the extra condition $\|h\|_\infty \leq c(g)c_k$
for some constant $c_k$ depending on $k$.\\
 Therefore, if $m$ is the
particle density associated to the canonical measure $\nu$ and thanks to the
equivalence of ensembles (see lemmas \ref{equi} and \ref{vespe}),
for any disorder configuration $\a$ we get
\begin{equation*}
\begin{split}
  & \bigl| \sum_{z\in\L_\ell^{(k)}} F_z(h_z;\nu)-\sum_{z\in\L_\ell^{(k)} }
  F_z(h_z;\mu_{\L_\ell}^{\la(m)})\bigr |\leq c(g)c_k,\\
  & \bigl|\sum_{z\in\L_\ell^{(k)} } F_z(h_z;\mu_{\L_\ell}^{\la(m)}) -
  \sum_{z\in\L_\ell^{(k)}} F_z(h_z;\mu^{\l_0(m)})\bigr|\leq
  c(g)c_k\,\ell^d \,\bigl|m-\mu^{\l_0(m)}(m_{\L_\ell})\bigr|.
\end{split}
\end{equation*}
Thanks to the previous observations we finally obtain
\begin{equation*}
 \Av_{|x|\leq \frac{1}{\ep}}
\sup_{\nu\in \cM(\L_{x,\ell})} V_{x,\ell}\bigl(g,\nu\bigr) \leq
 c \sqrt{k\ell^{-1}+k^{-\ov 2}}+
c_k \ell^{-d}+c_k \Theta_1 + c_{k,\ell} \Theta _2
\end{equation*}
where $c_{k,\ell}$ is a positive constant depending on $k,\ell$ such that
$\lim_{k\uto \infty,\ell\uto\infty} c_{k,\ell}=1$ and
\begin{equation*}
\begin{split}
& \Theta_1:=\Av_{|x|\leq\frac{1}{\ep}} \sup_{m\in[0,1]}
\bigl|m-\mu^{\l_0(m)}(m_{\L_{x,\ell}})\bigr| ,\\
& \Theta_2:=\Av_{|x|\leq\frac{1}{\ep}}\sup_{m\in[0,1]}\tau_x\Bigl(
\sup_{\underline{h}}\bigl\{(2k)^{-d}\Av_{z\in\L_\ell^{(k)}}
F_z(h_z,\mu^{\l_0(m)})\bigr\} \Bigr),\\
\end{split}
\end{equation*}
and $\sup_{\underline{h}}$ is as before.

\smallno It is clear that, by considering a fixed density $m$ in the
definition of $\Theta_1$ and $\Theta_2$, for almost all disorder
configurations $\a$, $\Theta_1$ is negligible as $\ell\uto\infty,\ep\dto
0$. Moreover, because of the ergodicity of the random field $\a$ and of
the variational characterization (\ref{maria}), it is also clear that
for almost all disorder configurations $\a$
\begin{equation*}
\limsup_{\ell\uto\infty, \ep \dto 0}\Theta _2 \leq \EE\bigl
(V_k(g;\mu^{\l_0(m)})\bigr)
\end{equation*}
To handle the supremum over $m\in [0,1]$ requires some simple additional
observations based on compactness of $[0,1]$ and lemma \ref{vespe} (see
e.g section 1.13  in \cite{AF}).
\end{proof}

\subsection{Back to the proof of proposition \ref{wannsee}}\label{urca}
Given the technical tools developed in the previous paragraph, let us
now complete the proof of proposition \ref{wannsee} modulo some non
trivial results to be proved later on.

\smallno The basic idea would be to benefit of the ergodicity of the
model by means of the \emph{localization technique} discussed in
subsection \ref{lete}. Unfortunately, the function $\bar
J_{b,a,\ell,\ep}^g$  appearing in (\ref{faust}) cannot be written as
$\Av_{x\in\toromi}f_x$ (or as a more complex spatial average) for
suitable functions $f_x$ having support independent of $\ep$. We will
need some subtle techniques developed for non-gradient systems in order
to approximate $\bar J_{b,a,\ell,\ep}^g$ with such a spatial average.
There is however one piece of $\bar J_{b,a,\ell,\ep}^g$, namely the
density ``gradient'' $
\bigl(2b/\ep\bigr)^{-1}\bigl[m_{x+\frac{b}{\ep}e',\frac{a}{\ep}}
-m_{x-\frac{b}{\ep}e',\frac{a}{\ep}}\bigr] $ which can be conveniently
written as a suitable spatial average. To this aim recall the
definition (\ref{dinoce}) of the spatial average $\Av^{\ell,s}_{z,y}$
and define for any
particle configuration $\h$, $m_\ell^{1,e},\,m_\ell^{2,e}$ and $m_\ell^{e}$
to be the particle density associated to $\h$ in the sets $\L^{1,e}_\ell,
\L^{1,e}_\ell$ and $\L_\ell^e$ defined in (\ref{blocchetto})
respectively. It is then simple to check the following identity (which motivates the
introduction of $\Av^{\ell,s}_{z,y}$):
\begin{equation}
\label{merlino}
\Av^{\ell,s}_{z,y} \t_z \frac{m^{2,e}_\ell-m^{1,e}_\ell}{\ell}=
\t_y\frac{m^{2,e}_s-m^{1,e}_s}{s}.
\end{equation}
Let now
$n,\frac{a}{\ep},\frac{b}{\ep}$ be odd integers such that
$\frac{a}{n\ep}\in\NN$ and $\frac{b}{a}\in\NN$. Then, it is simple
to check that
\begin{equation}\label{maya}
\Av_{u=0}^{\frac{2b}{a}-1} \t_{x_u}\frac{m^{2,e}_\frac{a}{\ep}-
m^{1,e}_\frac{a}{\ep}}{a/\ep}=\frac{m_{x+\frac{b}{\ep}e,\frac{a}{\ep}}
-m_{x-\frac{b}{\ep}e,\frac{a}{\ep}}}{2b/\ep}
\end{equation}
where
$$
x_u:=x+\Bigl(
u\frac{a}{\e}-\frac{b}{\ep}+ \frac{1}{2}\bigl(\frac{a}{\ep}-1\bigr)
+1\Bigr)e.
$$
Therefore, if we define
\begin{equation}
\label{fogliegialle}
\Av_{z,x}^\ast\,f_z:=\Av_{u=0}^{\frac{2b}{a}-1}\Av_{z,x_u}
^{n,\frac{a}{\ep}} f_z
\end{equation}
(when  necessary we will also add
the versor $e\in\cE$ into the notation by writing  $\Av_{z,x}^{\ast,e}$),
 thanks to (\ref{merlino}) and (\ref{maya})  we obtain:
\begin{equation}
\label{perleaiporci.K.V.}
\Av_{z,x}^\ast\frac{m^{2,e}_n-m^{1,e}_n}{n}
=\frac{m_{x+\frac{b}{\ep}e,\frac{a}{\ep}}
-m_{x-\frac{b}{\ep}e,\frac{a}{\ep}}}{2b/\ep}\,.
\end{equation}
If the above conditions on $n,\frac{a}{\ep},\frac{b}{\ep}$ are not
satisfied, we extend the definition of $\Av_{z,x}^\ast$ by replacing in
(\ref{fogliegialle}) $\frac{a}{\ep}$, $\frac{b}{\ep}$, $\frac{2b}{a}$
with $r_1$, $r_2$ and $\frac{2r_2}{r_1}$ respectively, where $r_1$ is
the smallest odd number in $n \ZZ$ such that $\frac{a}{\ep}\leq r_1$ and
$r_2$ is the smallest odd number in $r_1\ZZ$ such that
$\frac{b}{\ep}\leq r_2$.

\smallno {\bf Warning.}  In the sequel, for the sake of simplicity we
will always assume $n,\frac{a}{\ep},\frac{b}{\ep}$ to be odd integers
such that $\frac{a}{n\ep}\in\NN$ and $\frac{b}{a}\in\NN$. The way to
treat the general case is shortly discussed right after section
\ref{apriticielo}.

\smallno
It is convenient to introduce also  $\Av_{z,x}^\star$ defined as
the dual average of $\Av_{z,x}^\ast$, \ie
\begin{equation}
\label{dualaverage}
\Av_{x\in\toromi}\Bigl( f_x \bigl (\Av_{z,x}^\ast g_z)\Bigr)=\Av_{x\in\toromi}
\Bigl(g_x \bigl (\Av_{z,x}^\star f_z)\Bigr).
\end{equation}
The explicit formula of $\Av_{z,x}^\star f_z$ can be easily computed and it is
similar to the formula of $\Av_{z,x}^\ast f_z$.
\\

We introduce at this point some special functions related to
the gradient of the density field.  Given two integers
$0\leq n \leq s$, $e\in\cE$  and a grand canonical measure $\mu$ on an
arbitrary set $\L$ containing $\L_s^e$, we write
\begin{equation}
\label{estasi}
m_n^{2,e}-m_n^{1,e} = \psi_{n,s}^e + \phi_{n,s}^e\,,
\quad \text{ with }\quad\,\phi_{n,s}^e :=\mu\bigl[\,m_n^{2,e}-m_n^{1,e} \tc
\cF^e_s\,\bigr]\, ,
\end{equation}
where $\cF^e_s$ is the $\s$--algebra generated by $m_s^e$.
Notice that, in absence of disorder, the function $\phi^e_{n,s}$ would be
identically equal to zero and that $\psi^e_{n,s}\in \cG$ for all $n<s$,
since $\nu(\psi^e_{n,s})=0$ for all $\nu \in \cM(\L)$ and all $\L$
containing $\L_s^e$.
Thanks to (\ref{montag}) with  $\D:=\L_n^e$ and
$s:=n$ and thanks to the equivalence of ensembles (see lemma
\ref{cleonimo}), given $\ell \geq n^2$ it is easy to check that
\begin{equation}
\label{kolja}
 V_\ell\bigl(\frac{\psi_{n,n}^e}{n}; \nu\bigr)\leq c \quad \forall \nu \in \cM (\L_\ell), \quad \qquad   V_\ell\bigl(\frac{\psi_{n,n}^e}{n}; \mu^{\la_0(m) } \bigr)\leq c\, m(1-m).
\end{equation}
Using  decomposition (\ref{estasi}) we can now write $\bar J_{b,a,\ell,\ep}^g$ as
$$
\bar J_{b,a,\ell,\ep}^g = \sum_{j=0}^5\Av_{x\in\toromi}J(\ep
x)\,\psi_x^{(j)}
$$
where (we omit in the notation the suffix $b,a,\ell,\ep,g$ )
\begin{equation*}
\begin{split}
 & \psi^{(0)}_x:=\Av_{y:|y-x|\leq \ell_1}\Bigl[
j_{y,y+e}+\t_y\cL g+\sum_{e'\in\cE}
D_{e,e'}(m_{x,\ell})\t_y\frac{\psi_{n,n}^{e'}}{n}\Bigr],\\
 & \psi^{(1)}_x:=\sum_{e'\in\cE} D _{e,e'}(m_{x,\ell})
 \Bigl[ \t_x\frac{\psi_{n,n}^{e'}}{n}-\Av_{y:|y-x|\leq \ell_1}\t_y\frac{\psi_{n,n}^{e'}}{n}\Bigr ]\\
 & \psi^{(2)}_x:=\sum_{e'\in\cE}\Bigl[D_{e,e'}
 (m_{x,\frac{a}{\ep}})-D_{e,e'}(m_{x,\ell})\Bigr]\t_x\frac{\psi_{n,n}^{e'}}{n}\\
& \psi^{(3)}_x :=\sum_{e'\in\cE}D_{e,e'}(m_{x,\frac{a}{\ep}})\Bigl[\Av
^{\ast,e'}_{z,x}\t_z\frac{\psi_{n,n}^{e'}}{n}-\t_x\frac{\psi_{n,n}^{e'}}{n}\Bigr]\\
& \psi^{(4)}_x:=\sum_{e'\in\cE}D_{e,e'}(m_{x,\frac{a}{\ep}})\Bigl[
\frac{m_{x+\frac{b}{\ep}e',\frac{a}{\ep}}-m_{x-\frac{b}{\ep}e',\frac{a}{\ep}}}{2b/\ep}
-\Av^{\ast,e'}_{z,x}\t_z\frac{m_n^{2,e'}-m_n^{1,e'}}{n}
\Bigr]\\
&\psi^{(5)}_x:=\sum_{e'\in\cE}D_{e,e'}(m_{x,\frac{a}{\ep}})\Av^{\ast,e'}
_{z,x}\t_z\frac{\phi_{n,n}^{e'}}{n}
\end{split}
\end{equation*}
and we define
\begin{equation*}
 \O_j:=\sup spec _{L^2(\mu_\ep)}\bigl\{\ep^{-1}\Av_x J (\ep x)\psi^{(j)}_x+\g\ep^{d-2}\cL_\ep\bigr\}\quad j=0,\dots,5\,.
\end{equation*}
Then, thanks to the subadditivity of ``$\sup spec$", proposition
\ref{wannsee} follows from the next result.
\begin{prop}\label{pedro}
  Let $d\geq 3$ and $\g>0$. Then, for almost any disorder configuration
  $\a$,
\begin{equation}
  \inf _{g\in\GG} \limsup_{n\uto \infty,\ell\uto\infty,\ep\dto 0}
  \sup _J\O_0 \leq 0
\label{bubu!!}
\end{equation}
and, for any $j=1,\dots,5$,
\begin{equation}
  \label{jojo}
   \limsup_{n\uto\infty,b\dto 0,a\dto 0,\ell\uto\infty,\ep\dto
  0}\;\sup_J\;\O_j \leq 0
\end{equation}
where $J$ varies in  $\{J\in C(\toroma)\,:\,\|J\|_\infty\leq 1\}$.
\end{prop}
The proof of proposition \ref{pedro} is best divided into several
pieces according to the value of the index $j$.

\subsection{The term $\mathbf{\O_0}$}
\label{pietro}
Let us first prove (\ref{bubu!!}).  By localizing on cubes of side
$2\ell+1$ (see (\ref{cappy})) and using the regularity of
$J(\cdot)$, it is enough to prove that for
almost any disorder configuration $\a$,
\begin{equation}
\label{cocktail}
\inf_{g\in\GG}\;\limsup _{n\uto\infty,\ell\uto\infty,\ep\dto 0}\;
 \Av_{x\in\toromi}\sup_{|\b|\leq 1}\,
\sup_{m}\underset{L^2(\nu_{\L_{x,\ell},m})}{\sup spec} \bigl\{
\ep^{-1}\b\,\Av_{y:|y-x|\leq\ell_1}\t_y\psi_m^{(n,g)} +c\ell^{-d}\ep^{-2}\cL_{\L_{x,\ell}}\bigr\}\leq 0
\end{equation}
where
$$
\psi_m^{(n,g)}:=j_{0,e}+\cL g+\sum_{e'\in\cE} D_{e,e'}(m)\frac{\psi_{n,n}^{e'}}{n}.
$$
Since $\ep\dto 0$ before $\ell \uto\infty$ and since for any $\ell$
large enough, any $|y-x|\le \ell_1$ and any
$\nu\in \cM(\L_{x,\ell})$, $\nu\bigl(\t_y\,\psi_m^{(n,g)}\bigr) =0$, we can appeal to
perturbation theory (see proposition \ref{jamme})
and conclude that it is
enough to show
that
\begin{equation}
\label{lowen}
\inf_{g\in\GG}\lim _{n\uto\infty,\ell\uto\infty,\ep\dto 0} \Av_{x\in\toromi}
\sup_{m\in[0,1]} V_{\L_{x,\ell}}\bigl(\psi_m^{(n,g)},\nu_{\L_{x,\ell},m})= 0
\end{equation}
where $V_{x,\ell}$ has been defined right after (\ref{V_ell}).
A minor modification of the proof of lemma \ref{carlosmontoya} shows
that (\ref{lowen}) follows from
\begin{equation}
\label{mormoni}
\inf_{g\in\GG}\limsup _{n\uto\infty}
\sup_{m\in[0,1]}V_m\bigl (\psi_m^{(n,g)}\bigr)=0 \qquad \forall d\geq3
\end{equation}
(see (\ref{V_m}) for the definition of $V_m$) which, in turn, follows
from theorem \ref{agognatameta}.

\subsection{The three terms $\mathbf{\O_1,\,\O_2,\,\O_3}$}\label{liberi}
Let us prove (\ref{jojo}) for $j=1,2,3$.
In what follows, by means of proposition \ref{reduce},
 we will reduce the eigenvalues estimate  $\O_1,\;\O_2$ and $\O_3$ to the
Two Blocks estimate (see subsection \ref{rachmaninov}). To this aim, by integrating by parts,
we can write
\begin{equation*}
 \ep^{-1}\Av_x J(\ep x)\psi_x^{(j)} =
\sum_{e'\in\cE}\ep^{-1}\Av_x\t_x \frac{\psi_{n,n}^{e'}}{n}\cdot B_x^{(j)}\quad\forall \,j=1,2,3
\end{equation*}
where
\begin{align*}
& B _x ^{(1)}: =J(\ep x)D_{e,e'}(m_{x,\ell})-\Av_{y:|y-x|\leq\ell_1}J(\ep y) D_{e,e'}(m_{y,\ell})\\
& B_x^{(2)}:=J(\ep x)\bigl[D_{e,e'}(m_{x,\frac{a}{\ep}})-D_{e,e'}(m_{x,\ell})\bigr]\\
& B_x^{(3)}:=\Av^{\star,e'}_{z,x}
J(\ep z) D_{e,e'}\bigl(m_{z,\frac{a}{\ep}}\bigr)-
J(\ep x)D_{e,e'}\bigl(m_{x,\frac{a}{\ep}}\bigr) .
\end{align*}
Notice that, for any $b\sset\L_{x,n}$,
$$
\nabla_b B_x^{(1)}=\nabla_b B_x^{(2)}=0, \quad
|\nabla_b B_x^{(3)}|\leq c\,n\frac{\ep}{a}\text{Osc}\bigl(D,c\frac{\ep^d}{a^d}
\bigr).
$$
Therefore, using proposition \ref{reduce}, it is enough to prove that for almost any disorder
configuration $\a$, given $\g>0$,
\begin{equation}
\label{bonaqa}
 \limsup_{b\dto 0,a\dto 0,\ell\uto\infty,\ep\dto 0}\sup_J\,\sup spec _{L^2(\mu _\ep)}\bigl\{\Av_{x\in\toromi} \bigl(B_x^{(j)}\bigr)^2+\ov2\g\ep^{d-2}\cL_\ep\bigr\}=0\qquad\forall j =1,2,3.
\end{equation}
Since $D$ can be approximated by Lipschitz functions and $J$ is smooth,
 (\ref{bonaqa})
can be derived from the Two Blocks estimate  (see subsection \ref{rachmaninov}). For simplicity of notation,
let us consider the case $j=2$ (the case $j=1$ is simpler, while $j=3$
is a slight variation) and $D$ Lipschitz continuous. Since
$\bigl(B_x^{(2)}\bigr)^2\leq c\,\bigl
|\,m_{x,\ell}-m_{x,\frac{a}{\ep}}\,|$, by introducing a scale parameter
$k$ such that $k\uto \infty$ after $a\dto 0,\ell\uto\infty$ and $\ep\dto
0$, we can estimate
\begin{equation*}
\bigl(B_x^{(2)}\bigr )^2\leq
c\,\Av _{|y|\leq\ell}\Av_{|z|\le \frac{a}{\ep}}|m_{x+y,k}-m_{x+z,k}| +
c\,\frac{k}{\ell}+c\frac{k}{a/\ep}.
\end{equation*}
At this point, by the subadditivity (\ref{koukaburra}) of ``$\sup spec$'', the thesis follows from
the Two Blocks estimate.

\subsection{The term $\mathbf{\O_4}$}\label{apriticielo}
The proof of (\ref{jojo}) for $j=4$ is based on the Two Blocks estimate.
Notice that, thanks to (\ref{perleaiporci.K.V.}), the function $\psi_x
^{(4)}$ entering in the definition of $\O_4$ is either identically equal
to zero if $n,\frac{a}{\ep},\frac{b}{\ep}$ are odd integers such that
$\frac{a}{n\ep}\in\NN$ and $\frac{b}{a}\in\NN$, or it can be written as
\begin{equation}\label{romania}
\psi_x^{(4)}=\sum_{e'\in\cE}D_{e,e'}(m_{Q_{x,\frac{a}{\ep}}})
\Bigl[\frac{m_{x+\frac{b}{\ep}e',\frac{a}{\ep}}-m_{x-\frac{b}{\ep}e',\frac{a}{\ep}}}
{2b/\ep}  -
\frac{m_{x+r_2 e',r_1}-m_{x-r_2 e',r_1}}{2r_2 }
\Bigr]
\end{equation}
where $r_1,r_2$ have been defined in subsection \ref{urca}.  By the Two
Blocks estimate it is simple to check that for any $\g>0$ and for almost
any disorder configuration $\a$
\begin{align}
& \lim_{a\dto 0, \ep\dto 0} \sup spec _{L^2(\mu_\ep)}\bigl\{\Av_{x\in\toromi}
\bigl | m_{x,\frac{a}{\ep} }-m_{x,r_1} |+\g\ep^{d-2}\cL_\ep\bigr\}=0
\label{vladimir}\\
& \lim_{a\dto 0, \ep\dto 0}\sup_{|w|\leq 2\frac{a}{\ep}} \sup spec
_{L^2(\mu_\ep)}\bigl\{\Av_{x\in\toromi} \bigl | m_{x,\frac{a}{\ep}
}-m_{x+w,\frac{a}{\ep} } |+\g\ep^{d-2}\cL_\ep\bigr\}=0
\label{vissotski}
\end{align}
(hint: introduce the scale parameter $k$ with $a\dto 0, k \uto
\infty,\ep\dto 0$ and write
\hbox{$ m_{x,s}= \Av_{y\in \L_{x,s}}m_{y,k} + O(k/s)$ for $s=\frac{a}{\ep}, r_1$}). \\
In (\ref{romania}) we can substitute $r_1$ by $ \frac{a}{\ep}$ (thanks
to (\ref{vladimir})) and after that in the numerators we can substitute
$r_2$ by $\frac{b}{\ep}$ (thanks to
(\ref{vissotski})).
In order to conclude is enough to observe that
$\ep^{-1}\bigl|\frac{1}{b/\ep}-\frac{1}{r_2}\bigr|\leq c\,\frac{a}{b^2}$
which goes to $0$.

\subsection{The term $\mathbf{\O_5}$}
The proof of (\ref{jojo}) for $j=5$ is based on the key results of section
\ref{Fluct} and it is one place where the restriction on the
dimension $d\ge 3$ is crucial for us.
We refer the reader to the beginning of
section \ref{Fluct} for an heuristic justification of the above
condition. Here it is enough to say that the main contribution to the
term $\O_5$ comes from the fluctuations in the density field induced by
the fluctuations of the \emph{disorder field}.

\smallno
By the subadditivity of "$\sup spec$" we only need to prove
that for almost all $\a$, given $e,e'\in \cE$ and $\g>0$,
\begin{equation}
\label{urania1000}
\limsup_{n\uto\infty,\,b\dto 0,\,a\dto 0,\,\ep\dto 0} \sup_{J}\;\sup
spec_{L^2(\mu_\ep)}\bigl\{\ep^{-1}\Av_x
J (\ep x)D_{e,e'}(m_{x,\frac{a}{\ep}})
\Av^{\ast,e'}_{z,x}\t_z\frac{\phi_{n,n}^{e'}}{n}+\g\ep^{d-2}\cL_\ep\bigr\}
\leq 0.
\end{equation}
Recall the definition of $\Av^{\ast,e'}_{z,x}$ and $x_u$ given in
(\ref{fogliegialle}). Then, thanks again to the subadditivity of "$\sup spec$",
the "$\sup spec$" in the  l.h.s. of (\ref{urania1000}) is bounded from above by
\begin{equation}
\label{detlef}
\Av_{u=0}^{\frac{2b}{a}-1}
\sup spec_{L^2(\mu_\ep)}\bigl\{\ep^{-1}\Av_x
J (\ep x)D_{e,e'}(m_{x,\frac{a}{\ep}})
\Av^{n,\frac{a}{\ep}} _{z,x_u}\t_z\frac{\phi_{n,n}^{e'}}{n}+\g\ep^{d-2}\cL_\ep
\bigr\}.
\end{equation}
Observe that $\Av^{n,\frac{a}{\ep}}
_{z,x_u}\t_z\frac{\phi_{n,n}^{e'}}{n}$ has support inside
$\L^{e'}_{x_u,\frac{a}{\ep}}$. We would like at this point to localize on
boxes of side length of order $O\bigl(\frac{a}{\ep}\bigr)$ in such a way
that $D_{e,e'}(m_{x,\frac{a}{\ep}})$ becomes a constant. To this aim,
given $u\in\{0,\dots,\frac{2b}{a}-1\}$ and $x\in\toromi$, we set
\begin{equation*}
\D_{x,u}:=
\begin{cases}
Q_{x,10\frac{a}{\ep}} & \text{if }Q_{x,\frac{a}{\ep}}\cap
\La_{x_u,2\frac{a}{\ep}}\not=\emptyset\\
Q_{x,\frac{a}{\ep}} &\text{otherwise}
\end{cases}.
\end{equation*}
and we observe that either
$\D_{x,u}$ is disjoint from or completely contains $\L_{x_u,2\frac{a}{\ep}}$.
Therefore, if in (\ref{detlef})
we could replace the term $D_{e,e'}(m_{x,\frac{a}{\ep}})$ by the new
term
$D_{e,e'}(m_{\D_{x,u}})$, then it would be simple to check (by localizing on
boxes $\La_{x_u, 2\frac{a}{\ep}}$) that all what is needed is that for $d\geq
3$, for all $T\in \bbN$ and for almost all $\a$,
\begin{equation}
\label{guccini}
 \limsup_{n\uto\infty,a\dto 0,\ep\dto 0}
\Av_{x\in\toromi}\sup_{|\b|\leq T}
\sup_{\nu\in \cM(\L_{x,2\frac{a}{\ep}})} \; \sup spec_{L^2(\nu)}\bigl\{\ep^{-1}
\b \Av^{n,\frac{a}{\ep}} _{z,x}\t_z
\frac{\phi_{n,n}^{e'}}{n}+\ep^{-2}\Av_{b\in\L_{x,2a/\ep}}\cL _b
\bigr\}\leq 0
\end{equation}
Section \ref{Fluct} is devoted to the proof of
(\ref{guccini}) (see theorem \ref{pte} there).
\\
Therefore,
it remains to prove that for $d\geq
3$, for almost all $\a$ and for any $\g>0$
\begin{equation}
\label{urania}
\begin{split}
\limsup_{n\uto\infty,b\dto 0,a\dto 0,\ep\dto 0} &\sup_J\,\Av_{u=0}
^{\frac{2b}{a}-1}\sup spec_{L^2(\mu_\ep)} \bigl\{\ep^{-1}\Av_{x\in\toromi}
J(\ep x) \\
\times \, &\bigl[D_{e,e'}(m_{x,\frac{a}{\ep}})-D_{e,e'}(m_{\D_{x,u}})\bigr]
\Av^{n,\frac{a}{\ep}}_{z,x_u}\t_z\frac{\phi_{n,n}^{e'}}{n}+\g\ep^{d-2}
\cL_\ep\bigr\}\leq 0.
\end{split}
\end{equation}
Notice that the only values of $u$ which contribute to the
$\Av_{u=0} ^{\frac{2b}{a}-1}$ above, in what follows called ``bad
values'', are those for which
$Q_{x,\frac{a}{\ep}}\neq \D_{x,u}$ for some $x\in\toromi$. It is easy to
check that the cardinality of the bad values of $u$ is of order
$O(1)$  for any fixed $x\in\toromi$. Thus we only need to bound the "$\sup
spec$" appearing in (\ref{urania}) by $o(\frac{b}{a})$, uniformly in $u$
in the bad set.
Thanks to (\ref{perleaiporci.K.V.}) and
(\ref{estasi}) we can write
\begin{equation}
\label{tunisia}
\Av^{n,\frac{a}{\ep}}_{z,x_u}\t_z\frac{\phi_{n,n}^{e'}}{n}
=\frac{m_{x_u+\frac{b}{\ep},\frac{a}{\ep}}- m_{x_u-\frac{b}{\ep},\frac{a}{\ep}} }{2b/\ep}
-\Av^{n,\frac{a}{\ep}}_{z,x_u}\t_z\frac{\psi_{n,n}^{e'}}{n}
\end{equation}
Then, the contribution in (\ref{urania}) coming from the first addendum
in the r.h.s. of (\ref{tunisia}) is not larger than $O(\frac{1}{b})$ and
therefore negligible.
\\
Let us consider the contribution of the second
addendum. An integration by parts shows that
$$
\Av_{x\in\toromi}
J(\ep x)\bigl(D_{e,e'}(m_{x,\frac{a}{\ep}})-D_{e,e'}(m_{\D_{x,u}}\bigr)
\Av^{n,\frac{a}{\ep}}_{z,x_u}\t_z\frac{\psi_{n,n}^{e'}}{n} = \Av_{x\in\toromi}
\t_x\frac{\psi_{n,n}^{e'}}{n} B_{x,u}
$$
where the functions $B_{x,u}$ satisfy $\|B_{x,u}\|\leq c$ together with
$$
|\nabla_b B_{x,u}|\leq c\frac{n\ep}{a} Osc (D_{e,e'},
c\frac{\ep^d}{a^d}) \quad \forall b \in \La_{x,n}^{e'}.
$$
Moreover, $B_{x,u}$ is a
particular spatial average (dual to $\Av^{n,\frac{a}{\ep}}_{z,x_u}$) of
$
J(\ep z)\bigl(D_{e,e'}(m_{z,\frac{a}{\ep}})-D_{e,e'}(m_{\D_{z,u}}\bigr)$.
Therefore, by proposition \ref{reduce} and the Two Blocks estimate  (see subsection \ref{rachmaninov}), the contribution of the second addendum is also negligible (see also the discussion at the end of subsection \ref{liberi}).


\subsection{Proof of the energy estimate}
\label{gatorade11}
In this subsection we prove lemma \ref{teti}.  It is simple to check
that
\begin{equation}
\label{sasha}
\text{spatial--time average in }(\ref{jerusalem})=
\sup _{H\in\cH_b } H^\star _{b,a,\ep}
\end{equation}
 where $\cH_b:=\{H\text{ smooth on } [0,T]\times\toroma,
\;\|H\|_\infty\leq\ov b \}$ and 
\begin{equation*}
\label{eq:energy_estimate_0}
H^\star_{b,a,\ep}:=\Av_{x\in\toromi}\int_0^T\Bigl(2 H(s,\ep x)
\Bigl[\frac{m_{x+\frac{b}{\ep}e,\frac{a}{\ep}}(s)-
 m_{x,\frac{a}{\ep}}(s)}{b}\Bigr]-H(s,\ep x)^2\Bigr)ds.
\end{equation*}
In what follows let $H$ belong to  $\cH_b$.  By the entropy inequality and
the Feynman-Kac formula (see (\ref{ent_ineq}) and (\ref{eq:entropy1000})\,), for any $\g>0$,
\begin{equation}
\label{corsica}
\begin{split}
 \EE^{\mu^\ep}\bigl( H^\star _{2b,a,\ep}\bigr)&\leq
 \frac{\k}{\g}-\g\Av_{x\in\toromi}\,\int_0^T ds H(s,\ep x)^2\\
&+ \int_0^T ds \,\underset{L^2(\mu_\ep)}{\sup spec}\bigl \{
\ep^{-1}\g\Av_{x\in\toromi}\, 2 H(s,\ep x)\Bigl[
\frac{m_{x+\frac{2b}{\ep}e,\frac{a}{\ep}}(s)-
 m_{x,\frac{a}{\ep}}(s)}{2b/\ep}\Bigr]+\ep^{d-2}\cL_\ep\bigr\}.
\end{split}
\end{equation}
It is convenient to introduce a free
scale parameter $n$, with
$n\uto\infty$ after $a\dto 0$ and $\ep\dto 0$, and write the gradient of masses
appearing in (\ref{corsica}) as
$\Av_{z,x}^\ast\t_z\bigl(\frac{\psi^{e}_{n,n}}{n}+\frac{\phi^{e}_{n,n}}{n}\bigr)$
(see (\ref{perleaiporci.K.V.}) and (\ref{estasi})).\\
By the definition of $\Av _{z,x}^\ast$, the subadditivity of $\sup spec$
and theorem \ref{pte},
\begin{equation*}
\limsup _{n\uto\infty,a\dto 0,\ep\dto
  0}\int_0^T ds\;\underset{L^2(\mu_\ep)}
{\sup spec}\bigl\{\ep^{-1}\g\Av_{x\in\toromi}\,2 H(s,\ep
x)\Av_{z,x}^\ast\t_z\frac{\phi^e_{n,n}}{n}+\ep^{d-2}\cL_\ep\bigr\}\leq 0.
\end{equation*}
Let us consider, for fixed $b,n,a$,
\begin{equation}
\label{floreancing}
 \sup spec _{L^2(\mu_\ep)}\bigl \{
\ep^{-1}\g\Av_{x\in\toromi}\,2 H(s,\ep
x)\Av_{z,x}^\ast\t_z\frac{\psi^e_{n,n}}{n}+\ep^{d-2}\cL_\ep\bigr\}.
\end{equation}
Thanks to the definition of the dual average $\Av_{z,x}^\star$ we can write
\begin{equation*}
\Av_{x\in\toromi}\, 2H(s,\ep
x)\Av_{z,x}^\ast\t_z\frac{\psi^e_{n,n}}{n}= \Av_{x\in\toromi}\,
a_x \t_x\frac{\psi^e_{n,n}}{n}
\end{equation*}
where
$a_x:=\Av_{z,x}^\star 2 H(s,\ep z)$.
Since $\Av_{z,x}^\star$ is translationally invariant w.r.t. $x$ and
$H$ is regular, we can proceed as at the very beginning of this
 section and safely replace $\t_x\frac{\psi^e_{n,n}}{n}$ by a local
 average $\Av_{|y-x|\leq \ell _1}
\t_y\frac{\psi^e_{n,n}}{n}$, $\ell \gg n$, to get
\begin{equation}
\label{quasimodo}
(\ref{floreancing}) \leq \sup spec _{L^2(\mu_\ep)}\bigl
\{\ep^{-1}\g\Av_{x\in\toromi}\, a_x\Av_{|y-x|\leq \ell _1}
\t_y\frac{\psi^e_{n,n}}{n} +\ep^{d-2}\cL_\ep\bigr\}+c(H)\g\ep \ell^2 .
\end{equation}
By the usual trick of localizing on boxes
$\L_{x,\ell}$ and
proposition \ref{jamme}, if $\ep$ is small enough then the first term in
the r.h.s. of (\ref{quasimodo}) is bounded from above by
\begin{equation*}
c\,\g ^2 \Av_{x\in\toromi}\,a_x^2 \sup _{\nu\in\cM(\L_{x,l})}
V_{x,\ell}\bigl(\frac{\psi^{e'}_{n,n}}{n};\nu\bigr)
\end{equation*}
which in turn, thanks to
(\ref{kolja}), is bounded from above by
$$
 c\, \g^2\Av_{x\in\toromi}\,a_x^2\leq c^\ast\,\g^2\Av_{x\in\toromi}\,H(s,\ep x)^2
$$
for some suitable positive  constant $c^\ast$.
Let us now choose $\g$ so small that $c^\ast\g^2-\g<0$.
Then, by the previous estimates, if $\ep$ is small enough,
\begin{equation}
\label{guatemala}
 \limsup _{n\uto\infty,a\dto 0,\ep\dto 0}\text{
r.h.s. of } (\ref{corsica})\leq
 \frac{\k}{\g}+(c^\ast\g^2-\g)
 \int_0^T\int_{\toroma}H(s,\theta)^2d\theta\,ds\leq\frac{\k}{\g}.
\end{equation}
 In order to conclude the proof it is enough to observe that
there exists a finite set  $\cH_b^\ast\sset\cH_b$
depending on $b$  such that
\begin{equation*}
\sup _{H\in\cH_b} H ^\star _{b,a,\ep}\leq 1 +\sup _{H\in
\cH_b^\ast}  H ^\star _{b,a,\ep}
\end{equation*}
so that
\begin{equation}
\label{dielangenacht}
\limsup_{n\uto\infty,a\dto 0,\ep\dto
0}\EE^{\mu^\ep}\bigl(\sup _{H\in
\cH_b} H^\star _{b,a,\ep}\bigr)\leq
1+\limsup_{n\uto\infty,a\dto 0,\ep\dto
0}\EE^{\mu^\ep}\bigl(\sup _{H\in
\cH^*_b} H^\star _{b,a,\ep}\bigr)
\leq 1 +\frac{\k}{\g}
\end{equation}
thus allowing to conclude the proof of (\ref{jerusalem}). 
\\ 
Let us now sketch the proof of (\ref{mirta}). Since
$C^1([0,T]\times\toroma)$ has a countable base, by Beppo--Levi theorem
it is enough to prove that there exists a constant $c_0$ such that,
given $H_1, \dots, H_n$ in $C^1([0,T]\times\toroma)$, then
\begin{equation}
\label{animemorte}
\int dQ(m)
\Bigl[\,\sup_{i=1,\dots,n}\int_0^T\int_{\toroma}\bigl(2\,m(s,\theta)
\frac{\partial}{\partial \theta_e} H_i(s,\theta)-H_i(s,\theta)^2\bigr)d\theta\,ds\,\Bigr]\leq c_0.
\end{equation}
By the  Lebesgue density theorem and the dominated convergence theorem, the l.h.s.
of $(\ref{animemorte})$ is equal to
$\lim_{a\dto 0}\EE_Q\bigl(\Theta^{(a)}\bigr)$ where, for any $\nu\in D([0,T],\cM_1)$,
\begin{equation}
\label{andreamoratto}
  \Theta^{(a)} (\nu):=\sup _{i=1,\dots,n}
  \int_0^T\int_{\toroma}\bigl(2\,\nu^{(a)}(s,\theta)
\frac{\partial}{\partial \theta_e}H_i(s,\theta)-H^2_i(s,\theta)\bigr)d\theta\,ds,
\end{equation}
with
$$
\nu^{(a)}(s,\theta):=\frac{1}{(2a)^d} \nu_s \bigl(
\{\theta'\in\toroma:\sup_{i=1,\dots ,d}|\theta_i'-\theta_i|\leq a
\}\bigr). 
$$
It is simple to prove (see  \cite{AF}, section 1.18) that 
\begin{equation*}
\begin{split}
& \lim_{a\dto 0}\int dQ(m)\bigl(\Theta ^{(a)}(m)\bigr)\leq
 \limsup_{a\dto 0,\ep\dto 0}\int dQ^{\a,\,\mu^\ep}(\nu)\bigl(\Theta ^{(a)}(\nu)
 \bigr)= \\
 \limsup_{b\dto 0,a\dto 0,\ep\dto 0}\EE^{\a,\mu^\ep}\Bigl(
&  \sup_{i=1,\dots,n}
\int_0^T\Av_{x\in\toromi}\bigl(\,2m_{x,\frac{a}{\ep}}(s,\ep x)
\bigl[\frac{H_i(s,\ep x +be)-H_i(s,\ep x)}{b}\bigr]-H_i^2(s,\ep x)\,\bigr)\,\Bigr).
\end{split}
\end{equation*}
By integrating by parts and observing that
$$
\sup _{H\in\cH_b} H^\star _{b,a,\ep}=\sup_{H\in C^1([0,T]\times\toroma)} H^\star
_{a,b,\ep},
$$
the thesis follows from (\ref{dielangenacht}). 

\subsection{Hydrodynamic limit without regularity of the diffusion matrix.}
\label{generale}
In this last paragraph we shortly discuss the hydrodynamic limit when
the regularity condition on the diffusion matrix is replaced by the two
conditions at the end of theorem \ref{HL}, in the sequel referred to as
assumptions $A(\rho)$. The main idea here is to
prove that one can safely introduce a density cutoff near the edges of
the interval $(0,1)$, and for this purpose the main technical tool is the
following result.
\begin{lemma}
\label{fabri11}
Assume that the sequence of initial probability measures $\mu^\ep$
satisfy $A(\rho)$. Then there exists a constant $0<\bar \rho \leq \rho$ such that,
for any $T>0$ and any disorder configuration $\a$,
\begin{equation}
\label{fabri10}
\lim_{\ell\uto\infty,\ep\dto 0}\EE^{\mu^\ep}\Bigl(\,\int_0^T ds\,
\Av_{x\in\toromi}
\,\bigl(\id _{\{m_{x,\ell}(s)<\bar\rho\}}+\id _{\{m_{x,\ell}(s)>1-\bar\rho\}}\bigr)\,\Bigr)=0.
\end{equation}
\end{lemma}  
\begin{proof}
  For simplicity, we consider in (\ref{fabri10}) only the contribution
  coming from $\id_{\{m_{x,\ell}(s)<\bar\rho\}}$, the other one being
  similar.  Given two probability measures $\mu_1,\,\mu_2$ on $\O_\ep$,
  we will write $\mu_1\leq \mu_2$ if $\mu_1(f) \le \mu_2(f)$ for any
  function $f$ which is increasing w.r.t. the partial order in $\O_\ep$
  given by $\h\le \h' \Leftrightarrow \h(x)\le \h'(x) \; \forall \,x\in
  \toromi$.  It is then simple to check that our model is
  \emph{attractive} \cite{Li} in the sense that $\mu_1 \le \mu_2$
  implies that $\PP_t^{\mu_1}\leq\PP_t^{\mu_2}$ for any $t>0$ and for
  any disorder configuration $\a$. Therefore, condition $A(\rho)$
  implies that there exists $\l<0$ such that $\mu^\l_\ep\leq\mu^\ep$ for
  any $\ep$. Let now $\bar \rho:=\ov 2 \min\bigl
  (\frac{e^{\l-B}}{1+e^{\l-B}},\rho\bigr)$. Then, given $\b>0$ and
  thanks to attractivity, the entropy inequality (\ref{ent_ineq}) and
  the identity
  $H[\mu^\ep|\mu^\ep_\ast]=H[\PP^{\mu^\ep}|\PP^{\mu^\ep_\ast}]$,
\begin{equation}
\label{elance}
\begin{split}
  & \EE^{\mu^\ep}\bigl(\,\int_0^T ds\, \Av_{x\in\toromi}\id
  _{m_{x,\ell}(s)<\bar\rho}
  \,\bigr)\leq \\
  \frac{1}{\b}H\bigl[\mu^\ep|\mu^\ep_\ast\bigr]&+\frac{1}{\b}\ln
  \Bigl(\,\mu_\ep^\l \bigl(\exp \bigl\{ \int_0^T
  ds\,\b\,\Av_{x\in\toromi}\id _{m_{x,\ell}(s)<\bar\rho}\bigr\}
  \bigr)\,\Bigr).
\end{split}
\end{equation}
Thanks to the Jensen's inequality and the reversibility of $\cL_\ep$
w.r.t.  $\mu^\l_\ep$ the second addendum in the r.h.s. of (\ref{elance})
can be bounded by 
\begin{equation}
\label{babbonatale}
\frac{1}{\b}\ln\bigl(\,\mu_\ep^\l \bigl(\exp
\{T\,\b\,\Av_{x\in\toromi}\id _{m_{x,\ell}<\bar\rho}\}\bigr)\,\bigr).
\end{equation}
Let us call $\nu^\l$ the  product measure on $\O_\ep$  such that
$\nu^\l(\eta_x)=\frac{e^{\l-B}}{1+e^{\l-B}}$. Then, it is simple to check that $
\nu^\l \leq\mu_\ep^\l$, which implies that
$$
(\ref{babbonatale}) \leq \frac{1}{\b}\ln\bigl(\,\nu^\l \bigl(\exp
\{T\,\b\,\Av_{x\in\toromi}\id _{m_{x,\ell}<\bar\rho}\}\bigr)\,\bigr).
$$
At this point, let us recall a general result based on the Herbst's
argument and the logarithmic Sobolev inequality (see \cite{CCirM} for a
complete discussion): for any $\g>0$ and any function $f$ on $\O_\ep$
$$
\nu^\l(e^{\g f})\leq e^{ c_f\, \gamma^2 + \g\nu^\l(f)}$$
where $c_f:= c\,\sum_{x\in\toromi}\|\nabla_x f\|^2_\infty$ and $c=c(B,\l)$ is a
suitable constant independent of $\ep$ ($c$ is related to the
logarithmic Sobolev
constant of the Bernoulli measure $\nu^\l$). \\
Thus
\begin{equation}
\label{stadtmitte1}
\frac{1}{\b}\ln\bigl(\,\nu^\l \bigl(\exp
\{T\,\b\,\Av_{x\in\toromi}\id _{m_{x,\ell}<\bar\rho}\}\bigr)\,\bigr) \leq
c\,T^2\,\b\,\ep^d\ell^d + T\,\Av_{x\in\toromi}\nu^\l\bigl
(\id _{m_{x,\ell} <\bar\rho}\bigr).
\end{equation}
 Since  $\bar \rho <\nu^\l(\eta_0)$,  by choosing $\b^2:=
H[\mu^\ep|\mu^\ep_\ast]/(T^2\,\ep^d\ell^d)$ the  r.h.s. of
 (\ref{stadtmitte1}) is negligible as $\ell\uto\infty,\ep\dto 0$.
Since $H[\mu^\ep|\mu^\ep_\ast]=o(\ep^{-d})$, the thesis
follows by collecting all the above estimates.
\end{proof}
Using the above result we are in position to discuss our density
cutoof. Let us recall first that, given a generic continuous extensions $\bar D$ of
$D$ outside the interval $[\rho,1-\rho]$, any weak solution $m(t,\theta)$ of the
Cauchy problem (\ref{eden}), where $D$ has been replaced by $\bar D$ and
$\rho\leq  m_0 (\theta) \leq 1-\rho$ for any $\theta \in\toroma$,
satisfies $\rho\leq m(t,\theta) \leq 1-\rho$ for any $0\le t\le T$ and
any $\theta\in \toroma$. Let $\bar D$ be defined as
\begin{equation*}
\bar D(m) :=
\begin{cases}
D(\bar \rho)&\text{ if } 0\leq m \leq \bar \rho\\
D(m) &\text{ if }  \bar \rho\leq m\leq1-\bar \rho\\
D(1-\bar\rho) &\text{ if } 1-\bar\rho\leq m\leq1.
\end{cases}
\end{equation*}
Let us explain next how one should modify the proof of theorem
\ref{zeus} in order to get the same result but with $D$ replaced by
$\bar D$ in the definition of $\bar H_{b,a,\ep}$ (in what follows this
replacement will be understood without further notice). To this aim it
is convenient to introduce the following shorter notation
$$
\chi_{x,\ell}:= \id _{m_{x,\ell}<\bar\rho}+\id _{m_{x,\ell}>1-\bar\rho}.
$$
Then, thanks to lemma \ref{fabri11}, equation (\ref{aiace})
 can be substituted by
\begin{equation*}
\begin{split}
\inf_{g\in\GG}\inf_{r\geq 0 } &\limsup_{b\dto 0,a\dto 0,l\uto \infty,\ep\dto
0}\EE^{\mu^\ep}\Bigl(\,\bigl|\int_{0}
^{T}{\ep}^{-1}\Av_{x\in\toromi}\nabla_e^\ep H(s,\ep x)\Bigr[\Av_{y:|y-x|\leq \ell
  _1}(j_{y,y+e}+\t_y\cL g)\\
&+ \sum_{e'\in\cE}\bar{D}_{e,e'}(m_{x,\frac{a}{\ep}})\bigl[
\frac{m_{x+\frac{b}{\ep}e',\frac{a}{\ep}}
-m_{x-\frac{b}{\ep}e',\frac{a}{\ep}}}{2b/\ep}\bigr] \Bigr]\,-r\int_0^T ds\,\Av_{x\in\toromi} \chi_{x,\ell}\,ds
\,\bigr|\Bigr)=0
\end{split}
\end{equation*}
and the main issue is to prove proposition \ref{wannsee} with $\bar
J^g_{b,a,\ell,\ep}$ replaced by
$$
\bar J^{g,r}_{b,a,\ell,\ep}:=\Bigl\{\text{r.h.s. of }(\ref{eq:kalinka})\text{
  with $D\to \bar D$}\,\Bigr\} -\ep
\,r\int_0^T ds\,\Av_{x\in\toromi} \chi_{x,\ell}\,.
$$
In turn the proof of the modified version of proposition \ref{wannsee}
is splitted into several steps, one for each term $\O^{(r)}_j$, $j=0,1,\dots
5$, where now
$$
\O^{(r)}_0:=\sup spec _{L^2(\mu_\ep)}\bigl\{\ep^{-1}\Av_x J (\ep x)\psi^{(0)}_x-r
\Av_{x\in\toromi} \chi_{x,\ell}+ \g\ep^{d-2}\cL\bigr\}.
$$
and all the other $\O_j$ are unchanged.
It thus remains to explain how the discussion in subsection
\ref{pietro} has to be modified in order to apply to
$\O_0^{(r)}$.
Because of the new definition
of $\O_0$, (\ref{cocktail}) has to be replaced by
\begin{equation}
\label{shipylov}
\begin{split}
& \inf_{g\in\GG}\;\inf_{r\geq 0}\limsup _{n\uto\infty,\ell\uto\infty,\ep\dto 0}\;
 \Av_{x\in\toromi}\sup_{|\b|\leq T}\,
\sup_{m\in[0,1]}\\
&\Bigl[
\underset{L^2(\nu_{\L_{x,\ell},m})}{\sup spec} \bigl\{
\ep^{-1}\b\,\Av_{y:|y-x|\leq\ell_1}\t_y\psi_m^{(n,g)} +\ep^{d-2}\cL_\ep\bigr\}
-r\id _{m<\bar\rho}-r\id _{m>1-\bar\rho}
\Bigr]\leq 0
\end{split}
\end{equation}
where $D\rightarrow \bar D$ in the definition of $\psi_m^{(n,g)}$.\\
We observe that, provided $\ep\ell^{d+2}\ll 1$, the $\sup spec$
inside the square bracket in (\ref{shipylov}) is bounded by $c_g \,
T^2$, for a suitable constant $c_g$ depending on $g$. That follows
immediately from perturbation theory (see proposition \ref{jamme}) and
the estimate (\ref{kolja}). Therefore, by choosing $r$ large enough, we only need to
prove (\ref{shipylov}) with $m\in [\bar\rho, 1-\bar\rho]$ where $D(m)$
and $\bar D(m)$ coincide. Similarly one shows that the two ``$\sup_{m\in
  [0,1]}$'' appearing in (\ref{lowen}) and (\ref{mormoni}) can be safley
replaced by ``$\sup_{m\in [\bar\rho, 1-\bar\rho]}$''.


\section{Disorder induced fluctuations in the averaged gradient density field}
\label{Fluct}
In this section we analyze a key term that, as we seen in section
\ref{pf_zeus}, arises naturally when one tries to approximate spatial
averages of the current with spatial averages of \emph{gradients of the
density profile}. Since the currents $j_{x,x+e}$ have, by construction,
zero canonical expectation with respect to \emph{any canonical measure}
on any set $\L \ni x,x+e$, in order to approximate $\Av_x \,j_{x,x+e}$
with suitable averages of gradients of the density field, one is forced
to subtract from these gradients appropriate canonical expectations.
Therefore, a key point in order to establish the hydrodynamical limit,
is to prove that these ``counter terms'' vanish as $\eps \downarrow 0$.
These kind of terms arise also in the hydrodynamical limit of
non--disordered lattice gases (see \cite{VY}, section $7$) with short
range interaction. In our context however their nature is quite
different and, as we will show next, they are basically produced by
fluctuations in the disorder field.

\smallno In order to be more precise recall first, for any given $e\in
\cE$, the notation $\L_n^{1,e}, \L_n^{2,e}$ and $\L_n^e :=
\L_n^{1,e}\cup \L_n^{2,e}$ described in section \ref{notation}, together
with the associated densities $m_n^{1,e}:= m_{\L_n^{1,e}},\,m_n^{2,e}:=
m_{\L_n^{2,e}},\, m_n^{e}:=m_{\L_n^e}$.

\smallno
Using the above notation and given two integers $n\leq s$ and a vector
$e\in\cE$,  the basic
object of our investigation is defined as (see (\ref{estasi})):
\begin{equation}
\label{phi(n,s)}
  \phi_{n,s} := \mu[ m_n^{2,e} - m_n^{1,e} \tc
  m^e_s ]
\end{equation}
Notice that if the disorder configuration $\a$ was identical in the two
cubes $\L_n^{1,e}$ and $\L_n^{2,e}$ then $\phi_{n,n}$ would be
identically equal to zero. Moreover $\bbE\bigl(\phi_{n,s} \bigr)=0$
and $\bbE\bigl([\phi_{n,s}]^2 \bigr)= O(n^{-d})$ uniformly in $s\ge
n$.
\begin{remark}
The fact that $\phi_{n,s}$ is small (on some average sense) with $n$ and
not with $s$ is one of the main differences with non disordered lattice gases
where, instead, the analogous term goes very fast to zero as $s\uparrow
\infty$ (see \cite{VY}, section $10$).
\end{remark}
The main result of this section is the proof that the contribution to
the hydrodynamical limit of suitable spatial averages of
$\frac{\phi_{n,n}}{n}$ is negligible as $\eps \downarrow 0$ at
least in dimension $d\ge 3$.

\medno In order to be more precise let us introduce the following
equivalence relation.
\begin{definition}
\label{equiv}
Given two families of functions $f_{x,n,a,\eps}(\a,\h)$ and
$g_{x,n,a,\eps}(\a,\h)$ with $x\in \toromi, n\in \bbN, a>0, \eps >0$ we
will write $ f_x \approx g_x$ if, for any given $T>0$ and for almost all
disorder configurations $\a$,
\begin{equation*}
 \limsup_{n \uparrow \infty , a \downarrow 0 ,\ep \downarrow 0}
  \Av_{x \in \toromi } \sup_{ |\b | \leq  T}
\sup_\nu \; \sup spec_{L^2(\nu)}\bigl \{ \ep^{-1} \b (f_x - g_x)  +
\ep^{-2} \Av_{b\in\L_{x,2a/\eps}}
\cL_b
 \bigr \} \leq 0
\end{equation*}
where $\sup_\nu$ is the supremum over $\nu$ in the set
$\cM(\L_{x,2\frac{a}{\ep}})$ of all the canonical measures on
$\L_{x,2\frac{a}{\ep}}$.
\end{definition}
We are now in a position to state our main result.  Assume that a given
direction $e$ has been fixed once and for all and, given two integers
$\ell \le s$ with $\frac{s}{\ell}\in \bbN$ and $x\in \toromi$, recall
the definition of the spatial average $\Av_{z,x}^{\ell,s}$ given in (\ref{dinoce}).
\begin{theorem} For any $d\ge 3$
\label{pte}
\begin{equation*}
\Av_{z,x}^{n,\frac{a}{\ep}} \t_z \frac{\phi_{n,n}}{n} \approx 0
\,.\end{equation*}
\end{theorem}
Before discussing the plan of the proof of the theorem we would like to
justify the restriction $d\ge 3$. If we pretend that the particle
density is constant everywhere, say equal to $m$, then
\begin{gather*}
  \sup spec_{L^2(\nu)}\bigl\{\ep^{-1}\b\Av_{z,x}^{n,\frac{a}{\ep}}\t_z
  \frac{\phi_{n,n}}{n}
  + \ep^{-2} \Av_{b\in\L_{x,2a/\eps}}\cL_b \bigr\} \le \\
  \ep^{-1}\b\Av_{z,x}^{n,\frac{a}{\ep}}\t_z \frac{\phi_{n,n}(m)}{n}
\end{gather*}
Since the typical fluctuations (in $\a$) of the quantity
\begin{equation*}
  \eps^{-1} \Av_{z,x}^{n,\frac{a}{\ep}}\t_z \frac{\phi_{n,n}(m)}{n}
\end{equation*}
are of the order of $\eps^{\frac{d-2}{2}}C(a,n)$, necessarily we must
assume $d \ge 3$ since $\eps\downarrow 0$ before $a\downarrow 0$ and
$n\uparrow \infty$.
\\

\subsection{Plan of the proof of theorem \ref{pte}.}
\label{road}
The main difficulty in proving theorem \ref{pte} lyes in the fact that
first $\eps \dto 0$ and only afterward $n\uto\infty$. In particular
there is no hope to beat the diverging factor $\eps^{-1}$ appearing in
definition \ref{equiv} with the typical smallness
$O(n^{-\frac{d+2}{2}})$ of $\frac{\phi_{n,n}}{n}$. The main idea is
therefore first to try to prove that
\begin{equation}
  \label{eq:n->s}
\Av_{z,x}^{n,\frac{a}{\ep}}\t_z\frac{\phi_{n,n}}{n}  \approx
\Av_{z,x}^{s,\frac{a}{\ep}}\t_z\frac{\phi_{s,s}}{s}
\end{equation}
where the new mesoscopic scale $s=s(\eps)$ diverges sufficiently fast
as $\eps \dto 0$. By standard large deviations estimates (see lemma
\ref{norma_uniforme}) it's simple to
verify that, given $ 0<\d \ll 1$ and $0<\g<1$, for
almost any disorder configuration $\a$ and $s=O(\ep^{-\g})$
\begin{equation}
\label{eq:phis-bound}
\sup_{x\in\toromi} |\t_x \,\phi_{s ,s}|\leq C s^{ -\frac{d}{2}+\d}
\end{equation}
for any $\eps$ small enough.
In particular, by a  trivial  $L^\infty$ estimate,
\begin{equation}
\label{trivial}
\Av_{z,x}^{s,\frac{a}{\ep}}\t_z\frac{\phi_{s,s}}{s} \approx 0 \quad \text{ if
} \quad \g > \frac{2}{d+2}.
\end{equation}
The above simple reasoning suggests to define a first mesoscopic critical scale
$s_\infty:=\ep^{-\frac{2}{d+2}}$ above which things become trivial. It is
important to outline that we will \emph{not} be able to prove
(\ref{eq:n->s}) with $s \gg s_\infty$ but only with $s=\bar s$ where
$\bar s := \ep^\d s_\infty$ and $0<\d
\ll 1$ can be taken arbitrarily small. \\
Once we have reached scale $\bar s$ we cannot simply use $L^\infty$ bounds but we need
to appeal to an improved version of the well known Two Blocks Estimate
(see proposition \ref{super-two-blocks} below) in order to conclude that
$\Av_{z,x}^{\bar s,\frac{a}{\ep}}\t_z \frac{\phi_{\bar s,\bar s}}{\bar s}
 \approx 0$.

\smallno
We now explain the main steps in the proof of (\ref{eq:n->s})
with $s =\bar s$.  As discussed in subsection \ref{lete}, a main tool for
estimating eigenvalues is given by localization together with
 perturbation theory. However, because
of proposition \ref{jamme}, it turns out that this technique can be
applied to prove (\ref{eq:n->s})  only if
  \begin{equation*}
\ep s^{d+2}\|\Av_{z,x}^{n,\frac{a}{\ep}}\t_z\frac{\phi_{n,n} }{n}
-
\Av_{z,x}^{s,\frac{a}{\ep}}\t_z\frac{\phi_{s,s}}{s}\|_\infty\leq {\rm const},
  \end{equation*}
  that is if $\ep s^{d+2} \le {\rm const}$. In particular we see
  immediately that this approach cannot be used directly
  to prove (\ref{eq:n->s}) for $s=\bar s$, but only up to a new critical
  mesoscopic scale $s_0:= \ep^{-\frac{1}{d+2}}$.
  \\
  Assuming that we have been able to replace
  $\Av_{z,x}^{n,\frac{a}{\ep}}\t_z\frac{\phi_{n,n} }{n}$ with
  $\Av_{z,x}^{s_0,\frac{a}{\ep}}\t_z\frac{\phi_{s_0,s_0}}{s_0}$, we face
  the problem to increase the mesoscopic scale from $s_0$ to $\bar s$.
  \\
  The main observation now is that the $L^\infty$ norm of the new
  quantity $\Av_{z,x}^{s_0,\frac{a}{\ep}}\t_z\frac{\phi_{s_0,s_0}}{s_0}$
  is at least smaller than $s_0^{-\frac{d+2}{2}}$ (see
  (\ref{eq:phis-bound})) almost surely (here and in what follows we
  deliberately neglect the correction $s^\d$ appearing in
  (\ref{eq:phis-bound})). This means that the limit scale beyond which
  perturbation theory cannot be applied, previously equal to $s_0$, is
  now pushed up to a new scale $s_1$ given by
\begin{equation*}
\ep s_1^{d+2}\,s_0^{-\frac{d+2}{2}} \le {\rm
  const} \; \Rightarrow \; s_1 = \ep^{-\frac{3}{2(d+2)}}
\end{equation*}
The above remark clearly suggests an inductive scheme on a sequence of
length scales $\{s_k\}_{k\ge 0}$ given by
\begin{equation*}
s_0:= \ep^{-\frac{1}{d+2}}\,; \qquad s_{k+1} := \ep^{-\frac{1}{d+2}} \sqrt{s_k}
\end{equation*}
in which one proves recursively, by means of localization on
scale $s_{k+1}$ combined together with perturbation theory, that
\begin{equation*}
  \Av_{z,x}^{s_k,\frac{a}{\ep}}\t_z\frac{\phi_{s_k,s_k}}{s_k}
- \Av_{z,x}^{s_{k+1},\frac{a}{\ep}}\t_z\frac{\phi_{s_{k+1},s_{k+1}}}{s_{k+1}}
\approx 0.
\end{equation*}
Notice that $\lim_{k\to \infty} s_k = s_\infty$ where $s_\infty =
\ep^{-\frac{2}{d+2}}$ represents the limiting scale introduced at the
beginning of this section. \\
A large but finite number of steps of the inductive scheme
proves that
\begin{equation*}
  \Av_{z,x}^{n,\frac{a}{\ep}}\t_z\frac{\phi_{n,n} }{n}
- \Av_{z,x}^{\bar s,\frac{a}{\ep}}\t_z\frac{\phi_{\bar s,\bar s}}{\bar s}
\approx 0
\end{equation*}
where, as before, $\bar s=\ep^\d s_\infty$. We remark that for this part
of the proof we only need $d\ge 2$, while we will assume $d\geq 3$  when proving the improved version of the Two Blocks estimate (see proposition \ref{super-two-blocks}).

\subsection{Preliminary tools.}
In this section we collect some general techniques that are common to
all the steps of the proof of theorem \ref{pte}. We recall that  $\L^e_{z,\ell}$ denotes the translated by $z$ of the box $\L^e_{\ell}$.

\smallno
\begin{lemma}
\label{truc}
Let $\ell_0<\ell_1 <\ell_2$ be odd integers such that $\frac{\ell_2}{\ell_0} \in
\bbN$. Let $\nu$ be an arbitrary canonical measure on the cube
$\L_{\ell_2}$ and let $f$ be a function with
support in $\L^e_{\ell_1}$. Then
\begin{align*}
\label{truc1}
\sup spec_{L^2(\nu)}\{\Av_{z,0}
^{\ell_0,\ell_2} \t_z f &+ \Av_{b \in \L_{\ell_2}}\cL_b \} \leq \\
&\Av_{z,0}^{\ell_0,\ell_2} \sup_{\nu'} \, \sup spec_{L^2(\nu')} \{\t_z
f + c \Av_{b \in \L^e_{z,\ell_1}} \, \cL_b \}
\end{align*}
where $\nu'$ varies in $\mstor_\a(\L^e_{z,\ell_1})$ and $c$ is a
suitable constant.
\end{lemma}
\begin{proof}
It is sufficient to observe that
\begin{equation*}
  \Av_{b \in \L_{\ell_2}}\cL_b  \le c\,
  \Av_{z,0}^{\ell_0,\ell_2}\bigl(\Av_{b \in \L^e_{z,\ell_1}} \, \cL_b\bigr)
\end{equation*}
and localize in the box $\L^e_{z,\ell_1}$.
\end{proof}
At this point, it is convenient to observe the  factorization property
of the average $\Av_{z,x}^{\ell, s} $ defined in (\ref{dinoce}): given
odd integers  $\ell, \ell ', L$ such that $\frac{\ell'}{\ell},
\frac{L}{\ell'}\in\NN$, then
\begin{equation}
\label{bobo}
\Av_{z,x}^{\ell,L}f_z = \Av_{z,x}^{\ell',L}\bigl(
\Av_{w,z}^{\ell,\ell'}f_w \bigr).
\end{equation}

\begin{prop}
\label{tolmezzo}
Let $d \geq 2$, $0<\g\leq\g'<1$ and $\g'<\frac{1}{d+2}+\frac{\g}{2}$.
 If either $\ell=n $ and $s=O\bigl(\ep^{-\frac{1}{d+2} }\bigr)$ or
$\ell = O(\ep^{-\g})$ and $s=O(\ep^{-\g'})$, then
$$
\Av_{z,x}^{\ell, \frac{a}{\ep}}
\t_z \frac{\phi_{\ell,s}}{\ell} \approx \Av_{z,x}^{s,\frac{a}{\ep}}
\t_z \frac{\phi_{s,2s}}{s}.
$$
\end{prop}
\begin{proof}
  By the factorization property (\ref{bobo}), we have
\begin{equation*}
\Av_{z,x}^{\ell, \frac{a}{\ep}}
\t_z \frac{\phi_{\ell,s}}{\ell} -\Av_{z,x}^{s, \frac{a}{\ep}}
\t_z \frac{\phi_{s,2s}}{s} =
 \Av_{z,x}^{s, \frac{a}{\ep}} \Bigl[
 \Av_{w,z}^{\ell,s}\t_w \frac{\phi_{\ell,s}}{\ell} -
\t_z \frac{\phi_{s,2s}}{s}.
\Bigr]
\end{equation*}
Therefore, by lemma \ref{truc}, it is enough to prove that for any $T>0$
and for almost all disorder configuration $\a$
\begin{equation}
\label{tolm2}
 \limsup_{n \uparrow \infty , a \downarrow 0 ,\ep \downarrow 0 }
 \Av_{x \in \toromi} \sup_{|\b| \leq T} \sup_{\;\nu \in
 \mstor_\a (\L^e_{x,2s})} f_{x,\nu} \leq0
\end{equation}
where
$$
f_{x,\nu} := \sup spec_{L^2(\nu)} \Bigl\{ \ep^{-1} \b \Bigl[
\Av_{z,x}^{\ell,s} \t_z \frac{\phi_{\ell,s}}{\ell} - \t_x\frac{ \phi_{s,2s}}{s}
\Bigr] +c \ep^{-2} \Av_{b \in \L_{x,2s}^e} \cL_b
\Bigr\}
$$
for a suitable constant $c$. Notice that $\t_x\frac{ \phi_{s,2s}}{s}
= \nu\bigl( \Av_{z,x}^{\ell,s} \t_z \frac{\phi_{\ell,s}}{\ell}
\bigr)$ $\nu$ a.s..

\\
Because of lemma \ref{norma_uniforme}, given $0< \d \ll 1$, for almost
all $\a$ and $\ep$ small enough
\begin{equation*}
 \sup_{x\in\toromi} \bigl\|\t_x\frac{\phi_{\ell,s}}{\ell}\bigr\|_\infty \le
  \begin{cases}
    \ell^{-1} & \text{ if } \ell=n \\
     \ell^{-(d+2)/2 +\d} & \text{ if } \ell = O(\ep^{-\g})
  \end{cases}
\end{equation*}
Thanks to the above bound and to the choice
$\g'<\frac{1}{d+2}+\frac{\g}{2}$, for almost all
$\a$ and $\ep$ small enough, we can apply proposition \ref{jamme} together with lemma
\ref{cleonimo} to get
\begin{equation}
\label{sznitman}
\sup_{\nu \in\mstor_\a(\L^e_{x,2s})}f_{x,\nu} \le c\, T^2 \ell^{-2}s^{d+2}\sup_m F(x,m)
\end{equation}
where $m$ varies among all possible particle densities in $\L^e_{x,2s}$
 and
 \begin{equation*}
 F(x,m) := \Var_{\mu^{
 \l_x(m)}}\bigl(\Av^{\ell,s}_{z,x}\,\t_z \phi_{\ell,s} \bigr)
 \end{equation*}
and $\l_x(m):= \l_{\L^e_{x,2s}}(m)$.

\smallno
We claim that for almost all $\a$ and $\ep$ small enough
\begin{equation}
  \sup_{x\in \toromi}\sup_m F(x,m) \le c\, s^{-2d+2\d}
 \label{F}
\end{equation}
thus proving the proposition since $d\geq2$. The proof of (\ref{F})
follows exactly the same lines of the proof of proposition
\ref{bound:variance} with the main difference that it is necessery to
use lemma \ref{stima} in order to control the empirical chemical
potetials (see also section 4.7 in \cite{AF})
\\
\end{proof}

\begin{prop}
\label{paularo}
Let $d\geq2$, $\frac{1}{d+2}\leq  \g\leq \g' <1$ and
$\g'<\frac{1}{d+2}+\frac{\g}{2}$. Set  $s=O(\ep^{-\g})$ and
$s'=O(\ep^{-\g'})$. Then
$$  \Av_{z,x}^{s, \frac{a}{\ep}} \t_z  \frac{\phi
_{s,s}}{s} \approx \Av_{z,x}^{s, \frac{a}{\ep}} \t_z
 \frac{\phi_{s,s'}}{s} . $$
\end{prop}

\begin{proof}
\no By lemma $\ref{truc}$ it is enough to prove that for any $T>0$ and
for almost any disorder configuration $\a$
\begin{equation*}
 \limsup_{a\downarrow0,\ep\downarrow0}\Av_{x\in\toromi}\sup_{|\beta
 |\leq T}\sup_{\nu \in \cM_\a(\L^e_{x,s'})}f_ x\leq 0
\end{equation*}
where
\begin{equation*}
 f_x:=\sup spec_{L^2(\nu)}\Bigl\{\ep^{-1}\beta\t_x\Bigl[\frac{\phi_{s,s}}{s}
-\frac{\phi_{s,s'}}{s}\Bigr]+c \ep^{-2}\Av_{b\in\L_{x,s'}^e}\cL_b\Bigr\}
\end{equation*}
for a suitable constant $c$.
Notice that $\nu(\phi_{s,s})=\phi_{s,s'}$ $\nu$ a.s..

\\
Because of lemma \ref{norma_uniforme}, given $0<\d\ll 1$, for almost
all $\a$ and $\ep$ enough small
\begin{equation*}
\sup_{x\in\toromi}\bigl\|\t_x \frac{\phi_{s,s}}{s}\bigr\|_\infty\leq s^{-(d+2)/2 +\d}.
\end{equation*}
Thanks to the above bound and to the choice $\g'<\frac{1}{d+2}+\frac{\g}{2}$, for
almost any $\a$ and $\ep$ small enough, we can apply proposition
\ref{jamme} together with lemma \ref{cleonimo} to get
\begin{equation}
\sup_{\nu \in\cM_\a(\L^e_{x,s'})}f_x \le T^2 s^{-2}(s')^{d+2}\sup_m F(x,m),
\end{equation}
where $m$ varies among all possible particle densities in $\L^e_{x,s'}$,
 \begin{equation*}
 F(x,m) := \Var_{\mu^{ \l_x(m)}}\bigl(\t_x\phi_{s,s})
 \end{equation*}
and now $\l_x(m)=\l_{\L^e_{x,s'}}(m)$.\\
We claim  that for almost all $\a$ and $\ep$ small enough
\begin{equation}
  \sup_{x\in \toromi}\sup_m F(x,m) \le c\, s^{-2d+2\d}
\label{FF}
\end{equation}
thus proving the proposition because of the constraint on $\g,\g',d$.
The proof of (\ref{F}), requiring $d\geq 2$, follows exactly the same lines of the proof of proposition \ref{bound:variance} with the main difference that  it is necessery to use lemma  \ref{stima} in order to control the  empirical chemical potetials (see also section
4.6 in \cite{AF})
\end{proof}

\subsection{From scale $n$ to scale $s_0$.} 
Here we show how to replace the starting scale $n$ with our first
mesoscopic scale increasing with $\ep$, $s_0 =O\bigl(
\ep^{-\frac{1}{d+2}}\bigr)$.

\begin{prop}\label{pae12}
Let $d \geq 3$. Then
\begin{equation}
 \label{freud}
\Av_{z,x}^{n,\frac{a}{\ep}}\t_z\Bigl[\frac{\phi_{n,n}}{n}-\frac{\phi_{n,s_0}}{n}\Bigr]\approx 0.
\end{equation}
\end{prop}
\begin{proof}
Without loss of generality, we assume that $\frac{s_0}{n}\in\NN$ and similarly for $\frac{a/\ep}{n}$. \\
By the definition of $\Av_{z,x}^{n,\frac{a}{\ep}}$ and setting
$B=Q_{a/\ep}\cap n\,\ZZ^d$, in order to prove (\ref{freud}) it is
enough to show that
\begin{equation}
\label{eq:forster}
\limsup_{n\uto\infty, a \dto 0,\ep\dto 0}\Av_{x\in\toromi}\sup_{|\b|\leq
T}\sup_\nu\,\sup spec_{L^2(\nu)}\bigl\{\ep^{-1}\b\Av_{z\in B+x}\t_z
f_{n,s_0}+c\,\ep^{-2}\Av_{b\in\L_{x,3a/2\ep}}\cL_b\bigr\}\leq0
\end{equation}
where $f_{n,s_0}:=\frac{\phi_{n,n}}{n}-\frac{\phi_{n,s_0}}{n}$ and $\nu$ varies in $\cM ( \L_{x,3a/2\ep})$. The
proof is nothing more than a careful writing of the spatial average $\Av
_{z,x}^{n,\frac{a}{\ep}}$ together with the subadditivity property of
$\sup spec$.
\\
Setting  $B'=Q_{s_0}\cap n\,\ZZ^d $, $Y=Q_{a/\ep}\cap s_0\ZZ^d$ we can write
$B = \cup_{y\in Y} (B'+y)$ so that
$$
 \Av_{z\in B+x}\t_z f_{n,s_0} = \Av_{y\in Y+x}\Av_{z\in B'+y}\t_z
 f_{n,s_0}
$$
and
$$
\Av_{b\in\L_{x,3a/2\ep}}\cL_b \leq c\,\Av_{y\in
Y+x}\Av_{b\in\L_{y,2s_0}}\cL_b \,.
$$
By the subadditivity property of $\sup spec$, (\ref{eq:forster}) is
bounded from above by
\begin{equation*}
\limsup_{n \uparrow \infty, a \downarrow 0 , \ep \downarrow 0 }
\underset{x \in \toromi}{\Av }
\sup_{ |\b | \leq T} {\sup} _{\nu   } \sup spec_{ L^2 (\nu)}
\{ \ep^{-1}\b \Av^{(n)}_{y\in Q_{s_0}}\t_{x+y}f_{n,s_0}+c\,
\ep^{-2}\Av_{b\in\L_{x,2s_0}}\cL_b\}
\end{equation*}
where $\nu$ varies among $\cM(\L_{x,2 s_0})$
and $\Av^{(n)}_{y\in \L}:=\Av_{y\in \L\cap\, n\ZZ^d}$.

\smallno
At this point we can apply perturbation theory (see proposition
\ref{jamme}): since
$\lim_{n\uto \infty}\sup_{\eps>0} \ep s_0^{d+2}\|f_{n,s_0}\|_\infty = 0$, it is enough to prove that for
almost all disorder $\a$
\begin{equation}
  \label{lim-Psi}
\limsup_{n \uparrow \infty  , \ep \downarrow 0 }\frac{1}{n^2} \Av_{x \in
\toromi} \, {\sup}_{\nu \in \cM( \L_{x, 2s_0 } )  }
\Psi^{(\nu)}_{s_0}\bigl(\t_x\bigl[\phi_{n,n} - \phi_{n, s_0} \bigr] \bigr)=0
\end{equation}
where
$$
\Psi_{s_0}^{(\nu)}(f)=s_0^d\,\nu\bigl(\Av^{(n)}_{y\in Q_{s_0} }\t_y f,(-
\cL_{\L_{2s_0}})^{-1}\Av^{(n)}_{y\in Q_{s_0}}\t_y f \bigr).
$$
In order to prove (\ref{lim-Psi}) it is clearly sufficient to prove
it with $\phi_{n,n}$ replaced by $\phi_{n,n^4}$, provided one is able to
show that for almost any disorder $\a$
\begin{equation}
\label{n->n^4}
\limsup_{n \uparrow \infty  , \ep \downarrow 0 }\frac{1}{n^2} \Av_{x \in
\toromi} \, {\sup}_{\nu \in \cM( \L_{x, 2s_0 } )  }
\Psi^{(\nu)}_{s_0}\bigl(\t_x\bigl[\phi_{n,n} - \phi_{n,n^4} \bigr] \bigr)=0.
\end{equation}
We will concentrate only on the first step and refer the reader to
section 4.5 in \cite{AF} for the details of the proof of (\ref{n->n^4}). 
\\
Given $\nu\in\cM(\L_{2s_0})$ we
first estimate
$\Psi^{(\nu)}_{s_0}\bigl(\phi_{n,n^4}-\phi_{n,s_0}\bigr)$ as follows (a
similar bound will then be applied to any translation by $x$). \\
Assume, without loss of generality, that $s_0=N^4$
for some $N\in\NN$ and set $\ell_k :=k^4$ for any $k\in\NN$. Then,
given $0<\rho \ll 1$, by Schwarz inequality,
$$
\Psi^{(\nu)}_{s_0}(\phi_{n,n^4}-\phi_{n,s_0})\leq c_\rho
\sum_{k=n}^{N-1}k^{1+\rho}\,\Psi^{(\nu)}_{s_0}(\phi_{n,\ell_k}-\phi_{n,\ell_{k+1}}).
$$
In order to estimate
$\Psi^{(\nu)}_{s_0}(\phi_{n,\ell_k}-\phi_{n,\ell_{k+1}})$ we divide
$Q_{s_0}$ in cubes $\{Q_{i,k}\}_{i\in I_k}$ with side $\ell_k$ where,
without loss of generality, we assume that $s_0/\ell_k\in \bbN$ and
similarly for $\ell_k/n$. Let $\bar Q_{i,k}$ be the cube of side
$10\ell_k$ concentric to $Q_{i,k}$. Then by lemma \ref{miramare} with
\begin{equation*}
I:=I_{k+1}, \quad \L :=\L_{2s_0 },\quad \L_i:=\bar Q_{i,k+1},\quad
f_i:=\Av^{(n)}_{x\in Q_{i,k+1}}\t_x [\phi_{n,\ell_k} - \phi_{n,\ell_{k+1}} ]
\end{equation*}
we obtain (thanks also to lemma \ref{cleonimo})
\begin{equation}
\label{epancin}
\Psi_{s_0}^{(\nu)}(\phi_{n,\ell_k}-\phi_{n,\ell_{k+1}})\leq c\,\ell_{k+1}^{d+2}
\Av_{i \in I_{k+1 }}\Var_{\mu^{\l}}\Bigl(\Av^{(n)}_{x\in Q_{i,k+1}}
 \t_x [\phi_{n,\ell_k}-\phi_{n,\ell_{k+1}}] \Bigr)
\end{equation}
 where $\mu^{\l}$ is the grand canonical measure corresponding to $\nu$.
\\
Let now $J$ be the set of
possible densities on $\L_{2s_0}$.
Then, thanks to (\ref{epancin}), it is enough to prove that, for $\rho $
small enough and for almost any disorder $\a$,
\begin{equation}
\label{panipesci}
 \lim_{ n \uparrow \infty , \ep \downarrow 0 } \Av_{x \in \toromi}
\frac{1}{n^2} \sum_{k = n+1}^{N } k^{1+\rho } \ell_{k
}^{d+2} \Av_{i \in I_{k}}\sup_{m \in J }
 \Var_{\mu_{\L_{x,2s_0}}^{\l(m)}}
\Bigl(\Av^{(n)}_{y\in x+Q_{i,k} }\t_{y}\phi_{n,\ell_{k}} \Bigr)=0
 \end{equation}
and similarly with $\phi_{n,\ell_k}$ replaced by $\phi_{n,\ell_{k-1}}$. \\
Given $\g>0$ we set
$J_k=\{\ell_k^{-\g},2\ell_k^{-\g},\dots,1-\ell_k^{-\g}\}$. Then, using
(\ref{sam1}), the variance in (\ref{panipesci}) can
be bounded from above by
$$
\Var_{\mu_{x+\bar Q_{i,k}}^{\l(\bar m)}} \Bigl(\Av^{(n)}_{y\in
  x+Q_{i,k} }\t_{y}\phi_{n,\ell_{k}} \Bigr) + c \ell_k^{d-\g}
$$
provided that $\bar m \in J_k$ satisfies
$|\bar m -m|\leq \ell _k ^{-\g}$.\\
Therefore, by choosing $\g$ large enough, we can replace in
(\ref{panipesci}) $\mu_{\L_{x,2s_0}}^{\l(m)}$ by $\mu_{x+\bar
  Q_{i,k}}^{\l(m)}$ and $J$ by $J_k$. We can at this apply proposition
\ref{ulisse} to get that
\begin{equation}
  \label{ulisse1}
  \sup_{m \in J_k }
 \Var_{\mu_{\bar Q_{i,k}+x}^{\l(m)}}
\Bigl(\Av^{(n)}_{y\in x+Q_{i,k} }\t_{y}\phi_{n,\ell_{k}} \Bigr)
\le c\id_{\cA_{x,i,k}^c}(\a) {\ell_k}^{-2d+2\d}+\id_{\cA_{x,i,k}}(\a)\,,
\end{equation}
where $\cA_{x,i,k}$ is a set of disorder configurations in $x+ \bar
Q_{i,k}$ with $\bbP(\cA_{x,i,k}) \le \nep{-c\,\ell_k^\d}$, $\d>0$.
\\
Therefore
\begin{gather}
\text{ l.h.s. of (\ref{panipesci}) } \le
\lim_{ n \uparrow \infty , \ep \downarrow 0 }
\, n^{-2} \sum_{k = n+1 }^{N } k^{1+\rho}
\ell_{k }^{d+2} \Av_{x \in \toromi} \Av_{i \in I_{k}} \id_{
\cA_{x,i,k} } \nonumber \\
+ \lim_{ n \uparrow \infty , \ep \downarrow 0 } \,  n^{-2} \sum
_{k = n+1}^{N } k^{1+\rho } \ell_k^{2 - d + 2 \d } .
\label{pane1}
\end{gather}
The second addendum in the r.h.s. of (\ref{pane1}) is zero because of the
definition of $\ell_k$ and the condition $d \geq 3$.\\
Let us consider the first addendum  in the r.h.s. of (\ref{pane1}). By
Chebyschev inequality, for any $q>0$ and any $x,k$
\begin{align}
\label{flaminio1}
\PP\bigl(\,\Av_{i\in I_{k}}\id_{
\cA_{x,i,k}}\geq \ell_k^{-q}\,\bigr) &\leq
\PP\bigl(\,\exists i\in I_k \,:\, \id_{\cA_{x,i,k}}\geq
\ell_k^{-q}\,\bigr) \nonumber \\
&\leq s_0^d \,\ell_k^{\,q-d}\,e^{-c\,\ell_k^\d} .
\end{align}
Moreover, by setting
$\bar{\id}_{\cA_{x,i,k}}=\id_{\cA_{x,i,k}} - \PP \bigl (
\cA_{x,i,k} \bigr )$, we have for any $r\in \bbN$ and any $x,k$
\begin{equation}
\label{flaminio2}
 \PP \bigl (\,\Av_{i \in I_{k}} \id_{ \cA
_{x,i,k} } \geq \ell_k^{-q} \, \bigr ) \leq
 c_r  \ell_k^{2rq }  \EE \bigl [ \bigl (
\Av_{i \in I_{k}} \bar{\id }_{ \cA_{x,i,k}} \bigr )
^{2r} \bigr ] \leq c' _r  \ell_k^{2rq + dr } s_0^{-dr}
\end{equation}
By taking the geometric average of the two estimates (\ref{flaminio1})
and  (\ref{flaminio2}) we finally obtain
$$ \PP \bigl (\,\Av_{i \in I_{k}} \id_{ \cA
_{x,i,k} } \geq l_k^{-q} \, \bigr ) \leq c(q,r)s_0^{-d(r-1)/2}.
$$
It is enough at this point to choose $q$ and $r$ large enough, define
$$
\Theta_\ep:= \{\;  \exists x\in \toromi : \Av_{i \in I_{k}} \id_{
\cA_{x,i,k} }  \geq \ell_k^{-q} \quad \text{for some}\; k
\leq N\; \},
$$
and apply Borel-Cantelli lemma to get that also the first addendum in
the r.h.s. of (\ref{pane1}) is negligible.
\end{proof}

\subsection{From scale $s_k$ to scale $s_{k+1}$.}
Here we define precisely the sequence of length scales $s_k$ and
discuss the details of the inductive step $s_k\rightarrow s_{k+1}$
described section \ref{road}.
\\
Let $\{a_k\}_{k\geq0}$ be defined inductively by
$$
a_0= 1 \quad \text{ and } \quad a_{k+1}=1+\bigl(\ov2
-\frac{1}{2^{k+1}}\bigr)a_k
$$
It is easy to verify that the sequence $\{a_k\}_{k\geq 0}$ is
increasing with $\lim_{k\to \infty}a_k=2$. Let also $s_k :=
\ep^{-\frac{a_k}{d+2}}$.
\begin{prop}\label{scheme}
Let $d\geq 2$. Then
\begin{equation}
 \label{eq:scheme}
  \Av_{z,x}^{s_k, \frac{a}{\ep}} \t_z  \frac{\phi
_{s_k,s_k}}{s_k} \approx \Av_{z,x}^{s_{k+1}, \frac{a}{\ep}} \t_z
 \frac{\phi_{s_{k+1},s_{k+1}}}{s_{k+1}} \qquad
 \forall k \geq 0.
\end{equation}
\end{prop}
\begin{proof}
  In order to prove (\ref{eq:scheme}) observe that, by construction, the
  two exponents $\frac{a_k}{d+2}$ and $\frac{a_{k+1}}{d+2}$ satisfy the
  conditions of propositions \ref{tolmezzo} and \ref{paularo} with
$\g:=\frac{a_k}{d+2}$ and $\g':=\frac{a_{k+1}}{d+2}$. Therefore
  we have the following chain of equivalences:
\begin{equation*}
\Av_{z,x}^{s_k,\frac{a}{\ep}}\t_z\frac{\phi_{s_k,s_k}}{s_k}\approx
\Av_{z,x}^{s_k,\frac{a}{\ep}}\t_z\frac{\phi_{s_k,s_{k+1}}}{s_k}
\approx
\Av_{z,x}^{s_{k+1},\frac{a}{\ep}}\t_z\frac{\phi_{s_{k+1},2s_{k+1}}}{s_{k+1}}
\end{equation*}
Finally, using again proposition \ref{paularo} with $s=s_{k+1}$ and
$s'=2s$, we obtain (\ref{eq:scheme}).
\end{proof}

\subsection{Analysis of $\frac{\phi_{\bar s,\bar s}}{\bar s}$ via an
  improved Two Blocks Estimate}
Here we describe the final step in the proof of theorem \ref{pte},
namely we show that
$$
\Av_{z,x}^{\bar s,\frac{a}{\ep}}\t_z\frac{\phi_{\bar s,\bar s}}{\bar s}\approx 0
$$
where $\bar s = \ep^\d s_\infty$ and $s_\infty=\ep^{-\frac{2}{d+2}}$
 (see section \ref{road}).  The basic
tool is represented by the following improved version of the Two Blocks
Estimate (see e.g. \cite{KL}), whose  proof mainly relies on the same techniques used for proving proposition \ref{tvb} (see section 4.10 in \cite{AF}).
\begin{prop}
\label{super-two-blocks} (Improved Two Blocks Estimate)\\
  Let $d \geq 3$, $0<\gamma<\g'<1$ and set $s=\ep^{-\g}$, $\ell=\ep^{-\g
    '}$. Then, for any $r$ such that $0<r<\min\bigl(\frac{2(1-\g
    ')}{d+4}, \frac{\gamma}{2}\bigr)$ and for almost any disorder
  configuration $\a$
\begin{equation*}
   \limsup_{a\dto 0,\ep \dto 0}\Av_{x\in\toromi}
\sup_\nu \sup spec_{L^2(\nu)}\bigl\{\ep^{-r}
   \Av_{w,x}^{\ell,\frac{a}{\ep}}\,Av_{z,w }^{s, \ell}\,
|m_{z,s}^e-m_{w,\ell}^e|+\ep^{-2}
\Av_{b\in\L_{x,2\frac{a}{\ep}}}\cL_b\bigr\}\leq0
\end{equation*}
where $ \nu $ varies among $\cM(\L_{x,2\frac{a}{\ep}})$.
\end{prop}

\begin{cor}
  Let $d\geq 3$ and $0<\delta\ll 1$. Then
\begin{equation}
\label{eq:last_step}
 \Av_{z,x}^{\bar s,\frac{a}{\ep}}\t_z\frac{\phi_{\bar s,\bar s}}{\bar
   s}\approx 0.
\end{equation}
\end{cor}
\begin{proof}
  For simplicity of notation we omit the bar in $\bar s$ and we set $\De
  m:=m_s^{2,e}-m_s^{1,e}$ and $N:=N_{\L_s^e}$. \\
  Let $\hat\phi_{s,s}(\h) = \mu_{\L_s^e}^{\l(m_s^e(\h))}(\D m)$.
  Then, by the equivalence of ensembles (see proposition \ref{equi}), it is
  enough to prove (\ref{eq:last_step}) with $\phi_{s,s}$ replaced by
  $\hat \phi_{s,s}$.  Let $m$ be a particle density on $\L^e_s$ that,
  without loss of generality, we can suppose in $(0,\ov2)$ and set $\l:=
  \l_{\L_s^e}(m)$ and  $\l_0:=\l_0(m)$.
  Then, by Taylor expansion,
\begin{equation}
\label{taylor}
 \mu^{\l}(\De m)= \mu^{\l_0}(\De m)+
 \mu^{\l_0}(\De m;N)(\l-\l_0)+\mu^{\l'}(\De m;N;N)(\l-\l_0)^2
\end{equation}
where $\l'$ is between
$\l$ and $\l_0$.\\
Let us observe that $ |\mu^{\l'}(\De m;N;N)| \le c$, while by lemma \ref{vespe}
$$
|\l-\l_0|\leq c\, |1-\frac{\mu^{\l_0}(m_{s}^e)}{m}|.
$$
Moreover, $\EE\bigl[\mu^{\l_0}(m_{s}^e)\bigr]=m$ and $\EE\bigl
[\mu^{\l_0}(\Delta m_{s} ; N)\bigr]=0$. Therefore, thanks to the large
deviations estimate of lemma \ref{upupa} applied to the function
$f(\a):= \frac{\mu^{\l_0}(\h_0)}{m}-1$, for any $\b \in (0,1)$ and $\ep$
small enough
$$
\bbP(|\l -\l_0| \ge s^{-\frac{d}{2} +\frac{\b}{2}} ) \le
\bbP(|\Av_{x\in \L_s^e} \t_x\,f| \ge
\frac{1}{c}s^{-\frac{d}{2}+\frac{\b}{2}} ) \le \nep{-cs^\b}.
$$
A similar reasoning applies to the term $\mu^{\l_0}(\De m;N)$ if we
consider instead the function $f(\a):= \mu^{\l_0}(\h_0;\h_0) -
\bbE\bigl(\mu^{\l_0}(\h_0;\h_0)\bigr)$.  The above bounds together with
the fact that the number of possible choices of $m$ is polynomially
bounded in $s$ and together with Borel Cantelli lemma, implies in particular that for almost all disorder
configuration $\a$ and for $\ep$ small enough
\begin{equation*}
 \sup_{x\in\toromi}\|
\t_x( \,\hat\phi_{s,s}-\mu^{\l_0(m_s^e)}(\De m)\,)\|_\infty\leq s^{-d+\b}.
\end{equation*}
Thanks to the above estimate it is enough to prove (\ref{eq:last_step})
with $\phi_{s,s}$ replaced by $\mu^{\l_0(m_s^e)}(\De m)$, that is
\begin{equation}
\label{kolmogorov100}
 \Av_{z,x}^{s,\frac{a}{\ep}}\mu^{\l_0(m_{z,s}^e)}\Bigl(\t_z\frac{m^{2,e}_s-
m^{1,e}_s }{s}\Bigr) \approx 0.
\end{equation}
We assert that we only need to show that
\begin{equation}
\label{cantor}
\text{l.h.s. of }(\ref{kolmogorov100})\approx
\Av_{z,x}^{\ell,\frac{a}{\ep}}\mu^{\l_0(m_{z,\ell}^e)}\Bigl(\t_z
\frac{m^{2,e}_\ell-m^{1,e}_\ell}{\ell}\Bigr)
\end{equation}
where  $\ell=\ep^{1-\rho}$ is a new mesoscopic scale  with
$0<\rho<1$ so small that $s<\ell$ and $\ep^{-1} \ell^{-\frac{d+2}{2}}
\downarrow 0$ as $\ep \downarrow 0$. In fact, thanks to   lemma \ref{upupa} applied with $f(\a):=\mu^{\l_0(m_{z,\ell}^e)}(\eta_0-\eta_{\ell e})$,   given $0<\b\ll 1$ for almost any disorder configuration $\a$ and for $\e$ small enough the r.h.s. of (\ref{cantor}) is bounded by $\ell^{-\frac{d+2}{2}+\b}$. Because of our choice of $\ell$, the r.h.s. of (\ref{cantor}) is equivalent to $0$.\\
Let us prove (\ref{cantor}). To this aim, we observe that
thanks to (\ref{bobo}) and (\ref{merlino})
\begin{equation*}
\begin{split}
&  \text{l.h.s. of }(\ref{kolmogorov100})=
\Av_{w,x}^{\ell,\frac{a}{\ep}}\,\Av_{z,w }^{s,\ell}\,\mu^{
\l_0(m_{z,s}^e)}\Bigl (\t_z \frac{\D m }{s}\Bigr),\\
& \text{r.h.s. of }(\ref{cantor})=\Av_{w,x}^{\ell,\frac{a}{\ep}}\,\Av_{z,w }^{s, \ell}\,\mu^{ \l_0(m_{w,\ell}^e)} \Bigl(\t_z \frac{\D m}{s}
\Bigr).
\end{split}
\end{equation*}
Therefore, we only need to prove that
\begin{equation*}
\Av_{w,x}^{\ell,\frac{a}{\ep}}\,\Av_{z,w }^{s,\ell}\Bigl (\mu^{
\l_0(m_{z,s}^e)}\bigl (\t_z \frac{\D m }{s}\bigr)-
\mu^{ \l_0(m_{w,\ell}^e)}\bigl (\t_z \frac{\D m }{s}\bigr)\Bigr )\approx 0 .
\end{equation*}

Let us assume for the moment that, given $0<\b\ll 1$,
 for almost all disorder configuration $\a$ and $\ep$
small enough
\begin{equation}
 \label{eq:lip_TBE}
\sup_{x \in \toromi} |\mu^{\l_0(m)}\bigl(\t_x\De
m\bigr)-\mu^{\l_0(m')} \bigl(\t_x\De m\bigr)|\leq cs
^{-\frac{d}{2} +\b}|m-m'|+cs^{-\frac{d}{2}-\b } \quad \forall m,m' \in [0,1].
\end{equation}
Then it is simple to deduce (\ref{cantor}) from (\ref{eq:lip_TBE}) and
proposition \ref{super-two-blocks} with $\g=\frac{2}{d+2}-\d\,,\; \g ' =
1 -\rho$ and  $r=-\d\b +\frac{d+2}{2}\d+\frac{2}{d+2}\b$ by choosing suitable
$0<\b\ll\d\ll\rho\ll 1$.\\
It remains to  prove (\ref{eq:lip_TBE}). For simplicity of notation, let us consider only the case $x=0$ (the general case is a simple variation). By continuity, we may assume $0<m<m'<1$ and by Taylor expansion,
\begin{equation*}
\begin{split}
|\mu^{\l_0(m')}\bigl(\D m\bigr)-\mu^{\l_0(m)} \bigl(\D m\bigr)|& =
|\mu^{\l_0(\bar m)}(\D m;N)\l_0'(\bar m)(m'-m)|\\
& \leq c\, |\frac{ \mu^{\l_0(\bar m)}(\D m;N)}{\bar m}|(m'-m)
\end{split}
\end{equation*}
where $m<\bar m<m'$. If we could restrict the possible values of  $\bar m$ to
 $\{s^{-d},2s^{-d}, \dots, 1-s^{-d}\}$, then, by means of large
deviations estimate as in the first part of the proof, we would obtain
$\frac{1}{\bar m} |\mu^{\l_0(\bar m)}(\D m;N)|\leq c\, s^{-\frac{d}{2}+\b}$
  for almost any
disorder  $\a$ and for $\ep$ small enough, thus  implying  (\ref{eq:lip_TBE}). The complete proof requires some addional  straightforward computations (see also section $4.10$ in \cite{AF}).
\end{proof}

\section{Some technical results needed in section \ref{Fluct}}\label{Tech}
In this section we collect some technical results, mostly based on
estimates of large deviations in the disorder field $\a$, that are used
in the proof of theorem \ref{pte}. Our bounds mainly concern canonical
or grand canonical variances of suitable spatial averages of local
functions. Such variances arise naturally from eigenvalue estimates via
perturbation theory. We have seen in fact that, when perturbation theory
applies (see proposition \ref{jamme}), the maximal eigenvalue is bounded
by an expression containing an H$_{-1}$ norm that, in general, can be
bounded from above by:
\begin{equation}
\label{H-1}
\nu\bigl(f,-\cL^{-1}_{\L}f\bigr)\leq c\,\ell^{2}\Var_{\nu} (f) \le
c\, \ell^{2}\Var_{\mu}(f)
\end{equation}
where $\nu$ is a canonical measure on the cube $\L$ of side $\ell$ with
disorder $\a$, $\mu$ is the corresponding grand canonical measure (with
suitable empirical chemical potential) and $f$ is a (mean zero w.r.t.
$\nu$) function.  Above we used the spectral gap bound $\gap(\cL_{\L})
\ge c \ell^{-2}$ together with lemma \ref{cleonimo}.

When the function $f$ is the spatial average of local functions
$\{f_i\}_{i\in I}$ each with support much smaller than $\L$ it is
possible to do better than (\ref{H-1}). We have in fact:
\begin{lemma}\label{miramare}
  Let $\L$ be a box in $\ZZ^d$ and $\{\L_i\}_{i\in I} $ be a family of
  cubes $\L_i\subset \L$ with side $R$ satisfying
\begin{equation*}
\bigl|\,\{i\in I\,:\,x\in \L_i\}\,\bigr|\leq10^{10d}\quad \forall x\in \L.
\end{equation*}
Let $f=\Av_{i\in I}f_i$ where, for any $i\in I$ and for all $\a$, $f_i$
has support in $\L_i$ and has zero mean w.r.t. any canonical measure on
$\L_i$.  Then, for any canonical measure $\nu$ on $\L$ with disorder
configuration $\a$,
$$
\nu \bigl (f,-\cL_\L^{-1} f\bigr) \leq c\, R^2 |I |^{-1}
\Av_{i\in I}\nu\bigl(\,\Var_\nu(f_i \tc \cF_i) \, \bigr ).
$$
\end{lemma}
\begin{proof}
  Let $\cF_i :=\s(m_{\L_i},\h_x\;\text{with}\;x\not\in \L_i)$ and
  observe that
$$
\nu(f_i,g)=\nu\bigl(\nu(f_i;g\tc \cF_i)\bigr) \quad \forall \,g
$$
Thus, by Schwarz and Poincar\'e inequalities and the diffusive
scaling of the spectral gap
\begin{align*}
  | \nu (f,g) |&\leq c\, R \, \Av_{i \in I }\, \nu \bigl \{\, \bigl [
  \Var_\nu ( f_i \tc \cF_i ) \dstor_{ \L_i } \bigl ( g ; \nu ( \cdot
  | \cF_i )
  \bigr ) \bigr ] ^{1/2 }  \,\bigr\} \\
  &\leq c\, R\, |I|^{-1/2}\Bigl( \Av_{i \in I }\,\nu\bigl(\, \Var_\nu
  (f_i\tc \cF_i) \bigr)\Bigr)^{1/2} \dstor_\L ( g; \nu)^{1/2}.
\end{align*}
It is enough now to take $g=-\cL^{-1}_\L f$.
\end{proof}

\subsection{Variance bounds.}
\label{Var-Bounds}
One of the key issues is to provide sharp enough upper bounds (see
proposition \ref{bound:variance} below) on the variance
\begin{equation}
 \label{tropici}
   \Var_{\mu^{\l_0(m)}}\bigl(\Av_{x\in\L_k}\t_x\phi_{n,s}\bigr)
\end{equation}
where $ n,s,k$ are positive integers satisfying $n\leq s\leq k $ and
$m\in (0,\ov2)$ and $\phi_{n,s}$ has been defined in (\ref{phi(n,s)}).
Actually the method developed below is very general and it can be used
to estimate also other similar variances, like for example
(\ref{tropici}) with $\l_0(m)$ replaced by the empirical chemical
potential $\l_{\L_k}(\a,m)$.

It is convenient to define first some additional convenient notation
besides those already defined at the beginning of section \ref{Fluct}:
\begin{align}
\label{definition1}
\hat\phi_{n,s}(\h )&:=\mu_{\L_s^e }^{\l(m_s^e) }(m_n^{2,e}-m_n^{1,e}) \nonumber \\
 \xi_0(m)&:=\mu^{\l_0(m)}(\,m_n^{2,e}-m_n^{1,e}\,;N_{\L_n^e})\nonumber \\
\xi(m)&:=\mu_{\L_s^e}^{\l(m)}(\,m_n^{2,e}-m_n^{1,e}\,;N_{\L_s^e})
\nonumber \\
 \s^2_0(m)&:=\mu^{\l_0(m)}(\,m_s^e;N_{\L_s^e}) \nonumber \\
 \s^2(m)&:=\mu_{\L_s^e}^{\l(m)}(\,m_s^e;N_{\L_s^e}),
\end{align}
where $N_{\L_n^e}$, $N_{\L_s^e}$ denote  the particle number respectively
in the box $\L_n^e$ and $\L_s^e$. Let us recall  the definition of static
compressibility $\chi(m)=  \mathbb E\bigl(\mu^{\l_0(m)}(\,\h_0;\h_0\,)\bigr)$.\\
Moreover, given $0<\d\ll 1$ and a site $x$, we define the events:
\begin{align}
\label{definition2}
\cM_x(m)&:=\{|m^e
_{x,s}(\h)-m|\geq\sqrt{m}\,s^{-\frac{d}{2}+\frac{\d}{2}}\,\} \nonumber \\
\cA^{(1)}_x(m)&:=\bigl\{\,\frac{1}{m}|m-\mu^{\l_0(m)}(m^e_{x,s})|\geq
s^{-\frac{d}{2}+\frac{\d}{2}}\,\bigr\}\nonumber
\\
\cA^{(2)}_x(m)&:=\bigl\{\,\bigl|\t_x\frac{\s^2_0(m)}{\chi(m)}
-1\bigr|\geq s^{-\frac{d}{2}+\frac{\d}{2}}\,\bigr\}
\end{align}
\begin{remark}
  Notice that the first event is an event for the particles
  configuration $\h$ while all the others are events for the disorder
  field.
\end{remark}
\begin{lemma}
\label{montalbano} There exists $s_0(\d)$ such that the following holds for any $s\geq s_0
(\d)$. Assume $n \leq s$, $4s^{-d+\d}\leq m\leq 1/2$, $\h\not\in
\cM_x(m)$ and $\a\not\in\cA^{(1)}_x(m)\cup \cA^{(2)}_x(m)$.  Then, for
any site $y$,
\begin{equation}
 \label{gradient.phi}
 \Bigl|\nabla_y[\t_x\hat\phi_{n,s}](\h )-\frac{(1-2\h_y)}{2s^d}
\t_x\frac{\xi_0( m)}{\chi(m)}\Bigr|\leq c\,s^{-d}\Bigl\{\frac{s^{-d}}{m}+
\frac{1}{\sqrt{m}}s^{-\frac{d}{2}+\frac{\d}{2}}\Bigr\}.
\end{equation}
\end{lemma}
\begin{proof}
By Lagrange theorem we can write
\begin{equation}
\label{gradient.phi1} \nabla_y [\t_x\hat\phi_{n,s}](\h )=\int
_{m_{x,s}^e(\h)}^{m_{x,s}^e(\h^y)}\t_x\frac{\xi(m')}{
\s^2(m')}\, dm'.
\end{equation}
Assume $m'$ in the interval with end-points $m^e_{x,s}(\h)$ and
$m^e_{x,s}(\h^y)$. Then, by lemma \ref{turu},
$$
\xi_0 (m')\leq c\,m',\quad
\xi(m ')\leq c\,m', \quad \s^2_0(m')\geq  c\, m',\quad \s^2( m')\geq c\, m',
\quad \chi(m')\geq c\,m
.
$$
Moreover, since $\h\not\in \cM_x(m)$, $m'\geq c\, m $ if $s$ is large
enough depending on $\d$.  Therefore, by lemma $\ref{vespe}$

\begin{align}
  \Bigl|\t_x\frac{\xi(m')}{\s^2(m')}-\t_x\frac{\xi_0(m
    ')}{\s_0^2(m')}\Bigr |&\leq\frac{c}{m}|m
  '-\mu^{\l_0(m')}(m^e_{x,s})|,
\label{magoz1} \\
\Bigl|\t_x\frac{\xi_0 (m ')}{\s_0
 ^2(m')}-\t_x\frac{\xi_0(m)}{\s_0^2(m)}\Bigr
|&\leq\frac{c}{m}|m'-m|\leq \frac{c}{m}|m^e_{x,s}(\h)-m|+\frac{c}{m}s^{-d} , \label{magoz2} \\
\Bigl |\t_x \frac{ \xi_0(m)}{\s
  _0^2(m)}-\t_x\frac{\xi_0(m)}{\chi(m)}\Bigr|&
\leq c \Bigl|\t_x \frac{\s_0^2(m)}{\chi(m)}-1 \Bigr|
\label{magoz3} .
\end{align}
By lemma \ref{vespe} and the assumption $\a \not\in \cA^{(1)}_x(m)$, the
r.h.s. of (\ref{magoz1}) can be bounded from above by

\begin{equation}
  \label{errore1}
\frac{c}{m}|m-\mu^{\l_0(m)}(m^e_{x,s})|+\frac{c}{m}s^{-d} \le
c\bigl[\, s^{-\frac{d}{2}+\frac{\d}{2}} + \frac{1}{m}s^{-d}\,\bigr].
\end{equation}
Similarly, the contribution of the r.h.s. of (\ref{magoz2}) together
with (\ref{magoz3}) can be bounded from above by
\begin{equation}
  \label{errore2}
  c\Bigl[\frac{s^{-d}}{m}+
\frac{1}{\sqrt{m}}s^{-\frac{d}{2}+\frac{\d}{2}}\Bigr]
\end{equation}

\smallno The thesis follows immediately from (\ref{gradient.phi1})
together with (\ref{errore1}), (\ref{errore2}).
\end{proof}

\begin{lemma}
\label{stima}  
There exists $s_0(\d)$ such that the following holds for any $s\geq s_0
(\d)$. Let $n\leq s$, $m\in(0,\ov2)$, and let $\l(\a)$ be a bounded
measurable function such that for any disorder configuration $\a$
\begin{equation}\label{eq:stima1}
|\l(\a)-\l_0(m)|\leq s^{-\frac{d}{2}+\frac{\d}{4}}.
\end{equation}
Then, for any $s\geq s_0(\d)$ and any finite set $\Delta\subset\ZZ^d$,
\begin{equation}\label{eq:stima2}
\PP\bigl(\,\mu^{\l(\a)}\bigl(\cup_{x\in\Delta} \cM_x
(m)\bigr)\geq|\Delta|e^{-s^{\d/2}}\,\Bigr)\leq
c\,e^{-s^{\d/2}}.
\end{equation}
\end{lemma}
\begin{proof}
  By the Chebyshev inequality and the translation invariance of $\PP$,
  the l.h.s. of $(\ref{eq:stima2})$ can be bounded from above by
  $\exp(s^{\d/2 }) \EE
  \bigl[\mu^{\l(\a)}(\cM_0(m))\bigr].$ \\
  Let us bound the term
\begin{equation}
\label{eq:stima3}
e^{s^{\frac{\d}{2}}}\EE\bigl[\,\mu^{\l(\a)}\bigl(m^e_s-m\geq\sqrt{m}\,
s^{-\frac{d}{2}+\frac{\d}{2}}\bigr)\bigr].
\end{equation}
Thanks again to Chebyshev inequality, for any $0<t <1$ (\ref{eq:stima3})
can be bounded from above by
\begin{equation}
\label{velo}
e^{s^{\frac{\d}{2}}-2t\sqrt{m}s^{\frac{d}{2}+\frac{\d}{2}}}\EE\bigl[\prod
_{x\in\L_s^e} \mu^{\l(\a)}\bigl( e^{t( \h_x - m ) }
\bigr )\bigr].
\end{equation}
Using the basic assumption (\ref{eq:stima1}) and Lagrange theorem, it is
not difficult to see that
$$
\mu^{\l(\a)}\bigl(e^{t(\h_x-m)}\bigr)\leq(1+c\,tms^{-\frac{d}{2}+
  \frac{\d}{4}})\mu^{\l_0(m)}\bigl(e^{t(\h_x-m)}\bigr)
$$
so that (\ref{velo}) is bounded from above by
\begin{equation*}
  e^{s^{\frac{\d}{2}}-2t\sqrt{m}s^{\frac{d}{2}+\frac{\d}{2}}+
  c\,tm s^{\frac{d}{2}+\frac{\d}{4}}}\EE\bigl[\mu^{\l_0}
  \bigl(e^{t(\h_0-m)}\bigr)\bigr]^{2s^d}.
\end{equation*}
Since $e^x\leq 1+x+2\,x^2$ if $|x|\leq1$, the above expression is
bounded from above by
\begin{equation*}
 \exp\bigl(s^{\frac{\d}{2}}-2t\sqrt{m}\,s^{\frac{d}{2}+\frac{\d}{2}}+ c\,t\,m s^{
 \frac{d}{2}+ \frac{\d}{4}}+c\,t^2 m s^d \bigr ).
\end{equation*}
The thesis follows by choosing $t$ such that $t^2 m=s^{-d +\d/2}$.
\end{proof}
We are finally in a position to state our main bound on the variance
appearing in (\ref{tropici}).
\begin{prop}
\label{bound:variance} 
For $d\geq 2$ there exists $s_0(\d)$ such that the following holds for
any $s\geq s_0 (\d)$. Let $m\in(0,\ov2)$ and let $n\leq s\leq k\leq 1000 s
$. Then
there exists a measurable set $\cA$ with $\PP(\cA) \leq k^{2d} e^{
  -\,cs^{\d/2} }$ such that
\begin{equation}
\label{patti}
 \Var_{\mu^{\l_0(m)}}\left (\Av_{x\in\L_{k}}\t_x \phi_{n,s}\right)\leq
 c\id_{\cA^c}(\a) s^{-2d +2\d}+\id_{\cA}(\a).
\end{equation}
\end{prop}
\begin{proof}
  Let us consider first the case of ``low density'' $m\leq 4s^{-d+\d}$.
  \\
  Since $|\t_x \phi_{n,s}|\leq c\, m_{x,s}^e$, $|\Av
  _{x\in\L_k}\t_x\phi_{n,s}|\leq c\,m_{\L_{2k}}$ and therefore the
  l.h.s. of (\ref{patti}) can be bounded from above by
  $$
  \mu^{\l_0(m)}(m_{\L_{2k}}^2)\leq c (k^{-d}m+m^2) \leq c
  s^{-2d+2\d}.
 $$
 Let us now consider the ``high density'' case $m\geq 4s^{-d+\d}$. \\
 By the equivalence of ensembles (see proposition \ref{equi}), in the
 l.h.s. of $(\ref{patti})$ $\phi_{n,s}$ can be substituted by
 $\hat\phi_{n,s}$ with an error of order $s^{-2d}$. Therefore, by the
 Poincar\'e inequality
\begin{equation}
  \label{Poincare'}
  \Var_{\mu^{\l_0(m)}}(f) \le c \,m \,\mu^{\l_0(m)}(\sum_{y}|\nabla_yf|^2),
\end{equation}
it is enough to estimate
\begin{equation}
\label{itacaitaca}
  c\,m\,\mu^{\l_0(m)}\Bigl[\frac{1}{k^{2d}}
\sum_{y\in \L_{2k}}\bigl(\,\sum_{x\in\L_k\cap\L_{y,s}}\nabla_y[\t_x\hat
 \phi_{n,s}]\,\bigr)^2\Bigr].
\end{equation}
To this aim we set (recall (\ref{definition2}))
\begin{align}
  &\cM:=\cup_{x\in\L_k}\cM_x(m) \qquad \cA_0 :=\bigl\{\mu^{\l_0(m )}\bigl(\cM
  \bigr)\geq k^d\exp(-s^{\d/2 } )\,\bigr\}, \nonumber\\
  &\cA_1:=\cup_{x\in\L_k}\cA^{(1)}_x (m) \qquad
  \cA_2:=\cup_{x\in\L_k}\cA^{(2)}_x(m)\nonumber\\
  &\cA_3:=\cup_{y\in\L_{2k}}\Bigl\{\,\Bigl|\Av_{x\in\L_k\cap\L_{y,s}}
  \Bigl[\t_x\, \frac{\xi_0(m)}{\chi (m)}\Bigr]\,\Bigr|
  \geq|\L_k \cap \L_{y,s}|^{-\frac{1}{2}+\frac{\d}{2d} }\,\Bigr\},\nonumber\\
  &\cA :=\cA_0 \cup \cA_1\cup\cA_2 \cup \cA_3. \nonumber
\label{events}
\end{align}
We first estimate
\begin{equation}
\label{star} \id_{\cA^c }(\a)m \,\mu^{\l_0(m)}\Bigl[
\id_{\cM^c}\frac{1}{k^{2d}}\sum_{y\in \L_{2k}}
\bigl(\,\sum_{x\in\L_{k}\cap \L_{y,s} }\nabla_y[\t_x \hat \phi_{n,s}]\,\bigr)^2 \Bigr].
\end{equation}
By lemma $\ref{montalbano}$, for $s$ large enough (\ref{star}) can be
bounded from above by
\begin{equation}
 \label{star1}
 \frac{c}{k^{2d}}\id_{\cA^c} (\a)  \Av_{y\in \L_k}
 \Bigl[\,\frac{1}{s^{d}}\sum_{x\in \L_k \cap \L_{y,s}}\t_x\frac{\xi_0(m)}{
\chi(m)}\Bigr]^2 +
 c\,\frac{s^{-d +\d}}{k^d}
\end{equation}
By straightforward computations and the definition of $\cA_3$ the first
addendum in $(\ref{star1})$ can be bounded by $c\,k^{-d}s^{-d+\d}$.
Moreover, because of the definition of $\cA_0 $, expression
$(\ref{star})$ with $\id _{\cM^c}$ replaced by $\id_{\cM}$ can be
bounded by $c\,s^{2d}e^{-s^{\d/2}}$.
\\
In conclusion
\begin{equation*}
  \id_{\cA^c }(\a)m \,\mu^{\l_0(m)}\Bigl[\frac{1}{k^{2d}}\sum_{y\in \L_{2k}}
\bigl(\,\sum_{x\in\L_{k}\cap \L_{y,s} }\nabla_y[\t_x \hat
\phi_{n,s}]\,\bigr)^2 \Bigr] \le
c\,\bigl[\,k^{-d}s^{-d+\d} +  s^{2d}e^{-s^{\d/2}}\,\Bigr].
\end{equation*}
\medno
It remains to prove that $\PP(\cA)\leq k
^{2d}e^{-c\,s^{\d/2}}$. To this aim we set
\begin{align*}
  f_1 (\a)&:=1- \mu^{\l_0(m)}(\h_0)\,, \\
  f_2 (\a)&:=1-\mu^{\l_0(m)}(\h_0;\h_0)/\s_0^2(m) \,,\\
  f_3(\a)&:=\bigl(\mu^{\l_0(m)}(\h_{ne};\h_{ne})-\mu^{\l_0(m)}(\h_0;\h_0)\bigr)/\chi(m).
\end{align*}
By lemma \ref{stima} $\PP(\cA_0) \leq c\,e^{-s^{\d/2}}$ while
$\PP(\cA^{(1)}_x (m))$ and $\PP(\cA^{(2)}_x (m))$ can be bounded from
above by $e^{-c\, s^{\d}}$ by means of lemma $\ref{upupa}$ with $f=f_1$
and $f=f_2$ respectively.  Therefore
\begin{equation}
  \label{A1A2}
\PP(\cA_1)+\PP(\cA_2) \le
k^{d}e^{-c\,s^{\d}}
\end{equation}
In order to bound $\PP(\cA_3)$ we observe that
\begin{equation*}
\Av_{x\in\L_{k}\cap\L_{y,s}}\t_x\frac{\xi_0(m)}{\chi(m)}
= \Av_{z\in\L_n^{1,e}}\Av_{x\in\L_{k}\cap\L_{y,s}}\t_{x+z}f_3\,.
\end{equation*}
Thus
\begin{equation*}
\cA_3 \subset\cup_{y\in\L_{2k}}\cup_{z\in\L_n
^{1,e}}\cA_3 (y,z)
\end{equation*}
where
$$
\cA_3 (y,z)= \{|\Av_{x\in\L
  _{k}\cap\L_{y,s}}\t_{x+z}f_3|\geq|\L_{k}\cap\L_{y,s}|^{-\frac{1}{2}+\frac{\d}{2d}}\,\}\,.
$$
Using once more lemma \ref{upupa} we get
$$
\PP(\cA_3 (y,z))\leq
\exp (-c s^{\frac{d-1}{d}\d})
$$
and the proof is complete.
\end{proof}
We conclude this part with a slight modification of proposition
\ref{bound:variance}.

\begin{prop}
\label{ulisse}
  Let  $n\leq s$ be positive integers and let $0<\d\ll 1$. Let
  also $\g >0$ and set $J_s = \{1/s^\g, 2/s^\g,\dots 1-1/s^\g\}$. Then
  there exists a set $\cA$ of disorder configurations $\a$ in $\L_{2s}$
  satisfying
$$
\PP(\cA)\leq s^\g \,e^{-c \,s^\d}
$$
and such that, for $s$ large enough depending on $\d$,
\begin{equation}
\label{lebedev}
\sup_{m\in J_s}\Var_{\mu_{\L_{2s}}^{\l(m)}}\left(\Av^{(n)}_{x\in\L_{s}}\t_x\phi_{n,s}\right)
\leq c\id_{\cA^c}(\a) s^{-2d+2\d}+\id_{\cA}(\a)
\end{equation}
where $\Av^{(n)}_{x\in\L_s}:=\Av_{x\in \L_s\cap n\ZZ^d}$.
\end{prop}
\begin{proof}
The proposition can be proved as proposition \ref{bound:variance} with some slight modifications that we comment. For any $m\in J_s$ it is convenient to define  $\cM(m), \, \cA_1(m),$ and $\cA_2(m)$ as done  respectively for  $\cM, \, \cA_1,$ and $\cA_2$ in  the proof of proposition  \ref{bound:variance} and to set
\begin{align*}
\cA_0 (m)&:=\bigl\{\mu_{\L_{2s}}^{\l(m)}\bigl(\cM
\bigr)\geq s^d\exp(-s^{\d/2 } )\,\bigr\}, \\
\cA_3 (m)&:= \Bigl\{\,\Bigl|\Av^{(n)}_{x\in\L_s}\t_x\,
\frac{\xi_0(m)}{\chi(m)}\Bigr|\geq s^{-\frac{d}{2}+\frac{\d}{2}}\,\Bigr\}.
\end{align*}
Then one sets again $\cA(m) := \cA_0(m)\cup \cA_1(m) \cup \cA_2(m) \cup \cA_3(m)$, $\cA
:= \cup_{m\in J_s} \cA(m)$. By the same arguments as in the proof of proposition \ref{bound:variance} one obtains  (\ref{lebedev}).\\
Let us prove the estimate  $\PP(\cA)$ or, equivalently, that for any $m\in J_s$ $\PP(\cA(m))\leq e^{-c s^\d}$. For this purpose, given $m\in J_s$, it is convenient to define
$$
\cB(m) := \bigl\{\,\bigl|\l_{\L_{2s}}(m)-\l_0(m)\bigr|\geq
s^{-\frac{d}{2}+\frac{\d}{4}}\,\bigr\}
$$
and write
\begin{equation}
\bbP(\cA(m)) \le \bbP(\cB(m)) + \bbP(\cB^c(m) \cap \cA_0(m)) +
\bbP(\cA_1(m)) + \bbP(\cA_2(m)) + \bbP(\cA_3(m)).
\label{stepB}
\end{equation}
Let us suppose $0<m\leq \ov 2$. Then  lemma \ref{vespe} implies that
$$\bigl|\l_{\L_{2s}}(m)-\l_0(m)\bigr|\leq c\, \bigl |1-m^{-1} \mu
^{\l_0(m)}(m_{\L_{2s}})\bigr|.$$
Thanks to the above estimate and to lemma \ref{upupa} applied with $f:= 1-
 m^{-1} \mu^{\l_0(m)}(\h_0)$,
the first term in the r.h.s. of (\ref{stepB}) is smaller than
$\nep{-c\,s^{\d/2}}$. The second term is smaller
than $\nep{-c s^{\d/2}}$ by lemma \ref{stima}. Moreover, $\PP(\cA_1 (m))$ and
$\PP (\cA_2(m))$ can be bounded by $ s^d e^{-c\, s^\d}$ as in the proof
 of proposition \ref{bound:variance}.\\
Finally, let us consider $\PP(\cA_3 (m))$.  For simplicity of notation we restrict to the case  $d=1$ and  we write
\begin{equation}
\label{sligovitsch}
n \sum _{x\in\L_s\cap n\ZZ} \t_x\,\frac{\xi_0(m)}{\chi(m)} =
\Bigl(\sum_{x\in\L_s\cap 2n\ZZ}\sum_{z\in\L_{1,n}^e} \t_{x+z}f \Bigr) +\Bigl ( \sum _{x\in\L_s\cap (2n+n)\ZZ} \sum _{z\in\L_{1,n}^e}\t_{x+z}f \Bigr)
\end{equation}
where
$f:= \chi (m) ^{-1} \bigl (\, \mu^{\la_0(m)}(\h_0;\h_0)-
\mu^{\la_0(m)}(\h_{ne};\h_{ne})\,\bigr)$.
 We remark that  in both the addenda in the r.h.s. of (\ref{sligovitsch}) the appearing functions have disjoint support and form a set of cardinality $O(k^d)$, moreover $\EE(f)=0$. Therefore, by the same arguments used in the proof of lemma \ref{upupa}, we obtain that $\PP (\cA_3(m)) \leq e^{ -c s^\d}$.
\end{proof}

\subsection{An $L^\infty$ bound}
We conclude this section with a simple $L^\infty$ bound on
$|\t_x\phi_{s,s'}|$ when $s$ scales as an inverse power of $\ep$.

\begin{lemma}\label{norma_uniforme}
Let $0<\g<1 $ and $0<\d\ll 1$ and set  $s= O(\ep^{-\g})$. Then,
for almost all configuration disorder $\a $ and $\ep$ small enough,
\begin{equation}
\label{eq:norma_uniforme}
\sup_{x\in\toromi}|\t_x\phi_{s,s'}|\leq c \,s^{-\frac{d}{2}+\d} \quad
\forall s'\in [s,\ep^{-1}].
\end{equation}
\end{lemma}
\begin{proof}
By the equivalence of ensembles it is enough to prove
(\ref{eq:norma_uniforme}) with $\phi_{s,s'}$ replaced by
$\hat\phi_{s,s'}$. Using lemma $\ref{vespe}$ we get
\begin{equation}
 \label{fabrizio}
| \hat \phi_{s ,s' }-\mu^{\l_0(m_{s'}^e)}\bigl(m_s^{1,e}-m_s^{2,e}\bigr)|\leq
  c\,|m_{s'}^e(\h)-\mu^{\l_0}(m_{s'}^e)|
\end{equation}
and similarly upon translation by $x$.\\
Let us define
\begin{align*}
\cD_x(m)&:= \{ \, | m -\mu^{\l_0(m)}(m
_{x,s'}^e)|\geq(s')^{-\frac{d}{2}+\d}\,\}\\
\cD'_x(m)&:=\{\,|\mu^{\l_0(m)}\bigl(\t
_x(m_s^{1,e}-m_s^{2,e})\bigr)|\geq s^{-\frac{d}{2}+\d}\,\}\\
\cD&:=\cup_m\cup
_{x\in\toromi}\bigl(\cD_x(m)\cup\cD'_x(m)\bigr)
\end{align*}
where, in the last formula, $m$ varies among all possible
values of $m_{s'}^e$.\\
$\PP(\cD_x(m))$ and $\PP(\cD'_x(m))$ can now be estimated from above by
$e^{-cs^{2\d}}$ thanks to lemma $\ref{upupa}$ applied to $f(\a)=\mu^{\l
  _0(m)}(\h_0) -m$ and
$f(\a)=\mu^{\l_0(m)}(\h_0;\h_0)-\EE\bigl[\mu^{\l_0(m)}(\h_0;\h_0)\bigr]$
respectively. Therefore, $\PP(\cD)\leq \ep^{-2d}e^{-cs^{2\d}}$ and a
simple use of Borel-Cantelli lemma proves the thesis.
\end{proof}

\section{Central Limit Theorem Variance}
\label{CLTV}
In this section we investigate the structure of the space $\cG$ that we
recall was defined as (see (\ref{vialattea}))
\begin{equation*}
\cG:=\{g\in\GG\,:\,\exists \,\L\in\FF \; \text{ such that, }
\;\forall \a \text{ and } \forall \nu\in \cM^\a(\L)\,, \;\nu(g)=0\, \}
\end{equation*}
endowed with the non negative  semi-inner product
\begin{equation}\label{semi-inner_product}
V_m(f,g):=\lim_{\ell\uto\infty}\,V_m^{(\ell)}(f,g)
\end{equation}
where
\begin{equation*}
V^{(\ell)}_m(f,g):=(2l)^{-d}\,\EE\Bigl[
\mu^{\la_0(m)}\bigl(\,\sum_{|x|\leq \ell_1}\t_x f, (-\cL_{\L_\ell})^{-1}\sum_{|x|\leq \ell _1}\t_x g\bigr)\Bigr]\,,
\quad m\in(0,1)
\end{equation*}
with $\ell_1:=\ell-\sqrt{\ell}$. For $m=0,1$ we simply define $V^{(m)}(f,g)=
V_\ell^{(m)}(f,g)=0$

\smallno
In all what follows we fix a density $m\in(0,1)$ that, most of the
times, will not appear inside the notation and we denote by $\PP^\ast$
the annealed probability measure on $\tilde \Omega$ characterized by
\begin{equation*}
\PP^\ast(d\a,d\h)= \PP(d\a)\mu^{\a,\la_0(m)}(d\h).
\end{equation*}
We remark that $\PP ^\ast$ is translation invariant and we
write $\EE^\ast$ for the corresponding expectation.

\subsection{The pre-Hilbert  space \cG}
In what follows we prove that the semi--inner product $V$ is
well defined and that the subspace generated by the currents $j_{0,e}$, $e\in
\cE$, and by the fluctuations $\lstor g$, $g\in\GG$, is dense in $\cG$. To this
aim we need to generalize the standard theory (\cite{KL} and
references therein), based on closed and exact
forms, to the disordered case. The main new feature in the disordered
case is a richer structure of the space of closed
forms which requires a proper analysis.\\
We begin with a \emph{table of calculus} that can be easily checked as
in the non disordered case.
For any $f\in \cG$, $u\in\GG$ and $e\in\cE$ let
$$
t_e(f):=\sum_{x\in\ZZ^d}(x,e)\EE^\ast(\h _x, f), \quad
(f,u)_0:=\sum_{x\in\ZZ^d}\EE^\ast(\t_xf,u).
$$
\begin{lemma}
\label{aritmetica}
For any $f\in \cG$, $u\in\GG$ and $e,e'\in\cE$
\begin{align*}
& V(f,\cL u)=-(f,u)_0,\quad \quad V( \lstor u ,
\lstor u )=  \sum _{e\in \cE} \frac{1}{2}\EE^\ast\bigl (c_{0,e}(\nabla
_{0,e}\underline u)^2\bigr),\\
& V(j_{0,e},g)=-t_e(f),\;\;\quad \quad  V (j_{0,e},j_{0,e'})=
\frac{1}{2}\EE^\ast\bigl(\;c_{0,e}(\nabla_{0,e}\h_0)^2 \bigr )\delta _{e,e' },\\
&V(j_{0,e},\cL u)=-\frac{1}{2}\EE^\ast
\bigl(c_{0,e}\nabla _{0,e}\underline u\cdot\nabla_{0,e}\h_0\;\bigr).
\end{align*}
\end{lemma}
The main result of this paragraph is the following.
\begin{theorem}
\label{bingo}
\noindent
\begin{enumerate}[i)]
\item For any $f,g\in\cG$ the limit $V(f,g):=\lim_{\ell\uto\infty}
  V_\ell(f,g)$ exists, it is finite and it defines a non negative
  semi--inner product on $\cG$. In particular
  $V(f):=\lim_{\ell\uto\infty} V_\ell(f,f)$ is well defined.

\smallno
\item For any $f\in \cG$
\begin{equation}
\label{bambi}
\begin{split}
V(f)&=\sup_{a\in\RR^d}\sup_{u\in\GG}\,
\Bigl\{2V(f,\sum_{e\in\cE}a_e j_e +\cL u)-V(\sum_{e\in\cE}a_e j_e+\cL
u) \Bigr\}\\
&=\sup_{a\in\RR^d}\sup_{u\in\GG}\,
\Bigl\{\sum_{e\in\cE}2 a_e t_e (f)+2(f,u)_0-\sum_{e\in\cE}\frac{1}{2}
 \EE^\ast\bigl(c_{0,e}(a_e\nabla_e\h_0-\nabla_e\underline u)^2\bigr)\Bigr\}\,.
\end{split}
\end{equation}

\medno
\item The subspace
\begin{equation}\label{denso}
\bigl\{\,\sum_{e\in\cE}a_e j_{0,e}+\cL u\;:\;a\in\RR^d,\;u\in\GG \,\bigr\}
\end{equation}
is dense in $\cG$ endowed of the semi--inner product $V$.
\end{enumerate}
\end{theorem}
Notice that lemma \ref{aritmetica} proves that the two expressions
appearing in the r.h.s. of the first equality sign in (\ref{bambi}) are
equal.
\\
Before proving the theorem we need to introduce the notion of closed and exact
forms together with their generalization to the disordered case and
prove few preliminary results. We refer the reader to \cite{AF}
for a complete treatment.
\begin{definition}
\label{form}
A form on $\Omega $ is a family $\xi=\{\xi_b\}_{b\sset\ZZ^d}$ of
 functions  $\xi_b:\Omega\rightarrow\RR $. It is
called closed if, given $\eta \in \Omega$ and
bonds $b_1,\dots,b_n$ with $\eta=S_{b_n}\circ\dots\circ S _{b_2} \circ
S_{b_1}(\eta) $, then
\begin{equation*}
 \sum_{i=1}^n\xi_{b_i}(\eta_{i-1})=0\quad\text{where}\quad
\h_0:=\h,\; \eta_i:=S_{b_{i}}\circ\dots\circ S_{b_2}\circ
S_{b_1}(\eta)\;\; \forall i=1,\dots,n.
\end{equation*}
\end{definition}
The expression $\sum_{i=1}^n\xi_{b_i}(\eta_{i-1})$  can be
thought of as the integral of the form $\xi$ on the closed path
$\eta_0=\h$,  $\eta_1,\dots,\eta _n=\eta $. It can be proved, see
\cite{AF}, that a form on $\Omega $  is closed if and only if it
satisfies the following properties {\bf P.${\mathbf 1}$}, {\bf
P.${\mathbf 2}$} and {\bf P.${\mathbf 3}$}.

\\
{\bf P.${\mathbf 1}$}.
 Let $a,v,w\in \ZZ ^d$ with $|v | = |w| =1 $ and $ v \pm
w \not = 0$.  We set $ c= a +v +w$,  $x= a+v$, $ x'= a +w$, $b_1 = \{ a, x\}$, $ b _2 = \{ x, c \}$, $
 b _1 ' = \{ a, x ' \}$,  $ b _2  '= \{ x' , c \}$. Then
\begin{equation*}
  \xi_{b_1}\circ S_{b_2}\circ S_{b_1}+\xi _{b_2}\circ S_{b_1}+\xi_{b_1}=
\xi_{b'_1} \circ S_{b'_2}\circ S_{b'_1} +\xi_{b'_2}\circ S_{b'_1}+\xi_{b'_1}.
\end{equation*}
{\bf P.${\mathbf 2}$}.
 For any couple of bonds $b_1, b_2 \sset\ZZ^d $ such that $
b_1 \cap b_2 = \emptyset $,
\begin{equation*}
\xi _{b_2}\circ S_{b_1}+\xi_{b_1}=\xi_{b_1}\circ S_{b_2}+\xi_{b_2}.
\end{equation*}
{\bf P.${\mathbf 3}\,$}. For any bond $b\sset\ZZ^d $,
\begin{equation*}
\xi_b\circ S_b+\xi_b=0.
\end{equation*}
The above characterization allows us to generalize the definition of
closed forms to the disorder case.
\begin{definition}
A form in $L^2(\PP^\ast)$ is a
family of functions  $\xi=\{\xi_b\}_{b\sset\ZZ^d}$ with
$\xi_b\in L^2(\PP^\ast)$.\\
A form $\xi$ is  called \emph{closed} if it satisfies
properties  {\bf P}.${\mathbf 1}$, {\bf P}.${\mathbf 2}$
and {\bf P}.${\mathbf 3}$  where equalities are in
  $L^2(\PP^\ast)$.
A form  $\xi=\{\xi_b\}_{b\sset\ZZ^d}$ if called \emph{exact} is
 $\xi_b=\nabla_b\underline u$ for some  $u\in\GG$.
A form $\xi$ is called \emph{translation covariant} if
$\tau_x\xi_b=\xi_{b+x}$  for any $ x\in\ZZ^d$, $b\sset\ZZ^d$.
\end{definition}
It is easy to check that exact forms are automatically closed and
translation covariant.
Given a closed form $\xi$ in $L^2(\PP^\ast)$ the form
on $\O$ $\{\xi_b(\a,\cdot)\}_{b\sset\ZZ^d}$ is a
closed form on $\O$ for almost
any disorder configuration $\a$. \\
In what follows by a form we will always mean a form in
$L^2(\PP^\ast)$.
\begin{definition}
A family of functions $\xi=\{\xi_e\}_{e\in\cE}$, $\xi_e\in
L^2(\PP^\ast)$, is called the \emph{germ of the form}
$\xi'=\{\xi'_b\}_{b\sset\ZZ^d}$ if $\xi'_{x,x+e}=\t_x\xi_e$ for any
$x\in \bbZ^d$ and $e\in \cE$.
\end{definition}
It follows that $\xi'$ is automatically translation covariant as soon as
it is generated by a germ $\xi$.\\
Within the subset of closed and translation covariant forms we consider
the special family $\{\mathfrak{U}^e\}_{e\in \cE}$ defined by
\begin{equation*}
\mathfrak{U}^e _{x,x+ e'}(\eta):=\delta_{e,e'}(\eta_{x+e}-\eta_x)\,, \quad
\forall x\in \ZZ^d,\; e,e'\in\cE.
\end{equation*}
It is simple to check that the form $\mathfrak{U} ^e$ is not exact.
Finally, we define $\Xi_C$ as the set of germs of closed forms
and
\begin{equation*}
 \Xi_0  := \{\xi =\{ \xi _e \}_{e\in\cE }
 \, : \, \exists a \in \RR^d , u\in\GG \quad \text{with} \quad \xi_e = a_e\mathfrak{U}^e + \nabla_e
\underline u \quad \forall e\in\cE \,\}.
\end{equation*}
We remark that $ \Xi _0 \subset \Xi _C$ and that $\Xi_C$ is a closed
 subspace in $L^2(\otimes^d\PP^\ast)$. A deeper result is given by the
 following density theorem.
\begin{theorem}
\label{medina}
$\Xi _C  = \overline{\Xi}_0$ in $L^2(\otimes^d\PP^\ast)$.
\end{theorem}
\emph{Proof.}
  The proof follows closely the proof of theorem $4.14$ in appendix $3$
  of \cite{KL} with the exception of the last step. As in \cite{KL} it
  can be proved that for any $\xi \in \Xi _C$ there exists a germ $\o
  \in\Xi_C$ with the following properties:
  \begin{enumerate}[i)]
  \item $\xi-\o\in\overline{\Xi}_0$;
  \item $\o$ can written as $\o=\o_- + \o_+$ with $\o_\pm
    =\{\o_{\pm,e}\}_{e\in \cE}$, $\o_{\pm,e}(\a,\h)=\o_{\pm,e}(\a,\h_0,\h_e)$ such that
  $\forall e\in \cE$
\begin{equation}\label{follonica}
\begin{split}
&\o_{-,e} ( \alpha , \eta_0 , \eta_{2e} )-\o_{-,e} (\alpha ,
\eta_0 , \eta_{e} )=\o_{-,e} ( \alpha , \eta_{e} , \eta
_{2e} ),\\
& \o_{+,e} ( \alpha , \eta_{-e} , \eta_{e} )-\o_{+,e} (
\alpha , \eta _0 , \eta_{e} )=\o_{+,e}(\alpha,\eta_{-e}, \eta_{0})\,.
\end{split}
\end{equation}
\end{enumerate}
It remains to prove that $\o \in\overline{\Xi}_0$. Because of
(\ref{follonica}), $\forall e\in \cE$ there exists $a_{\pm,e}\in
L^2(\PP)$ such that $\o_{\pm,e} = a_{\pm,e}(\a)(\h_e-\h_0)$. Lemma
$\ref{tipologia}$ then completes the proof of the theorem.

\begin{lemma}
\label{tipologia}
Let $\o \in \Xi_C $ such that for any $e\in\cE$ there exists $a_e\in
L^2 (\PP)$ with $\o_{e}=a_e(\a)(\eta_e-\eta_0)$. Then
$\o\in\overline{\Xi}_0$.
\end{lemma}
\begin{proof}
  By subtracting $\sum_{e\in\cE}\EE(a_e)\,\mathfrak{U}^e$ from the germ
  $\o$, we can assume that $\EE(a_e)=0$ for any $e\in\cE$. In what
  follows we denote the form generated by
  the germ $\o$  by the same symbol $\o$.\\
  Given $x\in\ZZ^d$ let $\eta ^{(x)}\in \Omega$ the configuration with
  just one particle at $x$ and let $\{b_1,\dots,b_r\}$ be a sequence of
  bonds such that $ \eta ^{(x)}=S_{b_r}\circ\dots\circ S_{b_2}\circ
  S_{b_1}(\h^{(0)})$. Define
\begin{equation}\label{g_x}
g_x(\a)=\sum_{i=1}^r\o_{b_i}(\a,\eta_{i-1}) \quad\text{where} \quad
 \eta_{i} := S _{b _i } \circ \dots \circ S _{b_2}\circ
 S_{b_1}(\h^{(0)})\;\;\; \forall i=1,\dots,r.
\end{equation}
Notice that, since $\{\o_b(\a,\cdot)\}_{b\sset\ZZ^d}$ is a closed form
on $\O$ for almost any $\a$, the definition of $g_x$ does not depend on
the particular choice of the bonds $b_1,\dots,b_r$ and
the family $\{g_x\}_{x\in\ZZ^d}$ satisfies
\begin{equation*}
g_{x+e}- g _x = \tau _x a_e \quad \forall x \in \ZZ^d, \; e\in \cE\,.
\end{equation*}
Therefore, by setting $h_n:=-\sum_{x\in\La_n }g _x(\a)\eta _x$, we get
\begin{equation}
  \label{cp1}
\nabla_b h_n= \o _b \qquad \forall\, n \in \NN, \; b\in\La_n
\end{equation}
In order to conclude the proof it is enough to show that
\begin{equation*}
\lim _{n\uto\infty}\o^{(n)}_e=\o_e \quad \forall e\in \cE\quad
\text{where}\quad
\o^{(n)}_b:=\frac{1}{(2n)^d}\nabla_b\,\underline{h}_n\in\Xi _0.
\end{equation*}
By translation covariance and (\ref{cp1}) $\nabla_{0,e} \tau_x h_n = \o_e$ if $-x,-x+e
\in\L_n$. Thus, for any $e\in \cE$, we can write
\begin{equation}
\label{tokio}
\o^{(n)}_{e} = \frac{(2n+1)^{d-1}}{(2n)^{d-1}} \o_e +
\frac{1}{(2n)^d} \sum_{\substack{x\in\L_n \\
    x_e=n}}\t_{-x}\nabla_{x,x+e} h_n
+ \frac{1}{(2n)^d} \sum_{\substack{x\in\L_n-e \\ x_e=-n-1}}\t_{-x} \nabla_{x,x+e}h_n\,.
\end{equation}
and we are left with the proof that the second and third term in the
r.h.s.of (\ref{tokio}) tend to $0$ in $L^2(\PP^\ast)$. Let us consider the
second term (the third one being similar).
By Schwarz inequality and the identity
$$
\nabla_{x,x+e}h_n =-g_x(\a)(\h_{x+e}- \h_x )\quad\forall x\in\L_n \text{
  with } x_e =n ,
$$
it is enough to show that
\begin{equation}\label{annibale}
\lim_{n\uto\infty}\frac{1}{n^{d+1}}\sum_{x\in\L_n \atop x_e=n}\EE(g_x^2)=0.
\end{equation}
To this aim, given $x=(x_1,\dots,x_d)$, we choose the bonds $b_1, \dots
, b_r$ in the definition (\ref{g_x}) in such a way that $\eta _i = \eta
^{(y _i )} $ where $ y_0$ is the origin of $\ZZ ^d$, $y_r := x$ and in
general $y_0 , y_1 , \dots , y_r $ are the points encountered by moving
in $\ZZ ^d$ first from $(0, \dots , 0 ) $ to $ (x_1 , 0 , \dots , 0 ) $
in the first direction, then from $ (x_1 , 0 , \dots , 0 ) $ to $ ( x_1,
x_2, 0 , \dots , 0 ) $ in the second direction and so on until
arriving to $x$.\\
Given this choice, it is simple to verify that for any $x \in \La _n$
and $e\in \cE$ there exists $z_e\in \L_n$ and an integer $ k_e\in [0,n]$
such that
$$
 g_x^2\leq c\sum_{e\in \cE}\Bigl(\sum_{s =0}^{k_e}\t_{z_e+se}\,a_{e}\Bigr)^2.
$$
Therefore, in order to prove (\ref{annibale}), we need to show that
$$
\lim_{n\uto\infty}\sup _{k=0,1,\dots,n}\frac{1}{n^2}\EE\Bigl(\bigl(\sum
_{s=0}^k \t_{se}\,a_{e}\bigr)^2\Bigr) \quad \forall e \in\cE.
$$
To this aim, for simplicity of notation, we fix $e\in\cE$ and
we write $a_s$ in place of $\t_{se}a_e$. Moreover, for any $r\in\NN$
we set $a_{s}^{(r)}:=\EE[a_s\,|\,\a _{\La _{se,r}}] $. Since
 $a_s^{(r)}=\t_{se}a_0^{(r)}$ and  $\EE(a_s^{(r)})=0$, we have for any
 $0\le k\le n$
\begin{gather*}
 \frac{1}{n^2}\EE \Bigl ( \bigl ( \sum _{s=0} ^k  a_{s} \bigr )
^2 \Bigr ) \leq 2 \frac{1}{n^2}\EE \Bigl ( \bigl ( \sum _{s=0} ^k [ a
_{s} - a_{s}^{ (r) }] \bigr ) ^2 \Bigr ) +2\frac{1}{n^2} \EE \Bigl
(\bigl(\sum_{s=0}^k  a_{s}^{(r)}\bigr)^2\Bigr)\\
\leq 2\EE\Bigl(\bigl(a_0-a_0 ^{(r)})^2\bigr)+\frac{c(r)}{n}.
\end{gather*}
and the thesis follows.
\end{proof}
The connection between the forms and the space $\cG$ endowed with the
semi--inner product $V(f,g)$ is
clarified by next  proposition, which can be proved, following
 \cite{KL} and \cite{VY}, as explained in section $5.5$ of \cite{AF}.
\begin{prop}
\label{faber}
Given  $ f \in \gstor$ and $e\in\cE$ there exists  a function
$\phi_e\in\GG$ such that
\begin{equation}
\label{cp2}
\sup _{\xi\in\Xi_0}\The_f(\xi) \le \liminf_{\ell\uto\infty}V_\ell(f) \le
  \limsup_{\ell\uto\infty}V_\ell(f)\leq
  \sup_{\xi\in\Xi_C}\The_f(\xi)
\end{equation}
where
$$ \The_f(\xi):= \sum _{e\in\cE} 2\,\EE^\ast\bigl(
c_{0,e}\phi_e \xi _e)- \sum _{e\in\cE} \frac{1}{2}\EE^\ast
\bigl(c_e\xi_e^2\bigr).
$$
Moreover, given $a\in\RR^d$ and $u\in\GG$,
\begin{equation}\label{jago3}
 \Theta _f \bigl ( \sum _{e\in \cE}( - a_e \mathfrak{U}  ^e + \nabla _e
\underline u )\bigr ) = \sum _{e\in\cE}2a_e t_e (f)+2(f,u)_0 -
\sum _{e\in\cE}\frac{1}{2}\EE^\ast\bigl (
c_{0,e}(a_e\nabla_e\eta_0 -\nabla_e\underline u)^2 \bigr).
\end{equation}
\end{prop}
\noindent
We are finally in a position to prove theorem \ref{bingo}.
We first observe that theorem \ref{medina} proves that the inequalities
in (\ref{cp2}) are actually equalities so that $V(f)=\lim_{\ell\uto \infty}
V_\ell(f)$ exists and it is given by (\ref{bambi}).
Moreover, because of  (\ref{aglaja}), $V(f) < \infty$ so that, by polarization, $V(f,g)$
exists finite for any $f,g\in \cG$ and it defines a semi--inner product.
The density of the subspace (\ref{denso}) follows at once from the first equality in (\ref{bambi}).
\qed

\subsection{The method of long jumps revisited}
In this paragraph we consider, for any $e\in\cE$, a particular sequence
$\{W^e_n/n\}_{n\in \bbN}$ in the space $\cG$ which is asymptotically equivalent to the sequence
$$
2m(1-m)\la_0'(m)\frac{\psi_{n,n}^e}{n},\quad n\in\NN
$$
where $\psi_{n,n}^e$ has been defined in (\ref{estasi}) as
$$
\psi_{n,n}^e = m_n^{2,e}-m_n^{1,e} -\mu\bigl[m_n^{2,e}-m_n^{1,e}\tc m_n^e\bigr]
$$
The functions $W^e_n$ have been introduced in \cite{Q} in order to
depress the extra fluctuations produced by the disorder and
are defined as
\begin{equation*}
 W^e_n:=\Av_{x\in\L^{1,e}_n}\Av_{y\in\L^{2,e}_n}\,w_{x,y}\quad
\text{ where }\quad
 w_{x,y}:=\bigl(1+e^{-(\a_x-\a_y)(\h_x-\h_y)}\bigr)(\h_y-\h_x).
\end{equation*}
We remark that, for any bond $b=\{x,y\}$, the quantities
$c_{x,y}:=1+e^{-(\a_x-\a_y)(\h_x-\h_y)}$ are a possible choice of
transition rates compatible with our general assumptions (see section
\ref{model}). Therefore, for generic $x,y\in\ZZ^d$ $c_{x,y}$ can be
thought of as the rate of the (long) jump from $x$ to $y$ and
viceversa. In a sense the rates $c_{x,y}$, $x,y\in \bbZ^d$, define a new
process with arbitrarily long jumps but still reversibile w.r.t. the Gibbs
measure of the system.
\begin{remark}
  The role of the function $W_n^e$ here is very different from that
  indicated in \cite{Q}.  In our approach and for reasons that will
  appear clearly in the next subsection, we are interested in computing
  the asymptotic of the semi--inner product
  $V(j_{0,e'},\frac{\psi_{n,n}^e}{n})$ as $n\uto\infty$.  Our strategy
  to compute $V(j_{0,e'},\frac{\psi_{n,n}^e}{n})$ is to replace (in
  $\cG$) $\frac{\psi^e_{n,n}}{n}$ with $\frac{W^e_n}{n}$ and then to exploit some nice
  integration by parts properties pointed out in \cite{Q} (see below).\\
  In \cite{Q} instead, the main idea is first to approximate, as $\ep\dto
  0$, the microscopic current $j_{0,e}$ with a fluctuation term $\cL g$
  plus a linear combination of the $\frac{W^{e'}_{k}}{k}$, $e'\in \cE$, on a scale
  $k$ that must diverge as $\ep\dto 0$ like $\ep^{-\frac{2}{d+2}}$. The second step
  indicated in \cite{Q} is to replace $\frac{W^{e}_{k}}{k}$ with
$$
2m(1-m)\l'_0(m)\frac{(m_k^{2,e}-m_k^{1,e})}{k}
$$
Such a step is very similar to the main result of this subsection
described at the beginning but, at the same time, very different. The
first main difference is that our mesoscopic scale $n$ is not linked
with $\ep$. The second difference is that our functions $\psi^e_{n,n}$
represent (discrete) gradient of the density \emph{minus} their
canonical average. Such a counter term, discussed at length in section
\ref{Fluct}, is absent in the approach of \cite{Q}.
\end{remark}
Our main result is given by the following theorem
\begin{theorem}
\label{longjumps}
For any $e\in \cE$
\begin{equation}
\lim _{n \uparrow
\infty}V\Bigl(\frac{W^e_n}{n}-2m(1-m)\la_0'(m)\frac{\psi_{n,n}^e}{n}\Bigr)=0.
\end{equation}
\end{theorem}
We will use the above result only to compute the limit of
$V(j_{0,e'},\frac{\psi_{n,n}^e}{n})$ .  Indeed, as pointed in \cite{Q},
the function $w_{x,y}$ satisfies the following integration by parts
property: for any $\L\in\FF$ with $\L\ni x,y$ and any $\nu\in\cM(\L)$
\begin{equation*}
\nu(w_{x,y}g)=\nu((\h_x-\h_y)\nabla_{x,y}g).
\end{equation*}
By the above property and  lemma \ref{aritmetica} it is simple to
check that, for any $e,e'\in\cE$,
$V\Bigl(j_{0,e'},\frac{W^e_n}{n}\Bigr)=-2m(1-m)\d_{e,e'}$.
Therefore, by theorem \ref{longjumps}, we get
\begin{equation}
 \label{macchiavelli}
\lim _{n\uto\infty} V\Bigl(j_{0,e'},\frac{\psi^e_{n,n}}{n}\Bigr)=-\chi(m)\d_{e,e'},\quad
\forall e,e'\in\cE.
\end{equation}
\begin{proof}
In order to prove theorem \ref{longjumps}  it is convenient to introduce some
notation. First, we fix the vector $e\in\cE$ which will be often omit in
the notation.  Moreover we introduce the following
 equivalence relation.
\begin{definition}
\label{def-approx}
 Given two sequences of functions $\{f_n\}_{n\in \bbN}$ and $\{g_n\}_{n\in \bbN}$
 such that $f_n$ and $g_n$ have support in $\L^e_n$, we write $f_n\approx g_n$ if
\begin{equation*}
\lim_{n\uto\infty}V\Bigl (\frac{f_n-\mu[f_n\,|\,m^e_n]}{n} -\frac{g_n-\mu[g_n\,|\,m^e_n]}{n}\Bigr)=0.
\end{equation*}
\end{definition}
\noindent
$\bullet$ {\bf Step 1}: $f_n:=
W_n-\mu[W_n|\,m^{1,e}_n,m^{2,e}_n]\approx 0$.

\smallno
For any $x\in\ZZ^d$, let  $\nu _{x,n}$ be the random
 canonical  measure  $\mu[\cdot \tc\cF_{x,n}]$ where
 $\cF_{x,n}$ is the $\sigma$--algebra generated by
$\t_x m^{1,e}_n$, $\t _x m^{2,e}_n$ and $\h_y$ with $y\not
\in\L_{x,n}^e$.
Let us observe that
\begin{enumerate}[i)]
\item
 $\mu^{\l_0(m)}(\tau _x f_n,g)= \mu^{\l_0(m)}\bigl(\,\nu_{x,n}(\t_x W_n;g)\,\bigr)$
for any function $g$;

\smallno
\item
$W_n$ can be written as sum of functions $f$  of the following  form
$$
f=\Av_{z\in\L_n^{1,e}}\t_z h\, \Av_{z'\in\L_n^{2,e}} \t_{z'} h'
$$
where $h$ and $h'$ depend only on $\a_0$ and $\h_0$.
\end{enumerate}
Because of $i)$ and $ii)$ and thanks to the
the variational characterization (\ref{maria}) of
$V_\ell(\cdot,\mu^{\la_0(m)})$,
it is enough to prove that,
for a function $f$ as in $ii)$,
\begin{equation} \label{falco}
 \lim_{n\uto\infty}\lim_{k\uto\infty}\frac{1}{n^2 k^d}\,
\EE\bigl[\,\sup_{g\in \GG}\,\{\phi(g)/\cD_{\L_k}(g;\mu^{\l_0(m)})\}\,\bigr]=0
\end{equation}
where
\begin{equation*}
\phi(g):=\bigl[\sum_{|x|\leq
  k_1}\mu^{\l_0(m)}\bigl(\nu_{x,n}(\t_xf;g)\bigr)\bigr]^2,\quad k_1:= k-\sqrt{k}.
\end{equation*}
By proposition \ref{special_covariances}, for any $\d>0$  there exists
 $\ell_0\in\NN$  such that, if $n\geq \ell\geq \ell_0$, then
\begin{equation}
  \label{covarianze1}
\nu_{x,n}(\t_x f;g)^2\leq
\frac{c(\ell)}{n^d}\cD(g;\nu_{x,n})+\frac{\d}{n^d}\Var_{\nu_{x,n}}(g)+
\frac{c}{n ^d}\Var_{\nu_{x,n}}(g)\th_{\L_{x,n}^e,\ell}(\a)
\end{equation}
where, for any given $\g>0$ and $\ell\geq \ell_1(\g)\geq \ell_0$,
\begin{equation}
\label{ilsecolobreve}
\PP(\th_{\L_{x,n}^e,\ell}(\a)\geq \g)\leq e^{-c(\g,\ell)n^d}
\end{equation}
for a suitable constant $c(\g,\ell)$. Using the spectral gap estimate
(\ref{eq:gap}), the r.h.s. of  (\ref{covarianze1}) can be bounded by
$$
\cD(g;\nu_{x,n})\bigl(c(\ell)+c\,\d\, n^2+c\,n^2\th_{\L_{x,n}^e,\ell}(\a)\bigr)n^{-d}
$$
and therefore, by Schwarz inequality,
$$
\phi(g)\leq\cD_{\L_k}(g;\mu^{\l_0(m)})
\bigl(c(\ell)k^d+c\,\d k^d n^2+c\,n^2 \sum_{|x|\leq k_1}\th_{\L_{x,n}^e,
\ell}\bigr).
$$
By taking the limits $\d\dto 0, \,\ell\uto\infty,\,n\uto\infty,\,k\uto\infty$ (from
right to left), in order to prove (\ref{falco}) the thesis follows since
$\lim_{n\uto\infty}\EE\bigl(\th_{\L_n^e,\ell}(\a)\bigr)=0$
because of (\ref{ilsecolobreve}).

\\
$\bullet$ {\bf Step 2}: $\mu[W_n|\,m^{1,e}_n,m^{2,e}_n]\approx
2m(1-m)\la_0'(m)\psi_{n,n}^e.$

\smallno
The proof is based on the following lemma, which follows easily from
the variational characterization  of $V_\ell(\cdot, \mu^{\la_0(m)})$ given in
(\ref{maria}).
\begin{lemma}  
\label{criterio}
Let, for any  $n\in\NN$, $f_n,\,h_n\in\GG$ be such that
\begin{enumerate}[i)]
\item $\D_{f_n}\sset\L^e_n$\,;

\smallno
\item $\sup _n \|h_n\|_\infty <\infty  \quad\text{and}
\quad  \lim _{n\uto\infty}n^d\EE\bigl[\mu^{\l_0(m)}\bigl(h_n^2\bigr)\bigr]=0$;

\smallno
\item $ |f_n|\leq |h_n|$.
\end{enumerate}
Then   $f_n\approx 0$.
\end{lemma}
Thanks to the estimates given in the Appendix it can be proved
(see \cite{AF}) that condition \emph{ii)} of the lemma is
satisfied  by any
of the  following sequences:
$$
\{n^{-d}\}_{n\in\NN},\quad
\bigl \{\,(m-m_{\D_n})^2\,\}_{n\in\NN},
\quad
\bigl\{ \id_{\{|m-m_{\D_n}|\geq c\}}\bigr\}_{n\in\NN}, \quad
\bigl\{ \bigl(m-\mu^{\l_0(m)}(m_n^{i,e})\bigr)^2\bigr \}_{n\in\NN}
$$
where $i=1,2$, $c>0$ and $\D_n$ is either one of the sets $\L^e_n$, $\L _n^{1,e}$, $\L _n^{2,e}$.

\smallno
As in \cite{Q}  we define the (random w.r.t. $\a$) function $F_n(m_1,m_2)$ as
\begin{equation*}
  F_n(m_1,m_2) = \mu_{\L _n^{1,e}}^{\l(m_1)}\otimes \mu_{\L _n^{2,e}}^{\l(m_2)}\bigl(W_n\bigr)
\end{equation*}
It is not difficult to show that $F_n(m_1,m_2)$ has the explicit expression
\begin{equation*}
 F_n(m_1,m_2)=m_1-m_2+e^{\la_{1,n}(m_1)-\la_{2,n}(m_2)}(1-m_1)m_2-
e^{\la_{2,n}(m_2)-\la_{1,n}(m_1)}m_1(1-m_2).
\end{equation*}
The main reason to introduce $F_n(m_1,m_2)$ is that
\begin{equation}
\label{step2.1}
\mu[W_n|\,m^{1,e}_n,m^{2,e}_n]\approx  F_n(m^{1,e}_n,m^{2,e}_n).
\end{equation}
This equivalence follows at once from the equivalence of the ensembles
together with lemma \ref{criterio} applied to
$f_n=\mu[W_n|\,m^{1,e}_n,m^{2,e}_n]- F_n(m^{1,e}_n,m^{2,e}_n)$ and $h_n=
c\,n^{-d}$ for a large enough constant $c$.
\\
Next, again by lemma \ref{criterio} applied with $h_n =
\id_{\{|m-m_n^{1,e}|\geq c_m\}} + \id_{\{|m-m_n^{2,e}|\geq c_m\}}$, \hbox{$c_m=(m\mmin(1-m))/2$}, we get that
\begin{equation}
\label{step2.2}
  F_n(m^{1,e}_n,m^{2,e}_n) \approx
  F_n(m^{1,e}_n,m^{2,e}_n)\id_m
\end{equation}
where $\id_m := \id_{\{|m-m_n^{1,e}|\le c_m\}}\id_{\{|m-m_n^{2,e}|\le
  c_m\}}$.

\smallno
Next, by Taylor expansion around the arithmetic mean of $m_n^{1,e}$ and $m_n^{2,e}$, we write
\begin{gather*}
  F_n(m^{1,e}_n,m^{2,e}_n) = \\  F_n(m^e_n,m^e_n) +
\frac{\partial F_n}{\partial m_1}(m^e_n,m^e_n)(m_n^{1,e}-m^e_n)+
\frac{\partial F_n}{\partial m_1}(m^e_n,m^e_n)(m_n^{2,e}-m^e_n)
+ R_n(m^{1,e}_n,m^{2,e}_n)
\end{gather*}
Then, the zero order contribution $F_n(m^e_n,m^e_n)\id_m$ is negligible,
$F_n(m^e_n,m^e_n)\id_m \approx 0$, since $F_n(m^e_n,m^e_n)
\approx 0$ because of definition
\ref{def-approx} and $F_n(m^e_n,m^e_n)(1-\id_m) \approx 0$ again by
lemma \ref{criterio}.
\\
The second order error term, $R_n(m^{1,e}_n,m^{2,e}_n)\id_m$, is
negligible because of lemma \ref{criterio} applied with $h_n=
c\bigl[(m_n^{1,e}-m^e_n)^2 + (m_n^{2,e}-m^e_n)^2\bigr]$. Notice that it
is here that the characteristic function $\id_m$ plays an important
role since the second derivatives of $F_n(m_1,m_2)$ diverge as $m_i$ tends
to $0$ or to $1$.
\\
Let us now examine the relevant first order terms.
We claim that for $i=1,2$
\begin{equation}
  \label{statement1}
 \frac{\partial F_n}{\partial m_i}(m^e_n,m^e_n)(m_n^{i,e}-m^e_n)\id_m +
(-1)^{i}2m^e_n(1-m^e_n)\la_{i,n}'(m^e_n)(m_n^{i,e}-m^e_n) \id_m
\approx 0
\end{equation}
and
\begin{equation}
  \label{statement2}
2m^e_n(1-m^e_n)(\la_{i,n}'(m^e_n)-\la_0'(m))(m_n^{i,e}-m^e_n)\id_m\approx 0.
\end{equation}
where $\la _{i,n}:=\la_{\L^{i,e}_n}$.\\
Before proving (\ref{statement1}) and (\ref{statement2}) let us
summarize what we have obtained so far.
Thanks to (\ref{step2.1}), (\ref{step2.2}), the above discussion of the
Taylor expansion and (\ref{statement1}) together with (\ref{statement2})
\begin{equation*}
  \mu[W_n|\,m^{1,e}_n,m^{2,e}_n] \approx 2 m^e_n(1-m^e_n) \la'_0(m)(m_n^{2,e}-m_n^{1,e})\id_m
\end{equation*}
Using once more lemma \ref{criterio} it is now rather simple to remove the
 factor $\id_m$ and to replace $m_n^e$ with $m$, thus
concluding the proof.
\\
We are left with the proof of (\ref{statement1}) and (\ref{statement2}).

\smallno
Let us prove (\ref{statement1}) for $i=1$. By computing
$ \frac{\partial F_n}{\partial m_1}$  it is simple to check that the
 l.h.s. of  (\ref{statement1}) is equal to
\begin{equation}
 \label{camaleonte}
\begin{split}
&(e^{\la_{1,n}(m^e_n)-\la_{2,n}(m^e_n)}-1)m^e_n (\la_{1,n}'
(m^e_n)(1-m^e_n )-1 )(m_n^{1,e}-m^e_n)\id _m \, +\\
& (e^{\la_{2}(m^e_n)-\la_{1,n}(m^e_n)}-1)
(1-m^e_n) (\la'_{1,n}(m^e_n)m^e_n-1)(m_n^{1,e}-m^e_n)\id_m.
\end{split}
\end{equation}
It is enough to show that both addenda in (\ref{camaleonte}) are
equivalent to $0$ and for simplicity we deal with only with the first one.
Since $\sup _{n} \|\la_{1,n}'(m^e_n )\id_m \|_\infty \leq k_m$ for a
suitable constant $k_m$ depending on $m$, using the estimate
$|e^z-1|\leq e^{|z|}|z|$ valid  for any  $z \in \RR$ and thanks to lemma
\ref{vespe} we obtain
\begin{equation}
\label{circe}
\begin{split}
|\text{ first term in (\ref{camaleonte}) }|& \leq k_m\,  |\la_{1,n}(m^e_n)-\la_{2,n}(m^e_n)|\,|m_n^{(i)}-m^e_n|\id_m\\
&\le k'_m\Bigl(\,\sum_{i=1,2}\,\bigl(m^e_n-\mu^{\la_0(m^e_n)}(m_n^{i,e})\bigr)^2
+  \,(m_n^{1,e}-m^e_n)^2\,\Bigr)\,.
\end{split}
\end{equation}
The claim follows by applying lemma
\ref{criterio} with $h_n$ equal to the r.h.s. of (\ref{circe}).

\smallno
Let us prove (\ref{statement2}). By Schwarz inequality, it is enough to
apply lemma \ref{criterio} with
$h_n:=(\la_{i,n}'(m^e_n)-\la_0'(m))^2 \id _m + (m_n^{i,e}-m^e_n)^2$. In
order to verify condition $ii)$ of lemma \ref{criterio} for $h_n$, thanks to
the boundedness of
$ (\la_{i,n}'(m^e_n)-\la_0'(m))^2 \id _m$ uniformly in $n$,
we only need to prove that
\begin{equation*}
\lim _{n\uparrow \infty}n^d \EE\bigl[\mu^{\la_0(m)}\bigl
( (\la_{i,n}'(m^e_n)-\la_0'(m))^4 \id_m \bigr)\bigr]=0
\end{equation*}
or equivalently
\begin{equation}
\label{madrenotte}
\lim _{n\uto \infty}n^d \EE\Bigl[\mu^{\la_0(m)}\Bigl(\bigl\{
\Av_{x\in \L_n^{i,e}}\bigl[\,\mu^{\l_{i,n}(m_n^e)} (\h_x;\h_x)
-\EE\,\mu^{\la_0(m)}(\h_0;\h_0)\,\bigr]\,\bigr\}^4 \Bigl)\Bigr]=0
\end{equation}
Let $g_x(\l):=\mu^\l (\h_x;\h_x)$ and observe that
l.h.s. of (\ref{madrenotte}) is bounded from above by
\begin{equation}
 \label{rosa}
c\,\lim_{n\uto \infty}n^d\EE\Bigl[\mu^{\la_0(m)}\bigl(A_n^{(1)}+A_n^{(2)}+A_n^{(3)}\bigr)\Bigr]
\end{equation}
where
\begin{align*}
A_n^{(1)}&=\bigl\{\Av_{x\in \L_n^{i,e}} [\,g_x(\la_{i,n}(m^e_n))-
g_x(\la_{i,n}(m))\,]\bigr\}^4 \\
A_n^{(2)}&=\bigl\{\Av_{x\in \L_n^{i,e}} [\,g_x(\la_{i,n}(m))-
g_x(\la_0(m))\,]\,\bigr\}^4 \\
A_n^{(3)}&=\bigl\{\Av_{x\in \L_n^{i,e}} [\,g_x(\la_0(m))-\EE\,\mu^{\la_{0}(m)}(\h_0;\h_0)\,]\bigr\}^4.
\end{align*}
By lemma $\ref{vespe}$, $A_n^{(1)}\leq c\,(m^e_n-m)^4$ and $A_n^{(2)}
\leq c(m-\mu^{\la_0(m)}(m_n^{i,e}))^4$. At this point
(\ref{rosa}) follows by simple considerations for sum of centered independent
random variables.

\end{proof}

\subsection{The subspace orthogonal to the fluctuations}
Here we introduce a convenient Hilbert space $\cH$ containing $\cG$ and
we describe the orthogonal subspace in $\cH$ of the space of
fluctuations $\{\cL g\,:\,g\in\GG\}$.
\begin{definition}
Let $\cN:=\{g\in \gstor\,:\, V(g)=0\}$ and let $\cH$ be the
 completion of the  pre-Hilbert space $\cG/\cN$. With an abuse of notation, we
 write $V$ for the scalar product in $\cH$ induced by the semi-inner
 product $V$ in $\cG$.
\end{definition}
The sets
\begin{equation*}
\cL\GG:=\{\cL g\,:\,g\in\GG\},\quad\quad  \cL\cG:=\{\cL g\,:\,g\in\cG\}
\end{equation*}
can be considered as subsets of $\cH$ in a natural way.
Our main result proves that for any $e\in\cE$ the sequence
$\{\psi^e_{n,n}/n\}_{n\in\NN}$  converges in $\cH$ to some limit
point $\psi_e$ and that the set $\{\psi_e\}_{e\in\cE}$ forms a basis of
$\cL\GG^\perp$. The Cauchy property of the
sequence $\{\psi^e_{n,n}/n\}_{n\in\NN}$ follows by a telescopic estimate based on the
variance bounds discussed in subsection \ref{Var-Bounds}. To this aim
the following lemma is crucial.
\begin{lemma}\label{florence}
Given $k\in \bbN$ let $f\in\cG$  be such that $\D_f\sset\L_k$. Then
\begin{equation*}
V(f)\leq
c\,k^{d+2}\EE\bigl(\Var_{\mu^{\l_0(m)}}(\Av_{x\in\La_k}\t_xf)\bigl).
\end{equation*}
\end{lemma}
\begin{proof}
We first estimate $V_\ell(f)$ for $\ell\gg 1$ by means of  lemma
\ref{miramare}. To this aim we partition the cube  $\L_{\ell_1}$ into
non overlapping cubes $\{\L_{x_i,k}\}_{i\in I}$ of side
$2k+1$ and write
$$
\Av_{x\in\L_{\ell_1}}\t_x f=
\Av_{i\in I}\bigl(\Av_{x\in\L_{x_i,k}}\t_x f\bigr)
$$
Therefore, by applying lemma
\ref{miramare} with $\L=\L_{\ell}$ and $\L_i=\L_{x_i,2k}$, we obtain
\begin{equation*}
V_\ell(f)\leq c\,k^{d+2}\Av_{i\in I}\Var_{\mu^{\l_0(m)}}
\bigl(\Av_{x\in\L_{x_i,k}}\t_x f\bigr).
\end{equation*}
It is enough now to take the expectation w.r.t. $\a$ and then the limit
$\ell\uto\infty$.
\end{proof}
Lemma \ref{florence} and proposition \ref{bound:variance} allow us to prove
the key technical estimate of this subsection:
\begin{lemma}\label{baltimora}               
Let $d\geq 2$, $ n\leq s\leq k\leq 100s$ be positive integers and $0<\d\ll
1$. Then
\begin{equation}
\label{davide}
V(\phi^e_{n,s}-\phi^e_{n,k})\leq c\,s^{2-d+\d} \quad \forall e\in\cE
\end{equation}
for any $s$ large enough ($s\geq s_0(\d)$).
\end{lemma}
\begin{proof}
Since $\phi^e_{n,s}-\phi^e_{n,k}\in\cG$ has support in $\L_k$, by lemma
\ref{florence} we obtain
\begin{equation*}
V(\phi^e_{n,s}-\phi^e_{n,k})\leq c\,
k^{d+2}\sum_{r=s,k}\EE\bigl(\Var_{\mu^{\l_0(m)}}(\Av_{x\in\La_k}\t_x
\phi^e_{n,r})\bigl).
\end{equation*}
The thesis now follows from proposition \ref{bound:variance}.
\end{proof}
We also need a density result.
\begin{lemma}\label{budino}  
$\cL\cG$ and $\cL\GG$ have the same closure in $\cH$.
\end{lemma}
\begin{proof}
We fix  $g\in\GG$ and we prove that $\cL g=\lim_{s\uto\infty}\cL
(g-g_s)$ where $g_s = \mu[g\,|\,m_s]$, i.e. that $\lim_{s\uto\infty}V(\cL g_s)=0$.
To this aim we define $X_s:=\{x:s-1\leq |x|\leq s+1\}$. Then lemma
$\ref{aritmetica}$ implies that
\begin{equation}\label{coc12}
V(\cL g_s) =\sum _{e\in\cE}
\frac{1}{2}\EE^\ast\bigl(c_{0,e}(\sum_{x\in X_s}\nabla_{0,e}\tau _xg_s)
 ^2\bigr).
\end{equation}
Let $\hat{g}_s (\a,\h):=\mu_{\L_s}^{\l(m_s(\h))}(g).$ By the equivalence
of ensembles (see lemma \ref{equi}), in $(\ref{coc12})$
$g_s$ can be substituted by $\hat{g}_s$ with an error bounded by
$c\,s^{-2}$.  By lemma$\ref{vespe}$,
 $|\nabla_{0,e}\t_x\hat{g}_s|\leq c\,s^{-d}$  which, thanks to
$ (\ref{coc12})$ with $g_s$ replaced
with $\hat g _s $, implies that $V(\cL g_s)\leq c\,s ^{-2}$.
\end{proof}
We are ready for the first result about the structure of the space $\cL\GG^{\perp}$.
\begin{prop}  
\label{demetra}
Let $d\geq3$ and $e\in\cE$. Then the sequence
\begin{equation*}
\psi^e_{1,s}=\h_e-\h_0-\mu[\h_e-\h_0\,|\,m^e_s]
\end{equation*}
converges to some element $\psi_e\in \cL\GG^{\perp}$ as $s\uto\infty$.\\
Moreover,
\begin{equation}
  \label{superlimit}
\lim_{s\uto\infty}\frac{\psi^e_{n,s}}{n}=\psi_e \quad \forall n \in \NN.
\end{equation}
\end{prop}
\begin{proof}
We fix $0<\d\ll 1$. By lemma \ref{baltimora}, if $i\in\NN$ is large
 enough and  $i^3\leq s\leq(i+1)^3$,
\begin{equation*}
V(\psi^e_{1,i^3}-\psi^e_{1,s})\leq c\,i^{\,3(2-d+\d)}.
\end{equation*}
Since $d\geq 3$,  it is enough to prove that the sequence
$\{\psi^e_{1,i^3}\}_{i\in\NN}$ is Cauchy. This follows by applying
 again lemma \ref{baltimora} to get
\begin{equation*}
\sum_{i=1}^\infty V^{\frac{1}{2}}(\psi^e_{1,i^3}-\psi^e_{1,(i+1)^3})\leq
\sum_{i=1}^\infty c\, i^{\,\frac{3}{2}(2-d+\d)}<\infty.
\end{equation*}
Next we prove that $\psi_e$, the limit point of
$\{\psi^e_{1,s}\}_{s\in\NN}$, belongs to $\cL\GG^\perp$. To this aim, by
lemmas \ref{aritmetica} and \ref{budino}, we need to show that
\begin{equation*}
\lim _{s\uto \infty} \sum_{x\in\ZZ^d}
\EE[\,\mu^{\la_0(m)} (\psi^e_{1,s},\t_x g)\,]=0 \quad \forall g\in\cG,
\end{equation*}
or similarly (by translation invariance of the random field $\a$)
\begin{equation*}
\lim _{s\uto \infty} \sum_{x\in\ZZ^d}
\EE[\,\mu^{\la_0(m)} (\phi^e_{1,s},\t_x g)\,]=0 \quad \forall g\in\cG,
\end{equation*}
where we recall $\phi^e_{1,s} = \mu[\h_e-\h_0\,|\,m^e_s]$.
To this aim we set
\begin{equation*}
\D_s:=\{x\in\ZZ^d\,:\,(x+\D_g)\cap\L^e_s
\neq\emptyset\;\text{and}\;(x+\D_g)\cap(\L^e_s)^c\neq\emptyset\,\}.
\end{equation*}
Since $g\in\cG$,
\begin{equation}
 \label{temporale}
 \sum_{x\in\ZZ^d}\EE [\,\mu^{\la_0(m)} (\phi^e_{1,s},\t_x g)\,]=
\EE[\,\mu^{\la_0(m)}(\phi^e_{1,s}, \sum_{x\in\D_s}\t_x g)\,].
\end{equation}
We estimate the r.h.s. of (\ref{temporale}) by Schwarz
inequality. Let us  observe that
\begin{equation}
 \label{firstvariance}
\EE\,[\Var_{\mu^{\la_0(m)}}\bigl( \sum_{x\in\D_s}\t_x g\bigr)\,]\leq
 c_g\, s^{-d+1}.
\end{equation}
for some finite constant $c_g$. Therefore, in order to conclude the proof, it is enough to show
\begin{equation}
  \label{secondvariance}
\EE[\,\Var_{\mu^{\la_0(m)}}(\phi^e_{1,s})\,]\leq c\,s^{-d}.
\end{equation}
By the equivalence of ensemble (see lemma \ref{cleonimo})
and Poincar\'e inequality for Glauber dynamics, we obtain
\begin{equation}
 \label{28giugno1914}
\Var_{\mu^{\la_0(m)}}(\phi^e_{1,s})\leq c\,s^{-2d}+
c\,\Var_{\mu^{\la_0(m)}}(\hat \phi^e_{1,s}) \leq
c\,s^{-2d}+c\,s^d \mu^{\la_0(m)}\bigl((\nabla_0\hat \phi^e_{1,s})^2 \bigr)
\end{equation}
where $\hat \phi^e_{1,s}$ has been defined in (\ref{definition1}).\\
By lemma \ref{vespe} the last term in (\ref{28giugno1914}) is bounded by
$c\,s^{-d}$ thus proving (\ref{secondvariance}).\\
Finally we prove (\ref{superlimit}). To this aim, by writing
\begin{equation*}
\psi^e_{n,s}=\frac{1}{n}\sum _{v=0} ^{n-1} \Av _{x\in\L^{1,e}_n}
\bigl(\h_{x+(v+1)e}-\h_{x+ve}- \mu [\h_{x+(v+1)e}-\h_{x+ve}\,|\,m_s^e]
\, \bigr),
\end{equation*}
and by the observation that $\t_xf =f$ for any $f\in \cH$ and $x\in \bbZ^d$,
it is enough to prove that for any given $x\in\ZZ^d$
\begin{equation}
 \label{elke}
V(\mu[ \h_e-\h_0\,|\,m_{x,s}^e]-\mu[ \h_e-\h_0\,|\,m_s^e])
\end{equation}
goes to $0$ as $s\uto\infty$.
As in the proof of lemma \ref{baltimora} (\ref{elke}) is bounded from
above by $c(\delta)s^{2-d+\d}$ for any $0<\d\ll 1$.
\end{proof}
We are now able to exhibit a basis of $\cL\GG^\perp$ related to the
functions $\frac{\psi^e_{n,n}}{n}$ with $n\in\NN$ and $e\in\cE$.
\begin{theorem}
  \label{ortogonale} 
Let $d\geq 3$. Then
\begin{equation}
\label{policastro}
 \lim_{n\uto \infty}\frac{\psi^e_{n,n}}{n}=\psi_e \quad \forall e \in \cE
\end{equation}
where $\psi_e$ is as in proposition \ref{demetra}. Moreover,
\begin{equation}
 \label{aritmetica2}
 V(j_{0,e'},\psi_e)=- \chi(m) \d_{e',e}\quad \forall e,e' \in\cE
\end{equation}
and   $\{\psi_e\}_{e\in\cE}$ forms a
 basis of $\cL \GG^\perp$.
\end{theorem}
\begin{proof}
For any $n\in\NN$ let $k_n\in\NN$ be such that $(k_n-1)^3<n\leq k_n^3$.
Then, by lemma  \ref{baltimora}, $V(\psi^e_{n,n}/n -\psi^e_{n,k^3_n}/n)
\dto 0$ as $n\uto \infty$. Therefore, thanks to (\ref{superlimit}),
\begin{equation}
\label{gabbiano}
\lim_{n\uto \infty} V^{\ov2}(\psi_e-\frac{\psi^e_{n,n}}{n}) = \lim_{n\uto \infty}
V^{\ov2}(\psi_e - \frac{\psi^e_{n,k_n}}{n}) \le
\lim_{n\uto \infty} \frac{1}{n} \sum _{i=1}^\infty V ^\frac{1}{2}
 (\psi^e_{n,i^3}-\psi^e_{n,(i+1)^3} )
\end{equation}
and the last series is converging by lemma \ref{baltimora}.
Thus (\ref{policastro}) follows.\\
At this point, (\ref{aritmetica2}) follows from theorem
\ref{longjumps}. Let us prove that $\{\psi_e\}_{e\in\cE}$ forms a basis of
$\cL\GG^\perp$. Let $P$ be the
orthogonal projection of $\cH$ on $\cL\GG^\perp$. Then, $\cL\GG^\perp$
has dimension non larger than $d$ since, by theorem
\ref{bingo}, it is generated by  $\{Pj_{0,e}\}_{e\in\cE}$.
By (\ref{aritmetica2}) $\{\psi_e\}_{\in\cE}$ is a set of $d$  independent vectors
belonging to $\cL\GG^\perp$ and therefore a basis of $\cL\GG^\perp$.
\end{proof}
\begin{remark}
\label{omogeneo}
Let us make an observation which will reveal useful in the proof of the continuity of the diffusion matrix $D(m)$ (see next subsection). \\
Since the constant $c$ appering in (\ref{davide}) does not depend on the density $m$ and thanks to the estimate (\ref{gabbiano}), the statement (\ref{policastro}) in the above theorem can be strengthed as
$$ \lim_{n\uto\infty} \sup_{m\in (0,1)} V_m\bigl( \frac{\psi^e_{n,n}}{n}-\psi_e
\bigr) =0 \quad \forall e\in\cE.
$$
\end{remark}

\subsection{Decomposition of currents}\label{biodegradabile}
In this subsection we prove  the characterization and the
regularity of the diffusion matrix $D(m)$ stated in theorem \ref{matrix} and we prove   also theorem  \ref{agognatameta}, which is  crucial for the estimate of $\O_0$  (see subsection \ref{pietro}). In what follows, we assume $d\geq 3$.

\smallno
Denoting by  $P$  the orthogonal projection
of  $\cH$ on  $\cL\GG^\perp$, thanks to  theorem \ref{ortogonale},
for a suitable  $d\times d$ matrix $D(m)$  we can write
\begin{equation}
 \label{decomposition}
j_{0,e}=-\sum_{e'\in\cE}D_{e,e'}(m)\psi_{e'}+ (1-P)(j_{0,e})\quad \quad
\forall e\in \cE.
\end{equation}
By taking  the scalar product of both sides of (\ref{decomposition})
with  $j_{0,e'}$,  thanks to lemma \ref{aritmetica} and   (\ref{aritmetica2}),
 we obtain
\begin{equation*}
  D_{e,e'}(m)= \frac{1}{\chi(m)}V_m( Pj_{0,e}, Pj_{0,e'} )
\end{equation*}
thus proving that $D(m)$ is a non-negative symmetric matrix. In
particular, $D(m)$ can be characterized as the unique symmetric
$d\times d$ matrix such that
\begin{equation}
 \label{scozia}
(a, D(m) a) =\frac{1}{\chi (m)}V_m\bigl ( P(\sum _{e\in\cE}a_ej_{0,e})
  \bigr) \quad \forall a\in \RR^d.
\end{equation}
Since the r.h.s. of (\ref{scozia}) can be written as
$$ \inf _{g\in\GG}\frac{1}{\chi (m)} V_m \bigl (\sum
_{e\in\cE}a_ej_{0,e} -\lstor g \bigr ), $$
by lemma \ref{aritmetica} the matrix $D(m)$ corresponds to
the one described in proposition \ref{matrix}.

\smallno
In the following lemmas  we describe some properties of the diffusion
matrix $D(m)$.
\begin{lemma}\label{wedding10}
There exists $c>0$ such that  $c\id \leq D(m) \leq c^{-1}\id$ for any
$m\in(0,1)$.
\end{lemma}
\begin{proof}
Given $a\in\RR^d$ we set $w:=\sum_{e\in\cE}a_e\psi _e $ and $v:=\sum
_{e\in \cE} a_e Pj_{0,e}$. Then (\ref{scozia}) and lemma
\ref{aritmetica}
  imply the upper bound
\begin{equation*}
(a,D(m) a) = \frac{1}{\c(m)}V_m(v,v)\leq  \frac{1}{\c(m)}V_m(\sum _{e\in
\cE} a_e j_{0,e} )\leq c\, \|a\|^2.
\end{equation*}
In order to prove the lower bound we observe that, by theorem  \ref{ortogonale}, \hbox{$V_m(v,w)=-\chi (m)\|a\|^2$} while,  thanks to (\ref{kolja}),
 $V_m(w)\leq c\, m(1-m) \|a\|^2$. Therefore, by  Schwarz inequality,
\begin{equation*}
(a,D(m)a)\geq \frac{1}{\c(m)}\frac{V_m(v,w)^2}{V_m(w)}\geq c\,\|a\|^2
\end{equation*}
thus proving the lemma.
\end{proof}

\begin{lemma}
\label{cemutvadie}
$D(m)$ is a continuous function on $(0,1)$.
\end{lemma}
\begin{proof}
  Let $0<\b$ and $0<\d<\ov 2$.  We observe that the limit point $\psi_e$
  of the sequence $\frac{\psi_{n,n}^e}{n}$ depends on the closure of
  $\cG/\cN$ and therefore on $m$. Therefore, it is convenient to denote
  it by $\psi_e^{(m)}$. Moreover, thanks to remark \ref{omogeneo} and
  lemma \ref{wedding10}, there exists $n_0\in\NN $ such that
\begin{equation}
\label{biology}
\|D\|_\infty
\sup_{m\in(0,1)} V_m\bigl(\psi_e^{(m)}- \frac{\psi_{n,n}^e}{n} \bigr)^{\ov 2}
\leq \b \quad \forall e\in \cE \quad \forall n\geq n_0
\end{equation}
where,\ $\|D\|_\infty :=\sup_{e,e'\in\cE}\|D_{e,e'}\|_\infty$.\\
Together with (\ref{decomposition}), this implies that, for any given
$m\in (0,1)$, we can find $g_m\in\cG$ such that
\begin{equation}
\label{genetics}
V_{m}\bigl(j_{0,e}+\sum_{e'\in\cE}D_{e,e'}(m)\frac{\psi_{n_0,n_0}^{e'}}{n_0}+
\cL g_m\bigr)^{\ov 2}\leq2\b.
\end{equation}
Since $\l_0(m)$ is a smooth function of $m\in (0,1)$ and thanks to Lemma
\ref{aritmetica}, (\ref{genetics}) remains valid if $V_m$ is replaced by
$V_{m'}$, where $m'$ is arbitrary inside an open interval $I_m$
containing $m$.  In what follows we restrict to the density interval
$[\d,1-\d]$.  Thanks to compactness and interpolation and thanks to
(\ref{kolja}), there exists a continuos matrix $D^{(\b)}(\cdot)$ and a
family of functions $g_m^{(\b)}$, $m\in[\d,1-\d]$, such that
$\|D_{e,e'}^{(\b)}\|_\infty\leq\|D_{e,e'}\|_\infty$ and
$$
V_{m}\bigl(j_{0,e}+\sum_{e'\in\cE}D^{(\b)}_{e,e'}(m)\frac{\psi_{n_0,n_0}^{e'}}{n_0}+
\cL g_m^{(\b)}\bigr)^{\ov 2}\leq 3\b\quad \forall m\in[\d,1-\d]
$$
and therefore
\begin{equation}
\label{nevski}
V_{m}\bigl(j_{0,e}+\sum_{e'\in\cE}D^{(\b)}_{e,e'}(m)\psi _{e'}^{(m)}+
\cL g_m^{(\b)}\bigr)^{\ov 2}\leq 4\b\quad \forall m\in[\d,1-\d]
\end{equation}
From the above formula and (\ref{decomposition}),  we have
$$ Pj_{0,e}=-\sum_{e'\in\cE}D_{e,e'}(m)\psi _{e'}^{(m)} = -\sum_{e'\in\cE}D^{(\b)}_{e,e'}(m)\psi _{e'}^{(m)}+ \xi ^{(m)}_e
\quad \forall m \in[\d,1-\d]
$$
where $ V_m(\xi^{(m)}_e)^{\ov 2}\leq4\b$.
By taking the scalar product  with $j_{0,e'}$ we obtain (thanks to theorem \ref{ortogonale})
$$|\chi(m)\bigl(D_{e,e'}(m)-D_{e,e'}^{(\b)}(m)\bigr) |\leq4  V_m(j_{0,e})^{\ov 2}\b\quad \forall m\in[\d,1-\d], $$
that is $|D_{e,e'}(m)-D_{e,e'}^{(\b)}(m)| \leq c(\d) \b$, 
thus proving that $D_{e,e'}(\cdot)$ is continuous on $[\d,1-\d]$.
\end{proof}

We are now able to prove our main result.
\begin{theorem}
 \label{agognatameta}
Let $d\geq 3$. Then given $\d>0$
\begin{equation}\label{esodo99}
\inf_{g\in\GG}\limsup _{n\uto\infty}
\sup_{m\in[\d,1-\d]}V_m\Bigl(j_{0,e}+\cL g+\sum_{e'\in\cE}D_{e,e'}(m) \frac{\psi_{n,n}^{e'}}{n}\Bigr)=0.
\end{equation}
Moreover, if $D$ has continuous extension to $\{0,1\}$, (\ref{esodo99}) is valid with $\d=0$.
\end{theorem}
\begin{proof}
(\ref{esodo99}) is a simple consequence of the estimates exhibited in the
proof of lemma \ref{cemutvadie}. Let us observe that, given $\b>0$, by 
defining  $g_m^{(\b)}$ as in the above  proof, then
\begin{equation}\label{trizio}
\limsup _{n\uto\infty}
\sup_{m\in[\d,1-\d]}V_m\Bigl(j_{0,e}+\cL g_m^{(\b)}+\sum_{e'\in\cE}D_{e,e'}(m) \frac{\psi_{n,n}^{e'}}{n}\Bigr)\leq c\,\b.
\end{equation} 
In order to define a function $g$ independent of $m$, it is enough to proceed as in the proof of corollary 5.9,  chapter 7, \cite{KL}.
If $D$ has continuous extension to $\{0,1\}$ then it is simple to extend
 (\ref{trizio})  to all $[0,1]$. 
\end{proof}

\appendix
\section{}
In this final appendix we have collected several technical results
 used in the previous sections.

\subsection{Large deviations estimates}.

\begin{lemma}\label{upupa}
 Let $f=f(\a)$ be a mean-zero local function and $\La\in\FF$ be such that
\hbox{$(\D_f+x)\cap(\D_f+y)=\emptyset $} for any $x,y\in\La$.
 Then
\begin{equation*} \label{iride}
\PP [\,|\Av_{x\in\L}\t_xf\,|\geq\d\,]\leq
2e^{-\frac{\d^2|\La|}{4\|f\|_\infty^2}}\qquad \forall \d>0.
\end{equation*}
\end{lemma}
\begin{proof}
Given $t>0$,
since $\EE(f)=0$ and  $e^x-x\leq e^{x^2}$ for any $x\geq 0$,
$$
e^{tf}=\sum_{n=0}^\infty \frac{(tf)^n}{n!}\leq
e^{t\|f\|_\infty}-\|f\|_\infty t \leq e^{\|f\|_\infty^2}.
$$
Therefore, thanks  to the conditions on $f$ and $\L$,
\begin{equation*}
\PP [\,\Av_{x\in\L}\t_xf \geq\d\,]\leq  e^{-t\d}
\EE(e^{t\,\Av_{x\in\L}\t_xf })= e^{-t\d}[\EE(e^{tf|
\L|^{-1} })]^{|\L|}\leq e^{-t\d+t^2 \|f\|_\infty^2|\L|^{-1}}.
\end{equation*}
The thesis follows by taking  $t:=\d|\L|/(2\|f\|_\infty^2)$ and
by considering the above estimates  with $f$ replaced by $-f$.
\end{proof}

\subsection{Equilibrium bounds}
\begin{lemma}
 \label{turu}
Given  $\L\in\FF$ and $\l\in\RR$ we define $m:=\mu^\l(m_\L)$ and
 $a_m:=\min(m,1-m)$.
Then, for any $\D\sset\L$ and any function $f$ such that $\D_f\sset\L$,
\begin{align*}
& a) \; c|\D|m \leq\mu(N_\D)\leq c^{-1}|\D|m,\\
& b) \; c\,|\D|(1-m)\leq \mu(|\D|-N_\D)\leq c^{-1}|\D|(1-m), \\
& c) \; c |\D|a_m\leq \mu(N_\D;N_\D)\leq c^{-1}|\D|a_m,\\
& d) \; |\mu(f;N_\L)|\leq
 c\|f\|_\infty\min\Bigl(|\D_f|a_m,\sqrt{|\D_f|a_m}\,
\Bigr).
\end{align*}
\end{lemma}
\begin{proof}
In what follows we assume  $m\leq \ov2$.\\
$a)$ and $b)$ can be easily derived from the boundedness of the random
field $\a$. Let us  prove $c)$. The upper bound
follows by observing that $\mu(N_\D;N_\D)\leq \mu (\D)$ and by applying
$a)$.
 In  order to prove the lower bound, let us
introduce the set \hbox{$W:=\{x\in\L:\mu(\h_x)\leq\ov2\}$}. Since
 $|W|\geq|\L|/2$ and thanks to  $a)$,
\begin{equation*}
 \mu(N_\L;N_\L)\geq  \mu(N_W;N_W)\geq\ov2\mu(N_W)\geq c\,m|\L|
\end{equation*}
thus proving the lower bound in $c)$ with $\D$ replaced by $\L$. In
order to consider the general case, we define $m'=\mu(m_\D)$. Then by the
previous arguments, $\mu(N_\D;N_\D)\geq c\,m'|\D|$ which, by $a)$, is
bounded from below by  $c\,m|\D|$.\\
Let us prove $d)$. By Schwarz inequality and $c)$
\begin{equation*}
|\mu(f;N_\L)|\leq\mu(f;f)^{\ov2}\mu(N_{\D_f};N_{\D_f})^{\frac{1}{2}}
\leq c \,\mu(f;f)^{\ov2}\sqrt{m\,|\D_f|}
\end{equation*}
Since $\mu(f;f)\leq \|f\|_\infty^2$,
it remains to prove that $\mu(f;f)\leq c\,m\|f\|_\infty^2|\D_f|$.
To this aim let $\h^\ast$ be the configuration with no particle. Then,
thanks to $a)$,
\begin{equation*}
\mu(f;f)\leq\mu\bigl(\bigl(f-f(\h^\ast)\bigr)^2\bigr)
\leq\|f\|_\infty^2\mu(N_{\D_f})\leq c\,\|f\|_\infty^2|\D_f|.
\end{equation*}
\end{proof}

\begin{lemma}
 \label{vespe}
For any  $\l,\l'\in\RR$,  $\L\in\FF$ and any function
$f$  with $\D_f \sset \L$,
\begin{align}
& |\mu^{\l'}(f)-\mu^\l(f)|\leq c
\|f\|_{\infty}|\D_f|\,|\mu^{\l'}(m_\L)-\mu^{\l}(m_\L)|,\label{sam1}\\
& |\mu^{\l'}(\h_x;\h_x)-\mu^\l(\h_x;\h_x)|\leq
c\,|\mu^{\l'}(m_\L)-\mu^{\l}(m_\L)| \quad \forall x\in \L\label{sam2}.
\end{align}
For any  $m,m'\in(0,1)$ and any local function $f$,
\begin{align}
& |\mu^{\l_0(m')}(\h_0)-\mu^{\l_0(m)}(\h_0)|\leq c\,|m'-m|,
 \label{sam3}\\
& |\mu^{\l_0(m')}(\h_0;\h_0)-\mu^{\l_0(m)}(\h_0;\h_0)|\leq c\,|m'-m|,
  \label{sam4}\\
& |\mu^{\l_0(m')}(f)-\mu^{\l_0(m)}(f)|\leq c(|\D_f|)\,\|f\|_\infty |m'-m|
\label{zigeuner}
\end{align}
for a suitable constant $c(|\D_f|)$ depending on $|\D_f|$.\\
Moreover, for any $\L\in\FF$ and  any $m\in(0,1)$,
\begin{equation}
  \label{sam5}
|\la_\L(m)-\l_0(m)|\leq  \frac{c}{m(1-m)}|m-\mu^{\l_0(m)}(m_\L)|.
\end{equation}
\end{lemma}

\begin{proof}
It is simple to derive (\ref{sam2}), (\ref{sam4}) and (\ref{zigeuner})  from
(\ref{sam1}) and (\ref{sam3}).\\
Let us prove (\ref{sam1}). By setting  $\l(s):=\l_\L(s)$,
$m:=\mu^\l(m_\L)$ and $m':=\mu^{\l'}(m_\L)$, we have
\begin{equation*}
  \mu^{\l'}(f)-\mu^{\l}(f)=
 \int _m ^{m'}\frac{\partial}{\partial s}\mu^{\l(s)}(f)\,ds=
\int_m ^{m'}\mu^{\l(s)}(f;N_{\D_f})\l'(s)ds.
\end{equation*}
By lemma \ref{turu},
$|\mu^{\l(s)}(f;N_{\D_f})\l'(s)|\leq c\|f\|_\infty |\D_f|$, thus
concluding the proof of (\ref{sam1}).\\
In order to prove (\ref{sam3}), we observe that
\begin{equation*}
\mu ^{\l_0(m')}(\h_y)-\mu^{\l_0(m)}(\h_y)=\int_{m}^{m'}\frac{d}{ds}
\mu^{\l_0(s)}(\h_0)ds =\int_{m}^{m'}
\frac{\mu^{\l_0(s)}(\h_0;\h_0)}{\EE\, \mu^{\l_0(s)}(\h_0;\h_0)}ds.
\end{equation*}
Thanks to the boundedness of the random field $\a$, the last integrand
 is bounded, thus proving (\ref{sam3}).\\
Let us prove (\ref{sam5}). By Lagrange theorem
\begin{equation*}
m=\mu_\La^{\l(m)}(m_\L)=
\mu ^{\l_0(m)}(m_\L)+\mu^\l(m_\L;N_\L)\bigl(\l_\L(m)-\l_0(m)\bigl)
\end{equation*}
for a suitable $\l$ between $\l_\L(m)$ and $\l_0(m)$.
In order to conclude the proof, it is enough to apply lemma \ref{turu}.
\end{proof}

\subsection{Equivalence of ensembles}
In this paragraph we compare   multi-canonical and multi-grand canonical
expectations. The following results can be proved by the same methods
 developed in \cite{CM1} with strong simplifications since here the
 grand canonical measures are product (see \cite{AF} for a
complete treatment).

In what follows  we fix $\D\in\FF$ and we partition it as
$\D=\cup_{i=1}^k\D_i$. Moreover, chosen  a set
$\mathbf{N}=\{N_i\}_{i=1}^k$   of possible particle numbers in each atom
$\D _i $, we define the multi-grand canonical measure $\bar\mu$
and the multi-canonical measure $\bar\nu$ as
\begin{align*}
&\bar\mu:= \otimes_{i=1}^k
\mu_{\D_i}^{\la(m_i)}\quad \text{where}\quad m_i:=\frac{N_i}{|\D_i|},\\
&\bar\nu:=\mu(\cdot\tc N_{\D_i}=N_i\;\; \forall i=1,\dots,k ).
\end{align*}
Then we have the following main results (for the latter see also
proposition $3.3$ in \cite{CM2}).

\begin{prop}\label{equi} (Equivalence of ensembles)\\
Let  $\g,\d\in(0,1)$  and $f$ be a local function such that
 $|\D_i|\geq \d|\D|$, for any
$i=1,\dots,k$,  $\D_f\sset\D$  and  $|\D_f|\leq|\D|^{1-\g}$.\\
Then there exist constants $c_1,c_2$, depending respectively  on
$\g,\d,k$ and $\d,k$, such that
\begin{equation*}
|\D|\geq c_1\quad \Rightarrow \quad
 |\bar\nu(f)-\bar\mu(f)|\leq c_2\,\|f\|_{\infty}\frac{|\D_f|}{|\D|}.
\end{equation*}
\end{prop}
\begin{lemma}\label{cleonimo}
Let $\d\in(0,1)$ and $f$ be a local function such that $\D_f\sset\D$ and
 \hbox{$|\D_i\setminus\D_f|\geq\d|\D_i|$} for any $i=1,\cdots,k$.\\
Then there exist  constants $c_1,c_2$, depending respectively  on $k$ and
$k,\d$, such that
\begin{equation*}
|\D_i|\geq c_1 \quad \forall i =1,\dots,k\quad  \Rightarrow \quad \bar\nu(|f|)\leq c_2\, \bar\mu(|f|)\quad \text{and} \quad
\Var_{\bar\nu}(f)\leq c_2\Var_{\bar\mu}(f),
\end{equation*}
\end{lemma}

\subsection{Some special equilibrium covariances}
In this paragraph we estimate the canonical covariance between a
generic function and a function which can be written as the spatial
average of local functions. We observe that the bound we provide differs
from the standard Lu-Yau's Two Blocks Estimate (see \cite{LY}) by an additional
term depending on the random field $\a$ and satisfying a large deviations
estimate.

\\
In what follows we fix  functions $h,h'\in\GG$,
 depending only on  $\a_0$ and $\h_0$,   such that
 $\|h\|_\infty,\|h'\|_\infty \leq 1$. Moreover, for any positive integer $L$,
 we denote by $R_L$ the set of boxes with sides of length
in $[L, 100L]$.

\begin{prop}
 \label{special_covariances}
Given $0<\d<\ov2$ there exists $\ell _0\in\NN$ having the following
property.\\
 Let $\ell , L\in\bbN$ be such that $\ell_0\leq \ell\leq L$ and
let $V,W\in R_L$ with $V\cap W=\emptyset$.  Then, for any $\nu\in\cM(V)$
and any function $g\in\GG$,
\begin{equation}
  \label{eq:special_covariances1}
\nu(\Av_{v\in V}\t_v h;g)^2\leq\frac{c(\ell)}{|V|}\cD(g;\nu)+\frac{c\d}{|V|}
\Var_\nu(g)+\frac{c}{|V|}  \Var_\nu(g)\id_{\{m\in I_\d\}}\th_{V,\ell}(\a)
\end{equation}
where $m:=\nu(m_V)$ and  $I_\d:=[\d,1-\d]$. Moreover,  for any $\g>0$
 there exists $\ell_1 =\ell_1(\g)\geq \ell _0  $  such that
\begin{equation}
 \label{eq:special_covariances2}
\ell _1 \leq \ell\leq L\quad \Rightarrow  \quad \PP(\,\th_{V,\ell}(\a)\geq\g\,)\leq   e^{-c(\g,\ell) L^d}.
\end{equation}
Finally, for any $\nu\in\cM(V\cup W)$ and any function $g\in\GG$,
\begin{equation}
 \label{eq:special_covariances3}
\begin{split}
\nu(\Av_{v\in V}\t_v h\cdot &\Av_{w\in W}\t_w h';g)^2\leq\\
& \frac{c(\ell)}{|\La|}\cD
(g;\nu)+\frac{c\d}{|\La|}\Var_\nu(g)+\frac{c}{|\La|}
\Var_\nu(g)\bigl( \th_{V,\ell}(\a)+\th_{W,\ell}(\a)\bigr)
\end{split}
\end{equation}
\end{prop}

\begin{proof}
We first prove (\ref{eq:special_covariances1}) by referring, for many
steps, to the proof  of  proposition A.1 in \cite{CM2}. Let us fist
introduce some useful notation.\\
We fix a  partition  $V=\cup_{i\in I}Q_i$, with $Q_i\in R_\ell$, and
define $N_i:=N_{Q_i}$,  $m_i:=N_{Q_i}/|Q_i|$, $h_i:=\sum_{x\in Q_i}\t_x h$,
 $\cF:=\s (m_i\tc i\in I)$ and for  $s\in[0,1]$
\begin{equation*}
A_i(m):=\frac{\mu_V^{\l(m)}(h_i;N_i)}{\mu_V^{\l(m)}
(N_i;N_i)}-\frac{\EE\,\mu^{\l_0(m)}(h_0;\h_0)}{\EE\,\mu^{\l_0(m)}(\h_0;\h_0)}
\qquad B_i(s):=
\frac{\mu_{Q_i}^{\la(s)}(h_i;N_i)}{\mu_{Q_i}^{\la(s)}(N_i;N_i)}-\frac{\mu_V^{\l(m)}(h_i;N_i)}{\mu_V^{\l(m)}(N_i;N_i)}.
\end{equation*}
As in \cite{CM2},  if $m\not\in I_\d$ then  it is enough to apply
Schwarz inequality and lemma \ref{cleonimo} to obtain the thesis,
 otherwise it is convenient
to bound the l.h.s. of (\ref{eq:special_covariances1}) as
\begin{equation}
  \label{kafka}
 \nu(\Av_{v\in V}\t_v h;g)^2\leq 2\nu\bigl(\,\nu(\Av_{v\in V}\t_v h
 ;g\tc \cF)
\,\bigr)^2+2\nu\bigl(\,\nu(\Av_{v\in V}\t_v h \tc \cF);g\,\bigr)^2
\end{equation}
As in  \cite{CM2} we can bound the first addendum in the r.h.s. of
 (\ref{kafka})  by $c(\ell)\cD(g;\nu)$ and   the second one by
\begin{equation}
  \label{killer}
c\Var _\nu (g) \Bigl(\frac{1}{\ell ^d L^\d}+\frac{1}{L^d}
\sum _{i\in I}\Var_{\mu_V^{\l(m)}}(\xi^\g_i)\Bigr )
\end{equation}
where, for an arbitrarily fixed $\g$, $\xi^\g_i(\h):=\mu_{Q_i}
^{\l(m_i(\h))}(h_i-\g N_i)$.
 Let us explain how to proceed. Thanks to
 Poincar\'e  inequality for Glauber dynamics we obtain
\begin{equation}
 \label{north}
 \Var_{\mu_V^{\l(m)} }(\xi^\g_i)\leq c\,\sum_{x\in Q_i}\mu_V^{\l(m)}
 ((\nabla_x\xi_i^\g)^2).
\end{equation}
By choosing
 $\g=\frac{\EE\,\mu^{\l_0(m)}(h_0;\h_0)}{\EE\,\mu^{\l_0(m)}(\h_0;\h_0)}$
it is simple to check that
\begin{equation*}
\nabla_x\xi_i^\g=A_i+(-1)^{\h_x}\int_{m_i(\h)}^{m_i(\h^x)} B_i (s)ds .
\end{equation*}
By writing
\begin{multline*}
B_i(s)=  \frac{|Q_i|}{\mu_V^{\l_(m)} (N_i;N_i)}\Bigl(
\int _{\mu_V^{\l(m)}(m_i)}^{s}
\frac{\mu_{Q_i}^{\l(s')}(h_i;N_i;N_i)}{\mu_{Q_i}^{\l(s')}(N_i;N_i)}ds'
+\\
 \frac{\mu_{Q_i}^{\l(s)}(h_i;N_i)}{\mu_{Q_i}^{\l(s)}(N_i;N_i)}
\int _{\mu_V^{\l(m)}(m_i)}^{s}
\frac{\mu_{Q_i}^{\l(s')}(N_i;N_i;N_i)}{\mu_{Q_i}^{\l(s')}(N_i;N_i)}ds'\Bigr),
\end{multline*}
by lemma \ref{turu} and the condition $m\in I_\d$ we obtain that
$|B_i(s)|\leq \frac{c}{\d}|s-\mu_V^{\l(m)}(m_i)|$ and therefore
\begin{equation}
 \label{south}
  |\nabla_x\xi_i^\g|\leq A_i+\frac{c}{\d}
 \bigl |\,m_i(\h)- \mu_V^{\l(m)}(m_i)\bigr |
+\frac{c}{\d}\ell ^{-d}.
\end{equation}
By (\ref{killer}), (\ref{north}) and (\ref{south}) it is simple to
conclude the proof
if $\ell$ is large enough and
\begin{equation*}
 \th_{V,\ell}(\a):=\sup _{m\in M_V}\Av _{i\in I}A _i(m)  ^2 \quad \text{where}
\quad M_V=\bigl\{\frac{1}{|V|}, \frac{2}{|V|}, \dots,
1-\frac{1}{|V|}\bigr\}.
\end{equation*}

\smallno
By standard arguments (as for lemma $3.9$ in \cite{CM2})
 (\ref{eq:special_covariances3}) can be derived from
(\ref{eq:special_covariances1}).\\
 Let us prove (\ref{eq:special_covariances2}).
By lemmas \ref{turu} and \ref{vespe}
\begin{equation*}
\begin{split}
|A_i|\leq\frac{c}{m(1-m)}\Bigl (\, & |\mu^{\l_0(m)}(h_i;m_i)-
 \EE\,\mu^{\l_0(m)}(h_0;\h_0)|+\\
 &  |(\mu^{\l_0(m)}(N_i;m_i)- \EE\,\mu^{\l_0(m)}(\h_0;\h_0)|
  + | m-\mu^{\l_0(m)}(m_V)|\,\Bigr).
\end{split}
\end{equation*}
Therefore it is enough to prove that given a function $f=f(\a_0)$ with
 $\|f\|_\infty \leq 1$ then  for any $\g>0$ there exists
$\ell_1=\ell_1(\g)$ such that
\begin{equation*}
 \PP \bigl(\Av _{i \in I }(\Av _{x\in Q_i}\t_x f-\EE(f))^2\geq\g\,\bigr)\leq 2
e ^{- \frac{c\g^2 L^d}{\ell^d}} \qquad \forall \ell \geq \ell_1 .
\end{equation*}
To this aim  we define  $f_i:=(\Av _{x\in Q_i}\t_x f-\EE(f))^2$ and $\bar
f_i:=f_i-\EE(f_i)$. Then by lemma $\ref{upupa} $, for
any $0<\d<1 $,
$$
\EE(f_i)\leq\PP\bigl(\,|\Av _{x\in Q_i}\t_x f-
\EE (f)|\geq\d\,\bigr)+\d^2\leq
2 e ^{-c\,\d^2\ell^d}+\d^2.
$$
Therefore, by choosing $\d $  small enough and $\ell$ large enough,
 $\EE(f_i)\leq\frac{\g}{2}$ for any $i\in I$ and
(by applying again lemma \ref{upupa})

\begin{equation}
\PP ( \,\Av_{i\in I}f_i\geq\g\,)\leq
\PP (\,\Av_{i\in I}\bar f_i\geq\frac{\g}{2}\,) \leq 2 e^{-
c\,\g^2|I|}
\end{equation}
thus concluding the proof.
\end{proof}

\subsection{Moving Particle Lemma}
Given $x,y\in\ZZ^d$ we define
$$
z_i:=(y_1,y_2,\dots,y_i,x_{i+1},\dots,x_d)\qquad \forall i= 0,\dots,d
$$
and write $\g_{x,y}$ for the path connecting $z_0=x$ to $z_1$ by moving
along the first direction, then connecting  $z_1$ to $z_2$ by
moving along the second direction and so on until arriving to $z_d =y$.
We denote by $ | \gamma _{x,y}|$ the length of  the path $\gamma
_{x,y}$.

\begin{lemma}\label{MPL} (Moving Particles Lemma)\\
Given  a box $\L$ and $\nu\in\cM(\L)$  then
$$
\nu\bigl(\,(\nabla_{x,y}f)^2\,\bigr)\leq
c\, |\g_{x,y}|\sum_{b\in\g_{x,y}}\nu\bigl(\,(\nabla_b f)^2\,\bigr)\qquad
\forall x,y\in\L,\; f\in\GG.
$$
\end{lemma}
The above lemma is well known for non disordered systems (see
for example \cite{SpohnYau}).  We learned from J.Quastel  the
generalization to the disordered  case.

\subsection{An application of Feynman-Kac formula}
The following proposition can be derived from the
Feynman-Kac formula as explained in
\cite{KL}. We report only the statement.

\smallno
Let $X$ be a finite set on which it is
defined a probability measure $\nu$ and a Markov generator $\frL$ reversible
w.r.t. $\nu$. We denote by $\EE_\nu$  the expectation w.r.t. the
Markov process having infinitesimal generator $\frL$ and initial
distribution $\nu$ and by $x_t$ the configuration at time $t$.
\begin{prop}
\label{fk}
Let $V:\RR_{+}\times X\rightarrow\RR$ be a bounded measurable function
and let, for any $t>0$,
$$
\G_t:=\sup spec_{L^2(\nu)}\{\frL+V(t,\cdot)\}.
$$
Then
\begin{equation*}
\EE _{\nu}\Bigl[\,\exp\Bigl\{\int_0^tV(s,x_s)ds\Bigr\}\Bigr]\leq
\exp\Bigl\{\int_0^t\G_s\, ds\Bigr\} \quad \forall t>0.
\end{equation*}
\end{prop}

\subsection{Two Blocks Estimate}\label{rachmaninov} 

For a treatment of the Two Blocks estimate in non disordered systems see
\cite{KL} and reference therein. Let us state and prove a generalized
version.
\begin{prop}\label{tvb}
Given $\gamma >0$, for almost any disorder configuration
$\alpha$
\begin{equation}\label{tvbeq}
\limsup_{a\dto 0,k\uto\infty,\ep\dto 0}
\sup_{w:|w|\leq\frac{a}{\ep}}\sup spec _{L^2(\mu_\ep)}
\{\Av_{x\in\toromi}|m_{x,k}-m_{x+w,k}|+\g\ep^{d-2}\cL_\ep\}\leq 0.
\end{equation}
\end{prop}
\begin{proof}
  We extend to the disordered case the proof of the Two Blocks estimate
  of \cite{VY} thanks to the ergodicity of the random field $\a$.  To
  this aim let us introduce the scale parameter $\ell$ with
  $\ell\uto\infty$ after $k\uto\infty$.  Then, with a negligible error
  of order $O(\ell/k)$, for any $x\in \toromi$ we can substitute
  $m_{x,k}$ with $\Av_{y\in\La_k}m_{x+y,\ell}$. Therefore, thanks to the
  subadditivy of $\sup spec$, the l.h.s. of (\ref{tvbeq}) can be bounded
  from above (with an error $O(\d)$) by
\begin{equation*}
\sup_{w:|w|\leq\frac{a}{\ep}}\Av _{y \in\La_k}\Av _{y'\in\La_k\,:\,|w+y'-y|>
2\ell}\;\underset{L^2(\mu_\ep)}{\sup spec}
\{\Av_{x\in\toromi}|m_{x+y,\ell} -m_{x+w+y',\ell}|+\g\ep^{d-2}\cL_\ep\}
\end{equation*}
where the additional restriction $|w+y'-y|>2\ell$ is painless. By renaming the
index variables, it is enough to show that given $\g>0$, for
almost any disorder configuration $\a$,
\begin{equation}\label{fata}
 \limsup_{\ell\uto\infty ,a\dto 0,\ep\dto 0}\sup_{w:2\ell<|w|\leq\frac{a}{\ep}}\underset{L^2(\mu_\ep) }{\sup spec}\,\{\Av_{x\in\toromi}|m_{x,\ell}-m_{x+w,\ell}|+\g\ep^{d-2}\cL_\ep\}\leq 0.
\end{equation}
For any $u,v\in\ZZ^d$ let us define $\hat\cL_{u,v}=(1+
e^{-(\a_u-\a_v)(\h_u-\h_v)})\nabla_{u,v}$. It is simple to check that
$\hat\cL_{u,v}$ is self-adjoint w.r.t. Gibbs measures. Then, given $w$
as above, thanks to the Moving Particle lemma (see lemma \ref{MPL}) and
the properties of the transition rates, it is simple to prove that
\begin{equation}\label{di}
\Av _{x\in\toromi}\Av_{u\in\L_{x,\ell}}\Av_{v\in\L_{x+w,\ell}}(-\hat\cL_{u,v})\leq c\, a^2 \ep^{d-2}(-\cL_\ep).
\end{equation}
Therefore, by localizing as in (\ref{cappy}), the $sup spec$ in
(\ref{fata}) is bounded by
$$
\Av _{x \in \toromi }\sup_{\nu}\,\sup spec_{L^2
(\nu)}\{|m_{x,\ell}-m_{x+w,\ell}|+c\,\g a^{-2}\,
\Av_{u\in\L_{x,\ell}}\Av_{v\in\L_{x+w,\ell}}(-\hat\cL_{u,v})\}
$$
where $\nu$ varies in $\cM(\L_{x,\ell}\cup\L_{x+w,\ell})$. Thanks to
perturbation theory (see proposition \ref{jamme}) we only need to prove
that, for almost any disorder configuration $\a$,
\begin{equation}\label{pepe}
\limsup_{\ell\uto\infty,a\dto 0, \ep\dto 0}\sup_{w:2\ell<|w|\leq\frac{a}{\ep}}
 \Av_{x\in\toromi} \sup _{\nu}\, \nu (|m_{x,\ell}-m_{x+w,\ell}|).
\end{equation}
We observe that by lemma \ref{cleonimo} in the above expression we can
\item substitute $\nu$ with  the grand canonical measure $\mu$ 
such that $\mu(m_\L)=\nu(m_\L)$ where $\L:=\L_{x,\ell}\cup\L_{x+w,\ell}$.\\
  Let us introduce the scale parameter $s$ with $s\uto\infty$ after
  $\ell\uto \infty$. Then, by approximating $m_{x,\ell}$ with
  $\Av_{y\in\L_{x,\ell}}m_{y,s}$ and thanks to lemma \ref{vespe}
\begin{equation*}
  \mu (|m_{x,\ell}-m_{x+w,\ell}|)\leq 
c\,\Av_{y\in\L_{x,\ell}}\mu^{\la_0 (m)}(|m_{y,s}-m_{y+w,s}|)
+ c\,s^d\bigl |m-\mu^{\la_0(m)}(m_\L)\bigl|+O(s/\ell).
\end{equation*}
where $m=\mu(m_\L)=\nu(m_\L)$ and $\L $ is defined as above.  Therefore,
it is enough to prove that for almost all disorder configuration $\a$
\begin{equation*}
\begin{split}
  & \limsup _{s\uto\infty,\ell\uto\infty ,\ep\dto
    0}\Av_{x\in\toromi}\sup_{m}\Av_{y\in\L_\ell}\mu^{\la_0 (m)}(|m _{x+y,s}-m|)=0, \\
  & \limsup _{\ell\uto\infty ,\ep\dto 0}\Av_{x\in
    \toromi}\sup_{m}|m-\mu^{\l_0 (m)}(m_{x,\ell})|=0.
\end{split}
\end{equation*}
Since $\EE\mu^{\l_0 (m)}(m_{x,n})=m$ for any integer $n$ and any site
$x$, the above limits follow by straightforward arguments from the
ergodicity of the random field $\a$ and the technical estimate
(\ref{sam3}).
\end{proof}



\end{document}